\documentclass[a4paper,11pt]{article}

\usepackage{amssymb,amsmath,amsthm}
\usepackage{mathrsfs}
\usepackage{authblk}
\usepackage{indentfirst}
\usepackage[numbers,sort&compress]{natbib}
\usepackage[english]{babel}
\usepackage{caption}
\usepackage{hyperref}
\usepackage{bookmark}
\usepackage{tikz}
\usepackage{enumitem}
\usepackage{tikz-cd}
\usepackage[section]{placeins}
\usetikzlibrary{arrows.meta,positioning}

\newtheorem{theorem}{Theorem}[section]
\newtheorem{corollary}[theorem]{Corollary}

\newtheorem{lemma}[theorem]{Lemma}
\newtheorem{proposition}[theorem]{Proposition}

\theoremstyle{definition}
\newtheorem{definition}[theorem]{Definition}
\newtheorem{question}[theorem]{Question}
\newtheorem{example}[theorem]{Example}
\newtheorem{remark}[theorem]{Remark}

\DeclareMathOperator{\im}{im}
\DeclareMathOperator{\rank}{rank}

\DeclareMathOperator{\Lk}{Lk}
\DeclareMathOperator{\st}{St}
\DeclareMathOperator{\spann}{span}

\newcommand{\LH}{\mathrm{LH}}
\makeatletter
\renewcommand*\l@subsection{\@dottedtocline{2}{1.5em}{2.8em}}
\makeatother

\title{On the local homology of path complexes}

\author[1]{Jian Liu}
\author[1]{Kefeng Liu}

\affil[1]{Mathematical Science Research Center, Chongqing University of Technology, Chongqing 400054, China}

\makeatletter
    \renewcommand*{\@fnsymbol}[1]{\ensuremath{\ifcase#1\or \dagger\or *\or *\or
   \mathsection\or \else\@ctrerr\fi}}
\makeatother

\date{}

\begin{document}
\maketitle

\paragraph{Abstract}

The path homology introduced by Grigor'yan, Lin, Muranov and Yau is a relatively new theory that gives a topological characterization of digraphs and quivers. In this paper, we define local homology for path complexes over finite vertex sets. We construct star and link path complexes that generalize the classical closed stars and links from simplicial topology. The standard link formula no longer holds: local homology in degree $r$ can differ from the corresponding reduced link homology in degree $r-1$. We introduce the defect complex and the completed link chain complex to analyze this obstruction. We give criteria for recovering local homology from the link path complex through the associated exact sequences. Given a vertex $v$ and an integer $s\ge 1$, we show that the $s$-neighborhood of $v$ determines its local homology in all degrees $q<s$. This method further yields additional results, including stability of local homology after removing distant vertices and sufficient conditions for the inclusion of graphs to induce a path homology isomorphism. We define an interaction complex that measures the failure of local homology to be additive over disjoint marked sets. Applying this construction to singletons, we show that if each inclusion $G-v\hookrightarrow G$ for $v\in T$ induces a path homology isomorphism in each degree, then $G-T\hookrightarrow G$ does so if and only if the interaction complex is acyclic. Some other results include Mayer--Vietoris sequences, excision criteria, a local K\"unneth formula for rectangular marked sets, gluing formulas along matchings, and cubical link formulas for lattice digraphs.


\tableofcontents

\section{Introduction}\label{sec:introduction}

The path homology theory of Grigor'yan, Lin, Muranov, and Yau provides homological invariants of digraphs and path complexes \cite{GLMY,GLMYCohomology,GLMYOverview}. A path complex specifies allowed sequences of vertices and includes both the directed paths of a digraph and the increasing vertex sequences representing simplices. Its chain groups consist of linear combinations of allowed paths whose boundaries are again allowed. The theory admits homotopy invariance for digraphs \cite{GLMYHomotopy} and K\"unneth formulas for Cartesian products and joins \cite{GMYKunneth}. For regular path complexes, the Cartesian product is described by the shuffle construction of \cite{GLMYProduct}.

Relative homology and gluing results provide tools for defining and computing local homology. In \cite{GJMY,CDKOSW}, relative path homology is studied alongside the Eilenberg–Steenrod axioms and the cofibration category of digraphs. Gluing exact sequences, spectral sequences for path homology, and Mayer–Vietoris theorems within the GLMY framework appear in \cite{GrigorLiu,Asao,HepworthRoff}. Other work on path homology compares it with simplicial and cubical homology \cite{GMY,BGJW} and examines obstructions to topological realization \cite{FuIvanov}. Further developments include persistent path homology \cite{ChowdhuryMemoli}, discrete Morse theory for path homology \cite{LinWangYau}, and product formulas related to persistent homology \cite{GakharPerea}.

Classical simplicial homology is the main motivation for our work. For a simplicial complex $K$ and a vertex $v$, we have the standard isomorphism \cite{Hatcher}
\[
 H_r(K,K-v;\Bbbk)
 \cong H_r(\overline{\st}_K(v),\Lk_K(v);\Bbbk)
 \cong\widetilde H_{r-1}(\Lk_K(v);\Bbbk).
\]
The closed star contains all simplices related to $v$, and the link describes how these simplices connect around $v$. Persistent local Laplacians give a related description through harmonic spaces and spectra \cite{LiuFengLiu}. For path complexes, the difficulty lies in identifying the appropriate notions of star and link path complexes and understanding how the link path complex can be used to characterize local path homology.

Let $P$ be a regular path complex on a nonempty finite vertex set, and let $\Omega_*(P)$ be the $\partial$-invariant chain complex. For a marked set $T$, write $P-T$ for the induced path complex on $V(P)\setminus T$ and set $Q_*(P,T):=\Omega_*(P)/\Omega_*(P-T)$. The \textit{relative path homology} is defined by
\[
 \LH_r(P,T):=H_r(Q_*(P,T)),\quad r\geq 0.
\]
For a vertex $v$ of $P$, the \textit{local path homology} at $v$ is $\LH_r(P,v)=\LH_r(P,\{v\})$. 

Now, we write $\beta\preceq\alpha$ when $\beta$ is an ordered subsequence of $\alpha$. Define the \textit{star path complex} and \textit{link path complex} at $v$ by
\[
 \begin{aligned}
  S_v(P)&:=\{\beta\in P:\beta\preceq\alpha
  \text{ for some }\alpha\in P\text{ containing }v\},\\
  L_v(P)&:=S_v(P)-v.
 \end{aligned}
\]
These are path subcomplexes of $P$ with $P=S_v(P)\cup(P-v)$ and $L_v(P)=S_v(P)\cap(P-v)$. For the path complex from a simplicial complex, they are the usual closed star and link. For a digraph $G$, let
\[
 R_v=\{x\in V(G):x\leadsto v\text{ or }v\leadsto x\},
 \quad S_v=G[R_v],\quad L_v=S_v-v,
\]
where $\leadsto$ denotes directed reachability, allowing paths of length zero. Then we have
\[
 S_v(P(G))=P(S_v),\quad L_v(P(G))=P(L_v).
\]

The simplicial link formula need not hold for a general path complex. In the diamond with edges $a\to v$, $v\to b$, $a\to w$, and $w\to b$, one has (Example~\ref{ex:diamond})
\[
 H_1(S_v,L_v)\cong\Bbbk,\quad \LH_1(G,v)=0.
\]
The reduced homology of the star may also be nonzero, even for a directed acyclic graph (Example~\ref{ex:nonacyclic-star}). These two obstructions are described by the comparison maps
\[
 \LH_r(P,v)\xleftarrow{\ j_r\ }
 H_r(S_v(P),L_v(P))
 \xrightarrow{\ \delta_r\ }\widetilde H_{r-1}(L_v(P)).
\]

The \textit{star defect complex} at $v$
\[
 D_*(P,v):=
 \frac{\Omega_*(P)}{\Omega_*(S_v(P))+\Omega_*(P-v)}
\]
measures the failure of invariant chains to decompose over the cover $P=S_v(P)\cup(P-v)$. Put $U=V(S_v(P))$. If $S_v(P)=P|_U$ and $U$ is convex with respect to the allowed paths of $P$, then $D_*(P,v)=0$, and $j_r$ is an isomorphism for each $r\ge0$ (Proposition~\ref{prop:convex}). Here, convexity means that every allowed path with endpoints in $U$ lies entirely in $U$.

The defect also vanishes when $P$ is perfect, so each $j_r$ is an isomorphism (Corollary~\ref{cor:perfect-star-local}). If $P$ is both perfect and monotone, its star corresponds to a simplicial cone, and we recover
\[
 \LH_r(P,v)\cong\widetilde H_{r-1}(L_v(P)),
 \quad r\ge0.
\]

For a path complex $P$, we construct a \textit{completed link chain complex} $\widehat{\mathcal L}_*(P,v)$ and a long exact sequence
\[
 \cdots\longrightarrow\widetilde H_r(S_v(P))
 \longrightarrow\LH_r(P,v)
 \longrightarrow\widetilde H_{r-1}(\widehat{\mathcal L}_*(P,v))
 \longrightarrow\widetilde H_{r-1}(S_v(P))\longrightarrow\cdots.
\]
If $\widetilde H_r(S_v(P))=\widetilde H_{r-1}(S_v(P))=0$, the middle map is an isomorphism. For a simplicial complex, the completed link is the chain complex of the usual link. We also obtain a Mayer--Vietoris sequence for covers by induced path subcomplexes whenever each allowed path lies in one member of the cover.

Local homology can be determined by a bounded neighborhood in some sense. Let $T$ be a marked set. Let $P_s$ be the induced subcomplex on vertices whose undirected distance to $T$ is at most $s$. Here, distance is the minimum number of steps between vertices along allowed 1-paths, ignoring orientation. For each integer $s\ge1$, the finite-radius excision theorem (Theorem~\ref{thm:chain-locality}) shows that the inclusion of pairs $(P_s,P_s-T)\hookrightarrow(P,P-T)$ induces an isomorphism
\[
 H_q(P_s,P_s-T)\xrightarrow{\cong}\LH_q(P,T),
 \quad 0\le q<s.
\]
When $T=\{v\}$, the induced $s$-neighborhood of $v$ determines its local homology in all degrees less than $s$. At degree $q=s$, the comparison map is still surjective. At the chain level, the map is injective in all degrees and an isomorphism for $p\le s$. Both radius bounds are sharp, even for finite directed acyclic graphs.

Working with bounded neighborhoods also gives a stability property when we delete distant vertices. For $q\ge0$ and $A\subseteq V(G)\setminus N_{q+1}^G(v)$, the inclusion $G-A\hookrightarrow G$ induces an isomorphism $\LH_q(G-A,v)\xrightarrow{\cong}\LH_q(G,v)$ (Theorem~\ref{thm:deletion-local-stability}). This stability means local vanishing conditions still hold after we delete vertices far enough away.

Local homology for the union of marked sets is generally not the direct sum of the local homology of each individual set. Let $T_1,\ldots,T_n$ be pairwise disjoint nonempty sets with $n\ge2$, and set $T=\bigcup_iT_i$. The natural map from the local complex of $T$ to the direct sum of local complexes of each $T_i$ is injective. We define the \textit{interaction complex} $\mathcal I_*(P;T_1,\ldots,T_n)$ by the short exact sequence
\[
 0\longrightarrow Q_*(P,T)
 \xrightarrow{\rho}\bigoplus_{i=1}^nQ_*(P,T_i)
 \longrightarrow\mathcal I_*(P;T_1,\ldots,T_n)
 \longrightarrow0.
\]
Its homology measures how much local additivity fails. Let $\rho_r$ be the map induced by $\rho$ on homology. Let $d=\min_{i\ne j}d_P(T_i,T_j)$. Then $\rho_r$ is an isomorphism when $d>r+1$, and surjective when $d>r$ (Theorem~\ref{thm:separated-additivity}). These bounds are sharp. Interactions can involve arbitrarily many marked vertices. For example, take a graph made of $m\ge3$ paths $a\to x_i\to b$. The local homology vanishes for each proper subset of $\{x_1,\ldots,x_m\}$, but the whole set has nonzero local homology (Example~\ref{ex:higher-interaction}).

The interaction complex gives a simple condition for deleting multiple vertices at once. Suppose that for any $v\in T$, the inclusion $G-v\hookrightarrow G$ induces a homology isomorphism in each degree. Then the inclusion $G-T\hookrightarrow G$ induces a homology isomorphism in each degree if and only if the interaction complex built from these single-vertex sets is acyclic (Theorem~\ref{thm:parallel-interaction-criterion}).

Let $r\ge0$ and $T\subseteq V(G)$. Suppose that the local homology of $G$ at each $v\in T$ vanishes in degrees $0,\ldots,r+1$ and that
\[
 d_G(v,w)\ge r+3,
 \quad v,w\in T,\quad v\ne w.
\]
Then, for each $A\subseteq B\subseteq T$, the inclusion $G-B\hookrightarrow G-A$ induces a homology isomorphism
\[
 H_q(G-B)\xrightarrow{\cong}H_q(G-A),
 \quad 0\le q\le r
\]
(Theorem~\ref{thm:robust-parallel-deletion}). The vanishing condition at each $v$ can be checked in the induced ball $B_{r+2}^G(v)$.

Let $G$ be a finite digraph in the nearest-neighbor lattice $\mathbb Z^n$, with at most one orientation on each undirected edge. The unit cubes whose vertices and edges lie in $G$ and whose edges have a product orientation form a cubical complex $X(G)$. By Theorem~\ref{thm:grid-cubical-chains}, we have a natural chain isomorphism
\[
 C_*^{\mathrm{cub}}(X(G);\Bbbk)\xrightarrow{\cong}\Omega_*(G)
\]
with respect to subgraph inclusions. This isomorphism is compatible with induced vertex deletion. Thus, one has
\[
 \LH_q(G,v)\cong\widetilde H_{q-1}(\Lambda_v(G);\Bbbk),
 \quad q\ge0,
\]
where $\Lambda_v(G)$ denotes the cubical link of $v$ in $X(G)$.

Every cube of $X(G)$ passing through $v$ belongs to $B_n^G(v)$. Therefore, the inclusion $B_n^G(v)\hookrightarrow G$ induces the chain isomorphism $Q_*(B_n^G(v),v)\cong Q_*(G,v)$. Assume any $v\in T$ has a nonempty cubical link $\Lambda_v(G)$ whose reduced homology over $\Bbbk$ vanishes. If
\[
 d_G(v,w)\ge n+1,
 \quad v,w\in T,\quad v\ne w,
\]
then for any inclusion $G-B\hookrightarrow G-A$ with $A\subseteq B\subseteq T$, we have a homology isomorphism (Theorem~\ref{thm:grid-all-degree-deletion}).

Additional results consist of a local K\"unneth formula for rectangular marked sets (Theorem~\ref{thm:local-kunneth}) and local gluing formulas for digraphs glued along matching edges. The product formula holds because the Cartesian shuffle isomorphism commutes with the map induced by vertex deletion. This construction also yields examples that attain the optimal radius bounds. For gluing along a matching, the relative mapping cone provides suspension formulas and extends the gluing exact sequence to marked sets and interaction complexes. For the layered digraphs considered in \cite[Theorem~1.1]{GrigorLiu}, the prism decomposition establishes vanishing in degrees $q\ge2$, for both the boundary on all elementary paths and the regular boundary (Theorem~\ref{thm:temporal-vanishing}). We also investigate radius persistence and vertex deletion within filtrations.

The paper is organized as follows. Sections~\ref{sec:prelim}--\ref{sec:mv} introduce the chain complexes and prove locality properties, Mayer--Vietoris sequences, and excision. In Section~\ref{sec:star}, we study stars, links, and their comparison with local homology. Sections~\ref{sec:kunneth-radius}--\ref{sec:matching-gluing} discuss products, sharp radius bounds, geometric models, and gluing over matchings. In Sections~\ref{sec:separation} and~\ref{sec:parallel-deletion}, we study the interaction complexes and their applications to vertex deletion. The last section proposes some further problems.

\section{Path complexes and relative homology}\label{sec:prelim}

In this section, we review the basic concepts of regular path complexes and relative path homology. From now on, all digraphs considered are assumed to be finite with no loops. Opposite edges and directed cycles are permitted. The base field is set to be $\Bbbk$.

\subsection{Regular paths and path homology}

\begin{definition}\label{def:regular-paths}
Let $V$ be a vertex set and let $p\ge0$. An \textit{elementary $p$-path} on $V$ is a sequence $x_0\cdots x_p$ with $x_0,\ldots,x_p\in V$, denoted by $e_{x_0\cdots x_p}$. The path is called \textit{regular} if $x_i\ne x_{i+1}$ for each $0\le i<p$.
\end{definition}

\begin{definition}\label{def:regular-boundary}
For $p\ge0$, the \textit{regular path chain space} $\mathcal R_p(V)$ is the $\Bbbk$-vector space with basis the regular $p$-paths on $V$. We identify nonregular paths with zero. For $p\ge1$, the \textit{boundary operator} $\partial_p:\mathcal R_p(V)\to\mathcal R_{p-1}(V)$ is defined by
\begin{equation*}
 \partial e_{x_0\cdots x_p}
 =\sum_{i=0}^p(-1)^i e_{x_0\cdots\widehat{x_i}\cdots x_p}.
\end{equation*}
For the ordinary complex, we set $\partial_0=0$.
\end{definition}

The boundary operator is induced from the complex of all elementary paths by quotienting out its subcomplex of nonregular paths. Hence, one has $\partial^2=0$.

For a chain $c=\sum_\alpha c_\alpha e_\alpha$, its \textit{support} is the set of elementary paths $\operatorname{supp}(c)=\{\alpha:c_\alpha\ne0\}$. Its \textit{vertex support} is the set of vertices occurring in these paths. We say that $c$ is \textit{supported in} $U\subseteq V$ if its vertex support is contained in $U$.

\begin{definition}\label{def:regular-path-complex}
A \textit{regular path complex} $P$ on $V$ is a collection of regular paths on $V$ closed under deletion of the first vertex and deletion of the last vertex of every path of positive length. The paths belonging to $P$ are called \textit{allowed}. We denote by $\mathcal A_p(P)$ the $\Bbbk$-linear space generated by the allowed $p$-paths of $P$.
\end{definition}

In general, a path complex is not closed under deletion of an interior vertex. The path complex $P(G)$ of a digraph $G$ consists of all regular paths satisfying $x_{i-1}\to x_i\in E(G)$.

\begin{definition}
A regular path complex $P$ is called \textit{perfect} if, for any allowed path $x_0\cdots x_p$ with $p\ge1$ and any $0\le i\le p$, the path $x_0\cdots\widehat{x_i}\cdots x_p$ is allowed.
\end{definition}

The path complex associated with an ordered simplicial complex is perfect. Every induced path complex of a perfect path complex is also perfect.

\begin{definition}
Let $P$ be a regular path complex. The associated \textit{path chain complex} is defined by
\begin{equation*}
 \Omega_p(P):=\{c\in\mathcal A_p(P):\partial c\in\mathcal A_{p-1}(P)\},
 \quad \Omega_0(P)=\mathcal A_0(P).
\end{equation*}
Since $\partial^2=0$, we have $\partial\Omega_p(P)\subseteq\Omega_{p-1}(P)$. The \textit{path homology} of $P$ is defined by
\begin{equation*}
 H_p(P):=H_p(\Omega_*(P)).
\end{equation*}
For a digraph $G$, we write $\Omega_*(G)=\Omega_*(P(G))$ and $H_*(G)=H_*(P(G))$.
\end{definition}

This is the standard construction by $\partial$-invariant allowed chains \cite{GLMY}. The condition $\partial c\in\mathcal A_{p-1}(P)$ in the definition of $\Omega_p(P)$ applies to linear combinations of paths. An allowed path need not belong to $\Omega_p$.

\subsection{Relative complexes}

\begin{definition}\label{def:induced-digraph}
Let $G=(V,E)$ be a digraph and let $U\subseteq V$. The \textit{induced sub-digraph} of $G$ on $U$, denoted by $G[U]$, is the digraph with vertex set $U$ and edge set $E\cap(U\times U)$.
\end{definition}

\begin{definition}\label{def:induced-path-complex}
Let $P$ be a regular path complex on a vertex set $V$, and let $U\subseteq V$. The \textit{induced path complex} on $U$, denoted by $P|_U$, consists of all allowed paths of $P$ whose vertices lie in $U$.
\end{definition}

For a digraph $G$ and $U\subseteq V(G)$, we have $P(G)|_U=P(G[U])$.

\begin{lemma}\label{lem:intersection}
For any $p\ge0$, we have
\[
 \Omega_p(P|_U)=\Omega_p(P)\cap\mathcal R_p(U).
\]
Consequently, for $U,W\subseteq V$, one has $\Omega_*(P|_U)\cap\Omega_*(P|_W)=\Omega_*(P|_{U\cap W})$.
\end{lemma}
\begin{proof}
The inclusion $\Omega_p(P|_U)\subseteq\Omega_p(P)\cap\mathcal R_p(U)$ follows from the definition. Conversely, let $c\in\Omega_p(P)$ be supported on paths whose vertices lie in $U$. Its allowed paths belong to $P|_U$, and its boundary also uses only vertices of $U$. Since $\partial c\in\mathcal A_{p-1}(P)$, we obtain $\partial c\in\mathcal A_{p-1}(P|_U)$. Applying the equality to $U$, $W$, and $U\cap W$ gives the second assertion, since a regular path lies in both vertex sets exactly when it lies in their intersection.
\end{proof}

Let $P'\subseteq P$ be a sub-path complex. Then $\Omega_*(P')$ is a subcomplex of $\Omega_*(P)$. The \textit{relative chain complex} and \textit{relative path homology} are defined by
\[
 \Omega_*(P,P'):=\Omega_*(P)/\Omega_*(P'),\quad
 H_r(P,P'):=H_r(\Omega_*(P,P')).
\]
The short exact sequence yields the ordinary long exact homology sequence
\begin{equation*}
 \cdots\to H_r(P')\to H_r(P)\to H_r(P,P')
 \to H_{r-1}(P')\to\cdots .
\end{equation*}
Relative path homology and these exact sequences are part of the standard theory \cite[Theorem 3.11]{GJMY}.

For reduced path homology, we use the augmented complex with $\widetilde\Omega_{-1}(P)=\Bbbk$. Its augmentation is given by $\epsilon(\sum a_xe_x)=\sum a_x$. We retain this degree $-1$ term even when $P=\varnothing$. Thus, we have
\begin{equation*}
 \widetilde H_{-1}(\varnothing)=\Bbbk,\quad
 \widetilde H_{-1}(P)=0,\quad P\ne\varnothing.
\end{equation*}
The inclusion $P'\hookrightarrow P$ induces the identity on the degree-$-1$ terms of the augmented chain complexes, so $\widetilde\Omega_*(P)/\widetilde\Omega_*(P')=\Omega_*(P,P')$. Accordingly, both absolute terms in the long exact sequence of the pair may be replaced throughout by reduced homology, while relative homology remains unchanged. Ordinary chain complexes are zero in negative degrees.

\subsection{Envelope chain complex}

\begin{definition}\label{def:sigma}
For a regular path complex $P$, the boundary-completion chain complex is defined by
\begin{equation*}
 \Sigma_p(P):=\mathcal A_p(P)+\partial\mathcal A_{p+1}(P)
 \subseteq\mathcal R_p(V),\quad p\ge0.
\end{equation*}
This is also called the \textit{upper-envelope chain complex} of the allowed chain spaces.
\end{definition}

Since $\partial\Sigma_p(P)=\partial\mathcal A_p(P)\subseteq\Sigma_{p-1}(P)$, the spaces $\Sigma_*(P)$ form a chain complex. It is the smallest subcomplex containing all allowed chains, whereas $\Omega_*(P)$ is the largest subcomplex contained in the allowed chain spaces. The complex $\Sigma_*(P)$ is constructed inside the regular path chain complex by adjoining boundaries of allowed chains. This construction does not require a change in the edge set of a digraph.

\begin{theorem}\label{thm:sigma-quasi}
The natural inclusion
\begin{equation*}
 \Omega_*(P)\hookrightarrow\Sigma_*(P)
\end{equation*}
is a quasi-isomorphism. That is, it induces a homology isomorphism in each degree. For each sub-path complex $P'\subseteq P$, it also induces natural isomorphisms
\begin{equation*}
 H_r(P,P')\cong
 H_r\!\left(\Sigma_*(P)/\Sigma_*(P')\right).
\end{equation*}
\end{theorem}
\begin{proof}
Let $c=a+\partial b\in\Sigma_p(P)$ be a cycle, where $a\in\mathcal A_p(P)$ and $b\in\mathcal A_{p+1}(P)$. Then we have $\partial a=0$, and hence $a\in\Omega_p(P)$. Moreover, we have $c-a=\partial b$ and $b\in\Sigma_{p+1}(P)$. Thus, every homology class of $\Sigma_*(P)$ is represented by a cycle in $\Omega_*(P)$, which proves surjectivity.

If a cycle $a\in\Omega_p(P)$ is a boundary in $\Sigma$, write $a=\partial(b+\partial d)=\partial b$ with $b\in\mathcal A_{p+1}(P)$. Since $\partial b=a\in\mathcal A_p(P)$, we have $b\in\Omega_{p+1}(P)$. It follows that the inclusion $\Omega_*(P)\hookrightarrow\Sigma_*(P)$ also induces an injective map on homology. The construction commutes with inclusions. Finally, we compare the short exact sequences defining the relative complexes for $\Omega$ and $\Sigma$. The two absolute comparisons are quasi-isomorphisms, so the long exact sequences and the five lemma give the relative isomorphism.
\end{proof}

\section{Local path homology and finite-radius locality}\label{sec:local}

In this section, we define local path homology and prove that it can be computed in a neighborhood whose radius depends on the homological dimension.

\subsection{Local path homology}\label{sec:local-definition}

\begin{definition}[Local path homology at a vertex]\label{def:local}
Let $P$ be a regular path complex on $V$, and let $v\in V$. We write $P-v=P|_{V\setminus\{v\}}$. The \textit{local chain complex} and the \textit{local path homology} of $P$ at $v$ are defined by
\begin{equation*}
 Q_*(P,v):=\Omega_*(P)/\Omega_*(P-v),\quad
 \LH_r(P,v):=H_r(Q_*(P,v)).
\end{equation*}
For a digraph $G$, we apply this definition to $P(G)$ and use the notation $Q_*(G,v)$ and $\LH_r(G,v)$.
\end{definition}

For any vertex subset $T\subseteq V(P)$, we put $P-T=P|_{V(P)\setminus T}$ and use the notation
\begin{equation*}
 Q_*(P,T):=\Omega_*(P)/\Omega_*(P-T),\quad
 \LH_r(P,T):=H_r(Q_*(P,T)).
\end{equation*}
Thus, $Q_*(P,T)$ is the relative chain complex of the pair $(P,P-T)$. For $T=\varnothing$, this gives $Q_*(P,\varnothing)=0$. We discuss its relation to the local homology of geometric realizations in Section~\ref{sec:geometry}.

\begin{proposition}\label{prop:basic-local}
There is a natural long exact sequence
\begin{equation*}
 \cdots\to H_r(G-v)\to H_r(G)\to\LH_r(G,v)
 \to H_{r-1}(G-v)\to H_{r-1}(G)\to\cdots .
\end{equation*}
If $v$ is isolated, then we have $\LH_0(G,v)=\Bbbk$ and all other local homology groups vanish. If $v$ is not isolated, then we have $\LH_0(G,v)=0$.
\end{proposition}
\begin{proof}
The defining short exact sequence $0\to\Omega_*(P-v)\to\Omega_*(P)\to Q_*(P,v)\to0$ gives the stated long exact sequence. The space $Q_0$ is generated by $[e_v]$. Every edge incident to $v$ has relative boundary $\pm[e_v]$, so $H_0(Q)=0$ when $v$ is not isolated. If $v$ is isolated, no allowed path of positive length contains $v$, so $Q_p=0$ for all $p>0$.
\end{proof}

\begin{proposition}
Let $f:G\to G'$ be a digraph map such that $f(v)=v'$ and $f(V\setminus\{v\})\subseteq V'\setminus\{v'\}$. Then $f$ induces a chain map $Q_*(G,v)\to Q_*(G',v')$ and maps on local homology, compatible with identities and composition.
\end{proposition}
\begin{proof}
A digraph map may collapse an edge to a vertex. Applying $f$ vertexwise to regular chains and setting images with consecutive repeated vertices equal to zero gives the standard chain map $f_\#: \Omega_*(G)\to\Omega_*(G')$ \cite{GLMY,GJMY}. The stated hypothesis ensures that it maps $\Omega_*(G-v)$ into $\Omega_*(G'-v')$, and therefore descends to the quotient.
\end{proof}

For a map of pointed complexes, $f(v)=v'$ does not ensure a map of pairs: a second vertex may also map to $v'$. The stated naturality requires the complement condition.

\begin{example}\label{ex:interval}
Let $G$ be the interval $a\to v\to b$. The only allowed $2$-path is $avb$. Its boundary has the non-allowed term $ab$, and therefore $\Omega_2(G)=0$. The relative chain groups are
\[
 Q_1=\Bbbk[e_{av}]\oplus\Bbbk[e_{vb}],\quad Q_0=\Bbbk[e_v],
 \quad\bar\partial[e_{av}]=[e_v],\quad\bar\partial[e_{vb}]=-[e_v].
\]
Thus, we obtain $\LH_1(G,v)=\Bbbk$, generated by $[e_{av}+e_{vb}]$, and all other local homology groups vanish.

Adding the edge $a\to b$ gives $e_{avb}\in\Omega_2(G)$, with relative boundary $[e_{av}+e_{vb}]$. All local homology groups then vanish. Local path homology can therefore change when transitive edges are added.
\end{example}

\subsection{Finite-radius excision}

\begin{definition}
The \textit{distance} $d(x,y)$ between vertices $x,y\in V(P)$ is the length of a shortest path in the underlying undirected graph determined by the allowed $1$-paths of $P$. Vertices in different weakly connected components have distance $\infty$.
\end{definition}

\begin{definition}
For a nonempty vertex set $T\subseteq V(P)$, put $d(x,T)=\min_{t\in T}d(x,t)$. For $s\in\mathbb Z_{\ge0}$, we write
\[
 N_s(T)=\{x:d(x,T)\le s\},\quad P_s=P|_{N_s(T)},\quad P-T=P|_{V\setminus T}.
\]
For a digraph and $T=\{v\}$, we use the notation $B_s(v)=G[N_s(\{v\})]$, $B_s^\times(v)=B_s(v)-v$. The radius $s$ and the homological degree $r$ are distinct parameters.
\end{definition}

\begin{lemma}\label{lem:endpoints}
Let $p\ge1$, and let $c$ be an allowed $p$-chain. Write $c$ as the sum of its components indexed by the ordered pairs of endpoints $(a,b)$. Then $c\in\Omega_p(P)$ if and only if every endpoint component belongs to $\Omega_p(P)$. This holds for each regular path complex $P$.
\end{lemma}
\begin{proof}
For $p=1$, we have $\Omega_1(P)=\mathcal A_1(P)$, so the assertion follows. Suppose $p\ge2$, and write $c=\sum_{a,b}c_{a,b}$, where $c_{a,b}$ is the sum of all terms of $c$ with initial vertex $a$ and terminal vertex $b$. Let
\[
 \pi:\mathcal R_{p-1}(V)\longrightarrow \mathcal R_{p-1}(V)/\mathcal A_{p-1}(P)
\]
be the quotient map. For each ordered pair $(a,b)$, let $Q_{a,b}$ be the subspace of this quotient spanned by the classes of non-allowed regular $(p-1)$-paths with endpoints $(a,b)$. Since these classes form a basis of the quotient, we have
\[
 \mathcal R_{p-1}(V)/\mathcal A_{p-1}(P)=\bigoplus_{a,b\in V}Q_{a,b}.
\]
For any allowed $p$-path $x_0\cdots x_p$, closure under deletion of the first and last vertices gives
\[
 \pi\partial e_{x_0\cdots x_p}
 =\sum_{i=1}^{p-1}(-1)^i\pi e_{x_0\cdots\widehat{x_i}\cdots x_p}
 \in Q_{x_0,x_p},
\]
where nonregular paths represent zero. It follows that $\pi\partial c_{a,b}\in Q_{a,b}$. By the direct sum decomposition above, we obtain
\[
 \begin{aligned}
 c\in\Omega_p(P)
 &\Longleftrightarrow \pi\partial c=0\\
 &\Longleftrightarrow \pi\partial c_{a,b}=0\quad\text{for all }a,b\in V\\
 &\Longleftrightarrow c_{a,b}\in\Omega_p(P)\quad\text{for all }a,b\in V.
 \end{aligned}
\]
This completes the proof.
\end{proof}

\begin{theorem}[Finite-radius excision]\label{thm:chain-locality}\label{thm:mv-subset-radius}
Let $P$ be a regular path complex, and let $T\subseteq V(P)$ be a nonempty vertex set. Then the chain map induced by the inclusion of pairs $(P_s,P_s-T)\hookrightarrow(P,P-T)$
\[
 \Omega_*(P_s,P_s-T)\longrightarrow\Omega_*(P,P-T)
\]
is injective in each degree and is an isomorphism on chain groups in degrees $0\le p\le s$. Consequently, the induced homology map
\begin{equation*}
 H_r(P_s,P_s-T)\longrightarrow H_r(P,P-T)
\end{equation*}
is an isomorphism for $s\ge r+1$ and a surjection for $s=r$.
\end{theorem}
\begin{proof}
Put $C^s_*=\Omega_*(P_s,P_s-T)$ and $C_*=\Omega_*(P,P-T)$, and denote the chain map induced by $(P_s,P_s-T)\hookrightarrow(P,P-T)$ by $\iota_*:C^s_*\to C_*$. By Lemma~\ref{lem:intersection}, for any $p\ge0$, we have
\[
 \Omega_p(P_s)\cap\Omega_p(P-T)=\Omega_p(P_s-T).
\]
In degree zero, both $C^s_0$ and $C_0$ are canonically identified with the vector space generated by $T$, and $\iota_0$ is the identity under these identifications.

Let $1\le p\le s$, and let $c\in\Omega_p(P)$. Write $c=\sum_{a,b}c_{a,b}$ as in Lemma~\ref{lem:endpoints}, and put
\[
 I_p(T)=\{(a,b)\in V(P)^2:d(a,T)+d(b,T)\le p\},\quad
 c'=\sum_{(a,b)\in I_p(T)}c_{a,b}.
\]
By Lemma~\ref{lem:endpoints}, we have $c',c-c'\in\Omega_p(P)$. If an allowed path $x_0\cdots x_p$ satisfies $x_j\in T$ for some $0\le j\le p$, then we have
\[
 d(x_0,T)+d(x_p,T)\le j+(p-j)=p.
\]
Thus, every path occurring in $c-c'$ has all its vertices in $V(P)\setminus T$. Applying Lemma~\ref{lem:intersection}, we obtain
\[
 c-c'\in\Omega_p(P)\cap\mathcal R_p(V(P)\setminus T)
 =\Omega_p(P-T).
\]
For any path $x_0\cdots x_p$ occurring in $c'$ and any $0\le i\le p$, we have
\[
 \begin{aligned}
 d(x_i,T)
 &\le\min\{d(x_0,T)+i,\ d(x_p,T)+p-i\}\\
 &\le\frac{d(x_0,T)+d(x_p,T)+p}{2}\le p\le s.
 \end{aligned}
\]
It follows that
\[
 c'\in\Omega_p(P)\cap\mathcal R_p(N_s(T))=\Omega_p(P_s).
\]
Hence, we have $\iota_p([c'])=[c]$, which proves that $\iota_p$ is surjective for $0\le p\le s$.

Identify $C^s_*$ with its image in $C_*$ under $\iota_*$. For $s\ge r$, we have $C^s_r=C_r$. Since the differentials are compatible with inclusion and $\iota_{r-1}$ is injective, we obtain
\[
 Z_r(C^s_*)=Z_r(C_*),\quad
 \partial C^s_{r+1}\subseteq\partial C_{r+1}.
\]
For $r=0$, the equality of cycle spaces follows from $\partial_0=0$. The induced homology map is therefore the quotient map
\[
 \frac{Z_r(C_*)}{\partial C^s_{r+1}}
 \longrightarrow
 \frac{Z_r(C_*)}{\partial C_{r+1}},
\]
which is surjective. If $s\ge r+1$, then $C^s_{r+1}=C_{r+1}$, so this map is an isomorphism.
\end{proof}

\begin{corollary}\label{cor:radius}
Let $G$ be a digraph, and let $v$ be a vertex of $G$. For any $r\ge0$, the map induced by the inclusion of pairs $(B_s(v),B_s^\times(v))\hookrightarrow(G,G-v)$
\[
H_r(B_s(v),B_s^\times(v))\longrightarrow\LH_r(G,v)
\]
is an isomorphism for $r<s$ and a surjection for $r=s$.
\end{corollary}
\begin{proof}
The result follows from Theorem~\ref{thm:chain-locality} by taking $P=P(G)$ and $T=\{v\}$.
\end{proof}

\begin{remark}
The bounds on $s$ cannot be decreased in general, even when $P=P(G)$ for a finite digraph $G$ and $T=\{v\}$. More precisely, Corollary~\ref{cor:iso-radius-optimal} shows that, for any $r\ge0$, there exists a finite digraph $G$ with a vertex $v$ such that the map induced by $(B_r(v),B_r^\times(v))\hookrightarrow(G,G-v)$
\[
 H_r(B_r(v),B_r^\times(v))\longrightarrow\LH_r(G,v)
\]
is not injective. Corollary~\ref{cor:surj-radius-optimal} shows that, for any $r\ge1$, there exists a finite digraph $G$ with a vertex $v$ such that the map induced by $(B_{r-1}(v),B_{r-1}^\times(v))\hookrightarrow(G,G-v)$
\[
 H_r(B_{r-1}(v),B_{r-1}^\times(v))\longrightarrow\LH_r(G,v)
\]
is not surjective. Thus, $r+1$ and $r$ are the smallest nonnegative integer radii that guarantee, respectively, an isomorphism and a surjection in degree $r$ for each pair $(G,v)$. For a perfect regular path complex, Proposition~\ref{prop:perfect-radius-one} shows that radius one determines the entire local chain complex at each vertex.
\end{remark}

Let $r\ge0$, and put $Q_*=Q_*(G,v)$ and $Q^s_*=Q_*(B_s(v),v)$. By Theorem~\ref{thm:chain-locality}, the maps induced by the inclusions $B_s(v)\hookrightarrow G$ identify $Q^s_*$ with subcomplexes of $Q_*$. Under these identifications, we have
\[
 Q^r_*\subseteq Q^{r+1}_*\subseteq Q_*,\quad Q^r_r=Q^{r+1}_r=Q_r.
\]
The boundary spaces $\partial(Q^r_{r+1})$ and $\partial(Q^{r+1}_{r+1})$ are subspaces of the same chain space $Q_r$. Let
\[
 \mathcal K_r(v)=\frac{\partial(Q^{r+1}_{r+1})}{\partial(Q^r_{r+1})}.
\]

\begin{proposition}\label{prop:shell}
Then there is a short exact sequence
\begin{equation*}
 0\longrightarrow\mathcal K_r(v)
 \xrightarrow{\alpha_r} H_r(B_r(v),B_r^\times(v))
 \xrightarrow{\beta_r}\LH_r(G,v)\longrightarrow0.
\end{equation*}
For $b\in Q^{r+1}_{r+1}$ and $z\in Z_r(Q^r_*)$, the maps are
\[
 \alpha_r\colon [\partial b]\longmapsto[\partial b]_{B_r(v)},\quad
 \beta_r\colon [z]_{B_r(v)}\longmapsto[z]_G.
\]
\end{proposition}
\begin{proof}
In view of $Q^r_r=Q_r$, we have $Z_r(Q^r_*)=Z_r(Q_*)$. Denote $Z_r(Q^r_*)$ by $Z_r$. Theorem~\ref{thm:chain-locality} also gives $Q^{r+1}_{r+1}=Q_{r+1}$. It follows that
\[
 \partial(Q^r_{r+1})\subseteq\partial(Q^{r+1}_{r+1})=\partial(Q_{r+1})\subseteq Z_r.
\]
Thus, we obtain
\[
 H_r(B_r(v),B_r^\times(v))=\frac{Z_r}{\partial(Q^r_{r+1})},\quad
 \LH_r(G,v)=\frac{Z_r}{\partial(Q_{r+1})}.
\]
Under these identifications, $\alpha_r$ is the inclusion of the subspace $\partial(Q_{r+1})/\partial(Q^r_{r+1})$ into $Z_r/\partial(Q^r_{r+1})$, and $\beta_r$ is the quotient map
\[
 z+\partial(Q^r_{r+1})\longmapsto z+\partial(Q_{r+1}),\quad z\in Z_r.
\]
In particular, $\alpha_r$ is injective and $\beta_r$ is surjective. Moreover, we have
\[
 \ker\beta_r=\frac{\partial(Q_{r+1})}{\partial(Q^r_{r+1})}
 =\frac{\partial(Q^{r+1}_{r+1})}{\partial(Q^r_{r+1})}
 =\operatorname{im}\alpha_r.
\]
This proves the short exact sequence.
\end{proof}

\subsection{Weakly connected components}

\begin{definition}
Let $G=(V,E)$ be a digraph. The \textit{underlying undirected graph} of $G$ has vertex set $V$, with an undirected edge $\{x,y\}$ whenever $(x,y)\in E$ or $(y,x)\in E$.
\end{definition}

\begin{definition}
A \textit{weakly connected component} of $G$ is the induced sub-digraph $G[U]$, where $U$ is the vertex set of a connected component of the underlying undirected graph.
\end{definition}

For $v\in V$, we denote by $C_v$ the unique weakly connected component containing $v$.

\begin{theorem}\label{thm:component}
Let $C_v$ be the weakly connected component containing $v$. Then the map induced by the inclusion $C_v\hookrightarrow G$
\[
 Q_*(C_v,v)\xrightarrow{\cong}Q_*(G,v)
\]
is an isomorphism of chain complexes and therefore induces a local homology isomorphism in each degree.
\end{theorem}
\begin{proof}
Let $F$ be the union of the other components. Every allowed path lies in a single weakly connected component. The boundary of a path uses only vertices of that component, so the condition of having an allowed boundary can be checked separately in each component. Hence, we have
\[
 \Omega_*(G)=\Omega_*(C_v)\oplus\Omega_*(F),\quad
 \Omega_*(G-v)=\Omega_*(C_v-v)\oplus\Omega_*(F).
\]
On taking the quotient, the summand $\Omega_*(F)$ disappears from both sides. We obtain the required chain isomorphism.
\end{proof}

\begin{definition}
Let $G=(V,E)$ be a digraph. Two vertices $x,y\in V$ are \textit{mutually reachable} if there is a directed path from $x$ to $y$ and a directed path from $y$ to $x$, where paths of length zero are permitted.
\end{definition}

Mutual reachability is an equivalence relation on $V$. A \textit{strongly connected component} of $G$ is the induced sub-digraph on an equivalence class of this relation.

The conclusion of Theorem~\ref{thm:component} need not hold for the strongly connected component containing $v$.

\begin{example}
Let $G$ be the directed interval $a\to v\to b$, and let $D_v=G[\{v\}]$. Then $D_v$ is the strongly connected component containing $v$, but
\[
 \LH_1(D_v,v)=0,\quad \LH_1(G,v)\cong\Bbbk.
\]
Consequently, the inclusion $D_v\hookrightarrow G$ does not induce a local homology isomorphism in degree one.
\end{example}

Let $r\ge0$ and $s\ge r+1$. Recall that $B_s^\times(v)=B_s(v)-v$ is the induced sub-digraph on $N_s(\{v\})\setminus\{v\}$. By Corollary~\ref{cor:radius}, the inclusion of pairs $(B_s(v),B_s^\times(v))\hookrightarrow(G,G-v)$ induces an isomorphism
\[
 H_r(B_s(v),B_s^\times(v))\xrightarrow{\cong}\LH_r(G,v).
\]
Regarding the reduced long exact sequence of the pair $(B_s(v),B_s^\times(v))$, we obtain the exact sequence
\[
 \widetilde H_r(B_s(v))\longrightarrow\LH_r(G,v)
 \xrightarrow{\delta_{s,r}}\widetilde H_{r-1}(B_s^\times(v))
 \longrightarrow\widetilde H_{r-1}(B_s(v)).
\]
Here, $\delta_{s,r}$ is the connecting homomorphism. By exactness, the conditions
\[
 \widetilde H_r(B_s(v))=0,\quad \widetilde H_{r-1}(B_s(v))=0
\]
are sufficient for $\delta_{s,r}$ to be an isomorphism. No vanishing assumption in other degrees is required for this conclusion.

A vertex $w\ne v$ belongs to $B_s^\times(v)$ if and only if the underlying undirected graph contains a path of length at most $s$ from $v$ to $w$. There may be no directed path in $G$ from $v$ to $w$ or from $w$ to $v$. Deleting such a vertex $w$ can change the local homology at $v$.

\begin{example}
Consider the diamond, whose edges are $a\to v$, $v\to b$, $a\to w$, and $w\to b$. The edges $a\to v$ and $a\to w$ give $d(v,w)=2$, although no directed path joins $v$ and $w$ in either direction. In this case, we have
\[
 \LH_1(G-w,v)\cong\Bbbk,\quad \LH_1(G,v)=0.
\]
Indeed, the chain $e_{avb}-e_{awb}\in\Omega_2(G)$ has relative boundary $[e_{av}+e_{vb}]$ in $Q_*(G,v)$. After deleting $w$, this degree-one cycle represents a nonzero local homology class.
\end{example}

Therefore, restricting a ball to vertices joined to $v$ by a directed path in either direction does not in general preserve local homology.

\subsection{Computation on directed trees}

A \textit{directed tree} is obtained from an undirected tree by assigning exactly one direction to each edge.

\begin{proposition}\label{prop:tree}
Let $G$ be a directed tree. For each nonisolated vertex $v$ of $G$, we have
\[
 \LH_r(G,v)\cong
 \begin{cases}
 \Bbbk^{\deg(v)-1},&r=1,\\
 0,&r\ne1.
 \end{cases}
\]
\end{proposition}
\begin{proof}
Let $p\ge2$ and let $x_0\cdots x_p$ be an allowed path. The vertices $x_0$ and $x_2$ are distinct, and their unique undirected path is $x_0,x_1,x_2$. Deleting $x_1$ gives the regular non-allowed face $x_0x_2\cdots x_p$. Its only missing edge is $x_0\to x_2$, and $x_1$ is the unique vertex that can fill this gap. Hence, the coefficient of this face in the boundary of an allowed chain is the negative of the coefficient of $x_0\cdots x_p$. It follows that $\Omega_p(G)=0$ for $p\ge2$.

Put $Q_*=Q_*(G,v)$ and $d=\deg(v)$. Then $Q_0=\Bbbk[e_v]$, and $Q_p=0$ for $p\ge2$. The $d$ edges $f_1,\ldots,f_d$ incident to $v$ give a basis of $Q_1$. Set $\varepsilon_i=1$ when $f_i$ ends at $v$, and $\varepsilon_i=-1$ when $f_i$ starts at $v$. Then
\[
 \partial[f_i]=\varepsilon_i[e_v].
\]
The local chain complex in these bases is
\[
 0\longrightarrow\Bbbk^d
 \xrightarrow{\ D\ }\Bbbk\longrightarrow0,\quad
 D(t_1,\ldots,t_d)=\sum_{i=1}^d\varepsilon_i t_i,
\]
with $\Bbbk^d$ in degree one and $\Bbbk$ in degree zero. Since $v$ is not isolated, we have $d\ge1$. It follows that $D$ is surjective and $\dim_{\Bbbk}\ker D=d-1$. So one has
\[
 \LH_1(G,v)=\ker D\cong\Bbbk^{d-1},\quad
 \LH_0(G,v)=\operatorname{coker}D=0.
\]
The homology groups vanish for $r\ge 2$. This completes the proof.
\end{proof}

\section{Mayer--Vietoris sequences and excision}\label{sec:mv}

Let $P$ be a regular path complex, and let $U,W\subseteq V(P)$. Throughout this section, we set
\begin{equation*}
 X=P|_{U\cup W},\quad A=P|_U,\quad B=P|_W,\quad Z=P|_{U\cap W}.
\end{equation*}
In particular, for a digraph, a restriction is the path complex of the corresponding induced subgraph on the subset.

\subsection{Path covers}

Consider the two conditions
\begin{align*}
 (\mathrm{P})\quad&
 \mathcal A_p(X)=\mathcal A_p(A)+\mathcal A_p(B),
 \quad\forall p\ge0;\\
 (\mathrm{O})\quad&
 \Omega_p(X)=\Omega_p(A)+\Omega_p(B),
 \quad\forall p\ge0.
\end{align*}
Condition $(\mathrm{P})$ requires each allowed path to lie in $A$ or in $B$. Condition $(\mathrm{O})$ indicates that each $\partial$-invariant chain is a sum of $\partial$-invariant chains on the two sides. It is worth noting that the vertex cover $U\cup W$ alone ensures neither condition.

\begin{example}\label{ex:mv-independent}
Even for finite directed acyclic graphs, neither $(\mathrm{P})$ nor $(\mathrm{O})$ implies the other.

We first show that $(\mathrm{O})$ does not imply $(\mathrm{P})$. Take the directed interval $X:a\to b\to c$, without $a\to c$, and put $U=\{a,b\}$, $W=\{b,c\}$. Its only allowed $2$-path is $abc$, and
\[
 \partial e_{abc}=e_{bc}-e_{ac}+e_{ab}
\]
contains the non-allowed edge $ac$. Hence, we have $\Omega_2(X)=0$, and the higher chain groups also vanish. All vertices and edges are covered by the two sides, so $(\mathrm{O})$ holds in each degree. The path $abc$ lies in neither side, however, and $(\mathrm{P})$ fails.

To show that $(\mathrm{P})$ does not imply $(\mathrm{O})$, consider the diamond digraph
\[
 a\to x\to b,\quad a\to y\to b,
 \quad U=\{a,x,b\},\quad W=\{a,y,b\}.
\]
Every allowed path lies in one side, so $(\mathrm{P})$ holds. Both sides have $\Omega_2(A)=\Omega_2(B)=0$, whereas
\[
 e_{axb}-e_{ayb}\in\Omega_2(X)\setminus\{0\},
\]
Thus, $(\mathrm{O})$ fails.
\end{example}

\begin{proposition}\label{prop:mv-independent}
If $X$ is a perfect path complex, in particular if it arises from an ordered simplicial complex, then conditions $(\mathrm{P})$ and $(\mathrm{O})$ are equivalent.
\end{proposition}
\begin{proof}
If $X$ is perfect, then its induced path complexes $A$ and $B$ are perfect as well. Hence, we have $\Omega_\ast(Y)=\mathcal A_\ast(Y)$ for $Y=X,A,B$. Thus, the two conditions are equivalent.
\end{proof}

Example~\ref{ex:mv-independent} shows that the two cover conditions are independent in general.

\subsection{Defect complexes and excision}

\begin{definition}\label{def:mv-defect}
Without assuming $(\mathrm{P})$ or $(\mathrm{O})$, we define the \textit{sum complex} and the \textit{defect complex} of the cover by
\begin{equation*}
 M_*:=\Omega_*(A)+\Omega_*(B),\quad
 \mathcal E_*:=\frac{\Omega_*(X)}{M_*}.
\end{equation*}
The complex $M_*$ is a subcomplex of $\Omega_*(X)$. In general, it cannot be identified with $\Omega_*(A\cup B)$, where the latter union is the union of the sets of allowed paths.
\end{definition}

\begin{theorem}\label{thm:mv-defect}
We have the following short exact sequences
\begin{align*}
 0&\longrightarrow\Omega_*(Z)
 \xrightarrow{c\mapsto(c,-c)}\Omega_*(A)\oplus\Omega_*(B)
 \xrightarrow{(a,b)\mapsto a+b}M_*\longrightarrow0,
 \\
 0&\longrightarrow M_*\longrightarrow\Omega_*(X)
 \longrightarrow\mathcal E_*\longrightarrow0.
\end{align*}
\end{theorem}
\begin{proof}
For $p\ge 0$, the addition map $(a,b)\mapsto a+b$ has image $M_p=\Omega_p(A)+\Omega_p(B)$. Its kernel consists of the pairs $(a,-a)$ with $a\in\Omega_p(A)\cap\Omega_p(B)=\Omega_p(Z)$ by Lemma~\ref{lem:intersection}. This is exactly the image of the injection $c\mapsto(c,-c)$. In the second sequence, the quotient map $\Omega_p(X)\to\mathcal E_p$ is surjective with kernel $M_p$. All these maps commute with $\partial$, so both sequences are short exact sequences of chain complexes.
\end{proof}

\begin{corollary}\label{cor:mv-classical}
Suppose $H_*(\mathcal E_*)=0$. Then there is a long exact sequence
\begin{equation*}
 \cdots\to H_r(Z)\to H_r(A)\oplus H_r(B)
 \to H_r(X)\xrightarrow{\delta_r}H_{r-1}(Z)\to\cdots.
\end{equation*}
Here, for a cycle $c=a+b\in\Omega_r(X)$ with $a\in\Omega_r(A)$ and $b\in\Omega_r(B)$, the connecting map is given by $\delta_r([c])=[\partial a]$.
\end{corollary}
\begin{proof}
The second sequence of Theorem~\ref{thm:mv-defect} and the vanishing of $H_*(\mathcal E_*)$ give an isomorphism $H_r(i):H_r(M_*)\to H_r(X)$. Taking homology gives the desired long exact sequence. The connecting map is given by $\delta_r=\delta_r^M\circ H_r(i)^{-1}$, where $\delta_r^M:H_r(M_*)\to H_{r-1}(Z)$ is the original connecting map. For $c=a+b \in \Omega_r(X)$, we have $\partial(a,b)=(\partial a,-\partial a)$. The injection from $\Omega_{r-1}(Z)$ sends $\partial a$ to $(\partial a,-\partial a)$, giving $\delta_r([c])=[\partial a]$.
\end{proof}

\begin{remark}
This construction is compatible with existing conditional Mayer--Vietoris theorems for path homology. See
\cite[Theorem 3.25]{GJMY}.
The complex $\mathcal E_*$ allows us to express the obstruction to the ordinary sequence without assuming an invariant-chain cover.
\end{remark}

\begin{theorem}\label{thm:mv-excision}
The inclusion of pairs $(A,Z)\hookrightarrow(X,B)$ induces a chain map $\iota$ and leads to a short exact sequence
\begin{equation*}
 0\longrightarrow\Omega_*(A,Z)
 \xrightarrow{\iota}\Omega_*(X,B)
 \longrightarrow\mathcal E_*\longrightarrow0.
\end{equation*}
The map $\iota$ is a chain isomorphism if and only if $(\mathrm{O})$ holds. Moreover, there is a long exact sequence
\begin{equation*}
 \cdots\to H_{r+1}(\mathcal E_*)\to H_r(A,Z)
 \xrightarrow{\iota_r}H_r(X,B)
 \to H_r(\mathcal E_*)\to H_{r-1}(A,Z)\to\cdots.
\end{equation*}
\end{theorem}
\begin{proof}
By Lemma~\ref{lem:intersection}, we have $\Omega_p(A)\cap\Omega_p(B)=\Omega_p(Z)$. Thus, the kernel of $\iota$ is zero. Besides, the cokernel of $\iota$ is $\Omega_p(X)/(\Omega_p(A)+\Omega_p(B))=\mathcal E_p$. Thus, $\iota$ is a chain isomorphism when $(\mathrm O)$ holds. The short exact sequence induces the desired long exact sequence. 
\end{proof}

In the diamond of Example~\ref{ex:mv-independent}, we have $\mathcal E_2=\Bbbk$ and all other degrees vanish. In this example, $\widetilde H_q(A)=\widetilde H_q(B)=\widetilde H_q(X)=0$ for each $q\ge-1$, whereas $Z$ consists of two isolated vertices. An unconditional ordinary reduced Mayer--Vietoris sequence would therefore require exactness of $0\to0\to\Bbbk\to0$, which is impossible. The corresponding excision map in degree one is
\[
H_1(A,Z)=\Bbbk\to H_1(X,B)=0.
\]

\begin{definition}
Let $X$ be a path complex. A vertex subset $U\subseteq V(X)$ is called \textit{convex} with respect to the allowed paths of $X$ if, for any allowed path $x_0\cdots x_p$ with $x_0,x_p\in U$, we have $x_i\in U$ for each $0\le i\le p$.
\end{definition}

\begin{proposition}\label{prop:mv-convex-cover}
Let $X$ be a regular path complex, and let $U,W\subseteq V(X)$ satisfy $U\cup W=V(X)$. Put $A=X|_U$ and $B=X|_W$. Suppose that every allowed path of $X$ lies entirely in $U$ or entirely in $W$, so that the cover by $A$ and $B$ satisfies $(\mathrm{P})$. If $U$ is convex with respect to the allowed paths of $X$, then we have
\[
 \Omega_p(X)=\Omega_p(X|_U)+\Omega_p(X|_W)
 =\Omega_p(A)+\Omega_p(B),\quad p\ge0.
\]
Thus, the same cover satisfies $(\mathrm{O})$.
\end{proposition}
\begin{proof}
For $p\ge1$, decompose $c\in\Omega_p(X)$ by endpoints, and let $c_U$ be the sum of components with both endpoints in $U$. By Lemma~\ref{lem:endpoints} and convexity, we have $c_U\in\Omega_p(A)$. No path in the other components lies in $U$. As condition $(\mathrm{P})$ places each such path in $W$, we have $c-c_U\in\Omega_p(B)$. The assertion in degree zero follows from the vertex cover.
\end{proof}

In particular, let $G=(V,E)$ be a digraph and let $U\cup W=V$. The cover by $P(G[U])$ and $P(G[W])$ satisfies both $(\mathrm{P})$ and $(\mathrm{O})$ if either
\[
 E\cap\bigl(S\times(V\setminus S)\bigr)=\varnothing,\quad \forall S\in\{U,W\}
\]
or
\[
 E\cap\bigl((V\setminus S)\times S\bigr)=\varnothing,\quad \forall S\in\{U,W\}.
\]

Under the first condition, take an allowed path $x_0\cdots x_p$ and choose $S\in\{U,W\}$ with $x_0\in S$. As all edges lie in $S$, one has $x_i\in S$ for $0\le i\le p$. This proves $(\mathrm{P})$. Taking $S=U$ also proves convexity of $U$. Under the second condition, start with $S$ containing $x_p$. The same argument proves $(\mathrm{P})$ and convexity of $U$. Proposition~\ref{prop:mv-convex-cover} now gives $(\mathrm{O})$.

\subsection{Completed Mayer--Vietoris sequence}\label{sec:mv-completion}

We now use the boundary-completion model of Definition~\ref{def:sigma}.

\begin{definition}\label{def:completed-intersection}
With the notation fixed at the beginning of this section, we define the \textit{completed intersection complex} $K_*(A,B)$ and its quotient $\mathcal T_*$ by
\begin{equation*}
 K_*(A,B):=\Sigma_*(A)\cap\Sigma_*(B),\quad
 \mathcal T_*:=K_*(A,B)/\Sigma_*(Z).
\end{equation*}
The intersection and quotient are taken in the category of chain complexes.
\end{definition}

Since $Z=A\cap B$, we have $\Sigma_p(Z)\subseteq K_p(A,B)$. Each chain in $K_p(A,B)$ is supported in both $U$ and $W$, so it lies in $\mathcal R_p(U\cap W)$.

The restriction formula in Lemma~\ref{lem:intersection} can fail for $\Sigma_p$. For example, consider the digraph $a\to x\to b$ with $U=\{a,x,b\}$ and $W=\{a,b\}$. Then $e_{ab}=e_{ax}+e_{xb}-\partial e_{axb}$ lies in $\Sigma_1(X)$ and is supported in $U\cap W$. But one has $\Sigma_1(Z)=0$.

\begin{theorem}\label{thm:completed-mv}
Assume $(\mathrm{P})$. Then we have
\begin{equation*}
 \Sigma_*(X)=\Sigma_*(A)+\Sigma_*(B).
\end{equation*}
Hence, there is a short exact sequence of chain complexes
\begin{equation*}
 0\to K_*(A,B)\xrightarrow{k\mapsto(k,-k)}
 \Sigma_*(A)\oplus\Sigma_*(B)
 \xrightarrow{(a,b)\mapsto a+b}\Sigma_*(X)\to0.
\end{equation*}
Moreover, there are isomorphisms
\begin{align*}
 \Sigma_*(A)/K_*(A,B)&\cong\Sigma_*(X)/\Sigma_*(B),
 \\
 H_r\!\left(\Sigma_*(A)/K_*(A,B)\right)&\cong H_r(X,B).
\end{align*}
\end{theorem}
\begin{proof}
By $(\mathrm{P})$, each allowed chain in degree $p$ or $p+1$ is a sum of chains in $A$ and $B$. Applying $(\mathrm{P})$ to $\Sigma_p(X)=\mathcal A_p(X)+\partial\mathcal A_{p+1}(X)$, one has $\Sigma_p(X)=\Sigma_p(A)+\Sigma_p(B)$. Therefore, the addition map is surjective, and its kernel consists of the pairs $(k,-k)$ with $k\in K_p(A,B)$.

The same sum decomposition shows that $\Sigma_*(A)\to\Sigma_*(X)/\Sigma_*(B)$ is surjective, with kernel $K_*(A,B)$. It induces the stated chain isomorphism. The relative part of Theorem~\ref{thm:sigma-quasi} identifies the homology of its target with $H_*(X,B)$.
\end{proof}

\begin{remark}
Let $U\cup W=V(X)$, $A=X|_U$, and $B=X|_W$. No assumption $(\mathrm{P})$ is made. Let $X^{\mathrm{cov}}$ be the set of allowed paths contained in $U$ or in $W$. Deleting either endpoint keeps a path on the same side. Hence, $X^{\mathrm{cov}}$ is a path complex. Moreover, we have
\[
 X^{\mathrm{cov}}|_U=A,\quad X^{\mathrm{cov}}|_W=B,\quad
 X^{\mathrm{cov}}=X\ \Longleftrightarrow\ (\mathrm{P}).
\]
The cover of $X^{\mathrm{cov}}$ by $A$ and $B$ satisfies $(\mathrm{P})$. Theorem~\ref{thm:completed-mv} gives a Mayer--Vietoris sequence for $X^{\mathrm{cov}}$ with terms from $A$, $B$, and the completed intersection.
\end{remark}

\subsection{Characterization of the defect complex}

In this subsection, let $X$ be a regular path complex with vertex set $V=V(X)$, and let $U,W\subseteq V$ satisfy $U\cup W=V$. We keep the notation
\[
 A=X|_U,\quad B=X|_W,\quad Z=X|_{U\cap W}.
\]
We assume condition $(\mathrm{P})$, namely that every allowed path of $X$ lies entirely in $U$ or entirely in $W$. For any $p\ge0$, let $\mathcal N_p(X)$ be the subspace of $\mathcal R_p(V)$ spanned by the regular elementary $p$-paths on $V$ that are not allowed in $X$. Then we have the decomposition
\[
 \mathcal R_p(V)=\mathcal A_p(X)\oplus\mathcal N_p(X).
\]
Let $\nu_p:\mathcal R_p(V)\to\mathcal N_p(X)$ be the projection onto the second summand. For $Y\in\{A,B,Z\}$, regard $\mathcal A_p(Y)$ as a subspace of $\mathcal R_p(V)$ and define
\begin{equation*}
 F_p(Y):=\nu_{p-1}\partial\mathcal A_p(Y)\subseteq\mathcal N_{p-1}(X),\quad p\ge1.
\end{equation*}
By definition, we have $\ker(\nu_{p-1}\partial|_{\mathcal A_p(Y)})=\Omega_p(Y)$.

\begin{proposition}\label{prop:illegal-face-defect}
For $p\ge1$, there is a natural isomorphism of vector spaces
\begin{equation*}
 \mathcal E_p\cong
 \frac{F_p(A)\cap F_p(B)}{F_p(Z)}.
\end{equation*}
Moreover, we have
\[
 F_p(A)\cap F_p(B)\subseteq\mathcal N_{p-1}(X)\cap\mathcal R_{p-1}(U\cap W).
\]
Thus, every elementary path with nonzero coefficient in a chain belonging to $F_p(A)\cap F_p(B)$ has all its vertices in $U\cap W$ and is not allowed in $X$.

Use primes for a second induced cover satisfying $(\mathrm{P})$. Let $f:X\to X'$ be a path-complex morphism with $f(U)\subseteq U'$ and $f(W)\subseteq W'$. Under the above isomorphism, the map $\mathcal E_p\to\mathcal E'_p$ given by $[c]\mapsto[f_\#c]$ corresponds to
\[
 \bar f_p:\frac{F_p(A)\cap F_p(B)}{F_p(Z)}
 \longrightarrow\frac{F'_p(A')\cap F'_p(B')}{F'_p(Z')},
 \quad [z]\longmapsto[\nu'_{p-1}f_\# z],
\]
where $\nu'_{p-1}$ is the projection onto $\mathcal N_{p-1}(X')$.
\end{proposition}
\begin{proof}
Put $D=\nu_{p-1}\partial$. By $(\mathrm{P})$, every $c\in\Omega_p(X)$ admits a decomposition $c=a+b$ with $a\in\mathcal A_p(A)$ and $b\in\mathcal A_p(B)$. Since $Da=-Db$, define
\[
 \Phi_p:\Omega_p(X)\longrightarrow
 \frac{F_p(A)\cap F_p(B)}{F_p(Z)},\quad
 c\longmapsto Da+F_p(Z).
\]
Any two choices of $a$ differ by an element of $\mathcal A_p(A)\cap\mathcal A_p(B)=\mathcal A_p(Z)$, so $\Phi_p$ is a well-defined linear map. For $\eta\in F_p(A)\cap F_p(B)$, choose $a\in\mathcal A_p(A)$ and $b\in\mathcal A_p(B)$ with $Da=Db=\eta$. Then $a-b\in\Omega_p(X)$ and $\Phi_p(a-b)=\eta+F_p(Z)$, proving surjectivity.

Note that $\Omega_p(A)+\Omega_p(B)\subseteq\ker\Phi_p$. Conversely, if $c=a+b\in\ker\Phi_p$, choose $z\in\mathcal A_p(Z)$ with $Dz=Da$. Then we obtain
\[
 D(a-z)=D(b+z)=0,\quad
 c=(a-z)+(b+z)\in\Omega_p(A)+\Omega_p(B).
\]
Thus, one has $\ker\Phi_p=\Omega_p(A)+\Omega_p(B)$, and the first isomorphism theorem gives the asserted isomorphism.

Neither $\partial$ nor the projection $\nu_{p-1}$ enlarges the vertex support of a chain. Hence, one has
\[
 F_p(A)\cap F_p(B)
 \subseteq\mathcal N_{p-1}(X)\cap\mathcal R_{p-1}(U)\cap\mathcal R_{p-1}(W)
 =\mathcal N_{p-1}(X)\cap\mathcal R_{p-1}(U\cap W),
\]
as desired.

Put $D'=\nu'_{p-1}\partial$. Since $f_\#$ preserves allowed chains and commutes with $\partial$, we have
\[
 \nu'_{p-1}f_\#D
 =\nu'_{p-1}f_\#\nu_{p-1}\partial
 =\nu'_{p-1}f_\#\partial
 =D'f_\#.
\]
Applying this identity to allowed chains in $A$, $B$, and $Z$ gives $\nu'_{p-1}f_\#(F_p(Y))\subseteq F'_p(Y')$ for $Y=A,B,Z$, so $\bar f_p$ is well defined. For $c=a+b$ as above, we obtain $\bar f_p(\Phi_p(c))=[D'f_\#a]=\Phi'_p(f_\#c)$, which proves naturality.
\end{proof}

Proposition~\ref{prop:illegal-face-defect} identifies the degree-$p$ defect space $\mathcal E_p$, rather than its homology. More precisely, for any $p\ge1$, we have
\[
 \Omega_p(X)=\Omega_p(A)+\Omega_p(B)
 \quad\Longleftrightarrow\quad
 F_p(A)\cap F_p(B)=F_p(Z).
\]
An element $f\in F_p(A)\cap F_p(B)$ can be written as
\[
 f=\nu_{p-1}\partial a=\nu_{p-1}\partial b,
 \quad a\in\mathcal A_p(A),\quad b\in\mathcal A_p(B).
\]
Its class in the quotient is nonzero precisely when there is no $z\in\mathcal A_p(Z)$ such that $f=\nu_{p-1}\partial z$. Although $f$ is supported in $U\cap W$, it need not be the non-allowed part of the boundary of an allowed chain in $Z$.

To compute $H_p(\mathcal E_*)$, one must also use the quotient differential
\[
 d_{\mathcal E,p}:\mathcal E_p\longrightarrow\mathcal E_{p-1},\quad
 [c]\longmapsto[\partial c],
\]
so that
\[
 H_p(\mathcal E_*)=\ker d_{\mathcal E,p}/\operatorname{im}d_{\mathcal E,p+1}.
\]
In particular, a nonzero defect space $\mathcal E_p$ does not by itself imply that $H_p(\mathcal E_*)$ is nonzero. If the displayed quotient spaces are used to compute defect homology, their differentials must be transported from $d_{\mathcal E}$ through the isomorphisms of Proposition~\ref{prop:illegal-face-defect}.

\begin{theorem}\label{thm:completion-shift}
Under $(\mathrm{P})$, there is a natural chain isomorphism
\begin{equation*}
 \mathcal E_*[-1]\xrightarrow{\cong}\mathcal T_*,
\end{equation*}
where $(\mathcal E[-1])_p=\mathcal E_{p+1}$ and $d_{\mathcal E[-1]}=-\partial_{\mathcal E}$.
Consequently, there is a long exact sequence
\begin{equation*}
 \cdots\to H_{r+2}(\mathcal E_*)\to H_r(Z)\to H_r(K_*)
 \to H_{r+1}(\mathcal E_*)\to H_{r-1}(Z)\to\cdots.
\end{equation*}
In particular, condition $(\mathrm{O})$ holds if and only if $\mathcal T_*=0$.
\end{theorem}
\begin{proof}
For $c\in\Omega_{p+1}(X)$, write $c=a+b$ with $a\in\mathcal A_{p+1}(A)$ and $b\in\mathcal A_{p+1}(B)$, and define
\[
 \Phi_p([c]):=[\nu_p\partial a]\in K_p/\Sigma_p(Z).
\]
Put $\eta=\nu_p\partial a=-\nu_p\partial b$. Subtracting the allowed part from $\partial a$ shows that $\eta\in\Sigma_p(A)$. The same argument applies to the other side, so $\eta\in K_p$. Changing the decomposition changes the result by $\nu_p\partial z$, where $z\in\mathcal A_{p+1}(Z)$, and this belongs to $\Sigma_p(Z)$. Adding a sum of invariant chains from the two sides does not change the image. The map is therefore well-defined.

For $k\in K_p$, write $k=u+\partial a=w+\partial b$, where $u\in\mathcal A_p(A)$ and $w\in\mathcal A_p(B)$. Then we have $a-b\in\Omega_{p+1}(X)$ and $\Phi_p([a-b])=[\nu_pk]=[k]$. The last equality holds because the allowed part of $k$ is supported in the intersection. This proves surjectivity. If $\eta\in\Sigma_p(Z)$, its purely non-allowed support allows us to choose $z\in\mathcal A_{p+1}(Z)$ with $\nu_p\partial z=\eta$. Then we obtain $a-z\in\Omega_{p+1}(A)$ and $b+z\in\Omega_{p+1}(B)$, proving injectivity.

It remains to verify that $\Phi$ commutes with the differentials. The assertion in degree zero follows from the convention for negative degrees. Let $p\ge1$. Write $\partial a=\alpha+\eta$ and $\partial b=\beta-\eta$, where $\alpha,\beta$ are allowed chains on the respective sides. Then we have $\partial c=\alpha+\beta$ and $\partial\eta=-\partial\alpha$. Modulo $\Sigma_{p-1}(Z)$, the allowed part disappears, giving $\partial\Phi_p([c])=-\Phi_{p-1}(\partial[c])$. Thus, $\Phi$ commutes with the shifted differential $d_{\mathcal E[-1]}=-\partial_{\mathcal E}$. Taking homology of $0\to\Sigma_*(Z)\to K_*\to\mathcal T_*\to0$ and applying Theorem~\ref{thm:sigma-quasi} gives the long exact sequence. The chain isomorphism gives the last assertion since $(\mathrm{O})$ is equivalent to $\mathcal E_*=0$.
\end{proof}

\subsection{Low-dimensional criteria}

\begin{example}\label{ex:mv-completed-diamond}
Consider the diamond cover in Example~\ref{ex:mv-independent}. Although $a\to b$ is not an edge, we have
\[
 e_{ab}=e_{xb}+e_{ax}-\partial e_{axb}\in\Sigma_1(A).
\]
Similarly, one has $e_{ab}\in\Sigma_1(B)$. Thus, we obtain
\[
 K_0=\Bbbk\{e_a,e_b\},\quad K_1=\Bbbk e_{ab},\quad
 \partial e_{ab}=e_b-e_a,
\]
and $K_p=0$ for $p\ge2$. Hence, $K_*$ is isomorphic to the ordinary chain complex of a single edge, while $\Sigma_*(Z)$ is concentrated in degree zero with basis $e_a,e_b$. Consequently, we have $\mathcal T_1=\Bbbk[e_{ab}]$ and $\mathcal T_p=0$ for $p\ne1$. Since $\nu_1\partial e_{axb}=-e_{ab}$, the isomorphism of Theorem~\ref{thm:completion-shift} maps $[e_{axb}-e_{ayb}]\in\mathcal E_2$ to $-[e_{ab}]\in\mathcal T_1$.
\end{example}

\begin{example}\label{ex:mv-acyclic-defect}
Add to the preceding example a common vertex $t$ and edges $t\to a,t\to x,t\to y,t\to b$, and include $t$ in both $U,W$. The edge $a\to b$ is not included. Then the edges of $Z$ are $t\to a$ and $t\to b$. The only non-allowed paths that occur in $\Sigma_*(A)\cap\Sigma_*(B)$ are $ab$ and $tab$. Indeed,
\[
 e_{tab}=\partial e_{taxb}-e_{axb}+e_{txb}+e_{tax}\in\Sigma_2(A),
\]
and replacing $x$ by $y$ gives the corresponding expression in $\Sigma_2(B)$. Neither $ab$ nor $tab$ belongs to $\Sigma_*(Z)$. Hence, the only nonzero terms are
\[
 \mathcal T_2=\Bbbk[e_{tab}],\quad
 \mathcal T_1=\Bbbk[e_{ab}],\quad
 \partial[e_{tab}]=[e_{ab}-e_{tb}+e_{ta}]=[e_{ab}].
\]
Thus, we have $\mathcal T_*\ne0$, while $H_q(\mathcal T_*)=0$ for each $q\ge0$. By Theorem~\ref{thm:completion-shift}, condition $(\mathrm{O})$ fails, although ordinary Mayer--Vietoris and homological excision hold. Consequently, the condition $H_\ast(\mathcal E_*)=0$ does not imply $\mathcal E_*=0$.
\end{example}

Let $X=P(G)$. For $D\subseteq V(G)$, let $M(D)$ be the set of ordered pairs of distinct vertices $a,b\in D$ such that $a\not\to b$ and $a\to x\to b$ for some $x\in D$.

\begin{proposition}\label{prop:mv-missing-edges}
Under $(\mathrm{P})$, we have $\mathcal T_0=0$, and a basis of $\mathcal T_1$ is
\[
 \bigl\{[e_{ab}]:(a,b)\in (M(U)\cap M(W))\setminus M(U\cap W)\bigr\}.
\]
\end{proposition}
\begin{proof}
In degree zero, $K_0=\Sigma_0(A)\cap\Sigma_0(B)$ is spanned by the vertices in $U\cap W$. Hence, $K_0=\Sigma_0(Z)$ and $\mathcal T_0=0$.

For $Y=A,B,Z$, we have $\Sigma_1(Y)=\mathcal A_1(Y)+\partial\mathcal A_2(Y)$. Each allowed path $a\to x\to b$ gives
\[
 \partial e_{axb}=e_{xb}-e_{ab}+e_{ax},\quad e_{aa}=0.
\]
Since $e_{ax}$ and $e_{xb}$ are allowed, this boundary adds $e_{ab}$ to the edge span exactly when $a\ne b$ and $a\not\to b$. Thus, the basis of $\Sigma_1(Y)$ consists of the edges of $Y$ and the missing edges that occur as endpoints of two-step paths in $Y$.

The common edges of $G[U]$ and $G[W]$ lie in $G[U\cap W]$, so their classes vanish in $\mathcal T_1=K_1/\Sigma_1(Z)$. The remaining basis elements of $K_1$ are indexed by $M(U)\cap M(W)$. Among these, exactly those indexed by $M(U\cap W)$ lie in $\Sigma_1(Z)$.
\end{proof}

\begin{definition}\label{def:two-step-compatible}
Let $G=(V,E)$ be a digraph, and let $U,W\subseteq V$ with $U\cup W=V$. The pair $(G[U],G[W])$ is \textit{two-step compatible} if
\[
 M(U)\cap M(W)=M(U\cap W).
\]
\end{definition}

\begin{corollary}\label{cor:mv-low-degree}
Assume that condition $(\mathrm{P})$ holds. If $H_1(\mathcal T_*)=0$, then $H_1(A,Z)\to H_1(X,B)$ is an isomorphism, and there is an exact sequence
\begin{equation*}
 \begin{aligned}
 H_1(Z)&\to H_1(A)\oplus H_1(B)\to H_1(X)\to H_0(Z)\\
       &\to H_0(A)\oplus H_0(B)\to H_0(X)\to0.
 \end{aligned}
\end{equation*}
In particular, when $X=P(G)$, these conclusions hold if the induced cover $(G[U],G[W])$ is two-step compatible and satisfies $(\mathrm{P})$.
\end{corollary}
\begin{proof}
In degrees zero and one, all allowed chains are $\partial$-invariant. Condition $(\mathrm{P})$ therefore gives $\mathcal E_0=\mathcal E_1=0$. By the chain isomorphism $\mathcal E_*[-1]\cong\mathcal T_*$ of Theorem~\ref{thm:completion-shift}, the condition $H_1(\mathcal T_*)=0$ is equivalent to $H_2(\mathcal E_*)=0$. Thus, $H_1(M_*)\to H_1(X)$ and $H_0(M_*)\to H_0(X)$ are isomorphisms. The long exact sequence for the sum complex in Theorem~\ref{thm:mv-defect} gives the stated Mayer--Vietoris segment, and Theorem~\ref{thm:mv-excision} gives the degree-one excision isomorphism because $H_2(\mathcal E_*)=H_1(\mathcal E_*)=0$. For a two-step compatible induced cover, Proposition~\ref{prop:mv-missing-edges} gives $\mathcal T_1=0$, and hence $H_1(\mathcal T_*)=0$.
\end{proof}

\begin{table}[htbp]
\centering
\caption{Mayer--Vietoris sequences for an induced cover. Each row gives three consecutive terms in degree $r$.}
\label{tab:mv-conditions}
\begin{tabular}{lccccc}
\hline
Hypothesis & Intersection & & Summands for $A$ and $B$ & & Union\\
\hline
Any induced cover & $H_r(Z)$ & $\to$ & $H_r(A)\oplus H_r(B)$ & $\to$ & $H_r(M_*)$\\
Condition $(\mathrm{O})$ & $H_r(Z)$ & $\to$ & $H_r(A)\oplus H_r(B)$ & $\to$ & $H_r(X)$\\
Condition $(\mathrm{P})$ & $H_r(K_*)$ & $\to$ & $H_r(A)\oplus H_r(B)$ & $\to$ & $H_r(X)$\\
\hline
\end{tabular}
\end{table}

For an induced digraph cover satisfying $(\mathrm{P})$, two-step compatibility is equivalent to $\mathcal T_1=0$ by Proposition~\ref{prop:mv-missing-edges}. Condition $(\mathrm{O})$, however, is equivalent to $\mathcal T_p=0$ for each $p\ge0$. Moreover, since $\mathcal T_0=0$, the weaker condition $H_1(\mathcal T_*)=0$ is equivalent to $\partial\mathcal T_2=\mathcal T_1$. This equality holds with $\mathcal T_1\ne0$ in Example~\ref{ex:mv-acyclic-defect}.

\subsection{The local Mayer--Vietoris sequence}

\begin{proposition}\label{prop:mv-local}
Use the notation above. Let $v\in U\cap W$. Assume that condition $(\mathrm{O})$ holds for both the induced cover $(A,B)$ of $X$ and the induced cover
$(A-v,B-v)$ of $X-v$. Then we have a short exact sequence
\[
 0\to Q_*(Z,v)\to Q_*(A,v)\oplus Q_*(B,v)\to Q_*(X,v)\to0,
\]
which induces a long exact sequence
\begin{equation*}
 \cdots\to\LH_r(Z,v)\to\LH_r(A,v)\oplus\LH_r(B,v)
 \to\LH_r(X,v)\to\LH_{r-1}(Z,v)\to\cdots.
\end{equation*}
\end{proposition}
\begin{proof}
Apply the first short exact sequence of Theorem~\ref{thm:mv-defect} to $X-v$ and $X$. By the two cover hypotheses, the addition maps are onto $\Omega_*(X-v)$ and $\Omega_*(X)$, respectively. Inclusion gives a map from the first sequence to the second whose three components are injective. Their cokernels are $Q_*(Z,v)$, $Q_*(A,v)\oplus Q_*(B,v)$, and $Q_*(X,v)$, respectively. The snake lemma gives a short exact sequence of these cokernels, with maps $[c]\mapsto([c],[-c])$ and $([a],[b])\mapsto[a+b]$. The short exact sequence leads to the local Mayer--Vietoris sequence.
\end{proof}

\begin{corollary}\label{cor:mv-local-convex}
Assume that condition $(\mathrm{P})$ holds and that $U$ is convex with respect to the allowed paths of $X$. Then, for any $v\in U\cap W$, there is a natural long exact sequence
\[
 \cdots\to\LH_r(Z,v)\to\LH_r(A,v)\oplus\LH_r(B,v)
 \to\LH_r(X,v)\to\LH_{r-1}(Z,v)\to\cdots.
\]
\end{corollary}
\begin{proof}
Let $T\subseteq V(X)$. An allowed path avoiding $T$ lies in $A-T$ or $B-T$ by $(\mathrm{P})$. If its endpoints belong to $U\setminus T$, convexity puts each vertex in $U$, hence in $U\setminus T$. Proposition~\ref{prop:mv-convex-cover} applied to $X-T$ gives
\[
 \Omega_*(X-T)=\Omega_*(A-T)+\Omega_*(B-T).
\]
The choices $T=\varnothing$ and $T=\{v\}$ give the two hypotheses of Proposition~\ref{prop:mv-local}.
\end{proof}

\begin{remark}\label{rem:mv-delete}
Condition $(\mathrm{O})$ in Proposition~\ref{prop:mv-local} may fail after deleting $v$. Let $G$ be the digraph consisting of the three paths $a\to v\to b$, $a\to x\to b$, and $a\to y\to b$. Set $X=P(G)$, $U=\{a,v,x,b\}$, and $W=\{a,v,y,b\}$. Then $\Omega_2(X)$ is spanned by $e_{avb}-e_{axb}$ and $e_{avb}-e_{ayb}$. One may check that $(\mathrm{O})$ holds in $X$. After deleting $v$ and the associated edges from $G$, we have $\Omega_2(X-v)=\Bbbk(e_{axb}-e_{ayb})$ and $\Omega_2(A-v)=\Omega_2(B-v)=0$. Hence, $(\mathrm{O})$ fails in degree two.
\end{remark}

\section{Stars and links}\label{sec:star}

In this section, we will introduce the star and the link of a path complex. Our star and link constructions recover the closed star and link of a simplicial complex. For a digraph, they give the path complexes of the star digraph and its induced link. The comparison maps below describe when reduced link homology computes local path homology.

\subsection{Star and link path complexes}

\begin{definition}\label{def:star-path}
Let $P$ be a regular path complex, and let $v\in V(P)$. Write $\beta\preceq\alpha$ when $\beta$ is obtained by deleting vertices of $\alpha$ without changing their order. Define the \textit{star path complex} and \textit{link path complex} at $v$ by
\[
 \begin{aligned}
 S_v(P)&:=\{\beta\in P:\beta\preceq\alpha\text{ for some }\alpha\in P\text{ containing }v\},\\
 L_v(P)&:=S_v(P)-v.
 \end{aligned}
\]
\end{definition}

One can verify that the star path complex and link path complex defined above are indeed path complexes.

\begin{proposition}\label{prop:star-link-complexes}
Let $P$ be a regular path complex, and let $v\in V(P)$. The star $S_v(P)$ and the link $L_v(P)$ are regular path complexes.
\end{proposition}
\begin{proof}
Let $\beta\in S_v(P)$ be a path of positive length. Choose an allowed path $\alpha$ through $v$ such that $\beta\preceq\alpha$. Deleting either endpoint of $\beta$ gives a path $\beta'\in P$ with $\beta'\preceq\alpha$, so one has $\beta'\in S_v(P)$. All paths of $S_v(P)$ are regular since they belong to $P$. Thus, $S_v(P)$ is a regular path complex. The link $L_v(P)=S_v(P)-v$ is an induced restriction of this regular path complex.
\end{proof}

\begin{definition}\label{def:weakly-connected-path}
A regular path complex $P$ is called \textit{weakly connected} if, for any $x,y\in V(P)$, there is a sequence $x=x_0,\ldots,x_m=y$ such that $x_{i-1}x_i\in P$ or $x_ix_{i-1}\in P$ for any $1\le i\le m$.
\end{definition}

\begin{proposition}\label{prop:star-graph}
Let $P$ be a regular path complex. For any vertex $v$ of $P$, the star $S_v(P)$ is weakly connected. Moreover, we have
\[
P=S_v(P)\cup(P-v),\quad L_v(P)=S_v(P)\cap(P-v)
\]
and
\[
S_v(S_v(P))=S_v(P),\quad L_v(S_v(P))=L_v(P).
\]
\end{proposition}

\begin{proof}
By definition, an allowed path $\alpha$ passing through $v$ belongs to $S_v(P)$, since $\alpha\preceq\alpha$. Every vertex of $S_v(P)$ lies on such a path, whose consecutive pairs are $1$-paths in $S_v(P)$. Following these edges in the underlying undirected graph connects that vertex to $v$. Any path of $P$ avoiding $v$ lies in $P-v$, so $P=S_v(P)\cup(P-v)$. Furthermore, the intersection satisfies $S_v(P)\cap(P-v)=S_v(P)-v=L_v(P)$.

Now, let $\beta\in S_v(P)$. Choose $\alpha\in P$ passing through $v$ such that $\beta\preceq\alpha$. Since $\alpha\in S_v(P)$, we obtain $\beta\in S_v(S_v(P))$. The reverse inclusion follows from the definition of the star. Deleting $v$ from both sides, we obtain the identity for the link.
\end{proof}

\begin{definition}\label{def:star-link}
Let $G=(V,E)$ be a finite digraph, and let $v\in V$. We define
\[
 \begin{gathered}
 R_v=\{x\in V:x\leadsto v\text{ or }v\leadsto x\},\\
 S_v=G[R_v],\quad L_v=G[R_v\setminus\{v\}]=S_v-v.
 \end{gathered}
\]
We call $S_v$ the \textit{star digraph} and $L_v$ the \textit{induced link digraph}.
\end{definition}

We write $x\leadsto y$ if there exists a directed path from $x$ to $y$. In particular, $v\in R_v$. Proposition~\ref{prop:star-graph} implies that $S_v$ is weakly connected. Moreover, the $R_v$ computed in $S_v$ coincides with that computed in $G$.

\begin{proposition}\label{prop:star-link-digraph}
For any digraph $G$ and vertex $v$ of $G$, we have
\[
 S_v(P(G))=P(S_v),\quad L_v(P(G))=P(L_v).
\]
\end{proposition}
\begin{proof}
Any vertex on a path passing through $v$ belongs to $R_v$, so $S_v(P(G))\subseteq P(S_v)$. Conversely, let $\beta=(x_0,\dots,x_p)\in P(S_v)$. If $\beta$ contains $v$, then $\beta\in S_v(P(G))$. If $v\leadsto x_0$ or $x_p\leadsto v$, we extend $\beta$ by concatenation to a path passing through $v$. In the remaining case, $x_0\leadsto v$ and $v\leadsto x_p$. Choose the largest index $i$ such that $x_i\leadsto v$. Then $i<p$ and $v\leadsto x_{i+1}$. We replace the edge $x_i\to x_{i+1}$ with a walk $x_i\leadsto v\leadsto x_{i+1}$. This yields a path passing through $v$ that contains $\beta$ as a subsequence. Common endpoints are identified only once, so concatenation preserves regularity. Hence $P(S_v)\subseteq S_v(P(G))$. Restricting to paths that avoid $v$ yields the identity for the link.
\end{proof}

\begin{definition}\label{def:monotone-path}
A path complex $P$ is called \textit{monotone} if its vertices can be assigned distinct real numbers such that every allowed path $x_0\cdots x_p$ satisfies $x_0<\cdots<x_p$.
\end{definition}

A path complex comes from a simplicial complex precisely when it is perfect and monotone \cite{GLMY}. To form the path complex of $K$, label its vertices by distinct real numbers and list each nonempty simplex in increasing order. Denote the resulting path complex by $P(K)$ and use the same labels on subcomplexes.

\begin{proposition}\label{prop:simplicial-star-link}
For any simplicial complex $K$ and vertex $v$, we have
\[
 S_v(P(K))=P(\overline{\st}_K(v)),\quad
 L_v(P(K))=P(\Lk_K(v)).
\]
\end{proposition}
\begin{proof}
A path $\beta\in P(K)$ belongs to $S_v(P(K))$ precisely when its vertices form a face of a simplex through $v$. Such faces form $\overline{\st}_K(v)$. If $\beta$ avoids $v$, its vertex set $\sigma$ has this property precisely when $\sigma\cup\{v\}\in K$. These are the simplices of $\Lk_K(v)$.
\end{proof}

\begin{example}\label{ex:simplicial-star}
Let $K$ be the boundary of the $2$-simplex with vertex labels $a<v<b$. Then
\[
 \begin{aligned}
 P(K)&=\{a,v,b,av,vb,ab\},\\
 S_v(P(K))&=\{a,v,b,av,vb\},\quad
 L_v(P(K))=\{a,b\}.
 \end{aligned}
\]
The star consists of the two edges incident to $v$, and the link consists of two vertices. The path $ab\in P(K)$ is absent from the star because no simplex through $v$ has $\{a,b\}$ as a face. Thus, the star need not be induced in $P$.
\end{example}

\begin{remark}\label{rem:cut-link}
Let $G=(V,E)$ be a directed acyclic graph and fix $v\in V$. To construct the cut link at $v$, put
\[
 A_v=\{x\ne v:x\leadsto v\},\quad
 B_v=\{y\ne v:v\leadsto y\}.
\]
A vertex in $A_v\cap B_v$ would lie on a directed cycle through $v$. Thus, $A_v\cap B_v=\varnothing$. The \textit{cut link} $L_v^{\mathrm{cut}}$ has vertex set $A_v^-\sqcup B_v^+$ and the following edges:
\[
 \begin{array}{lll}
 x^-\to y^- &\Longleftrightarrow& x,y\in A_v,\ (x,y)\in E,\\
 x^+\to y^+ &\Longleftrightarrow& x,y\in B_v,\ (x,y)\in E,\\
 x^-\to y^+ &\Longleftrightarrow& x\to v\text{ and }v\to y.
 \end{array}
\]

Identify $A_v^-$ with $A_v$ and $B_v^+$ with $B_v$. Within either set the cut and induced links have the same edges. Between the sets, the cut link has edge set $I_v\times O_v$, where $I_v=\{x:x\to v\}$ and $O_v=\{y:v\to y\}$. The induced link has edge set $E\cap(A_v\times B_v)$ between the sets. Thus, the cut and induced links agree exactly when $I_v\times O_v=E\cap(A_v\times B_v)$. For transitive $G$, each side equals $A_v\times B_v$.

Let $G$ be the directed tree with center $v$, incoming edges $a_i\to v$ for $1\le i\le p$, and outgoing edges $v\to b_j$ for $1\le j\le q$, where $p,q\ge1$. The induced link $L_v$ consists of $p+q$ isolated vertices. Since $\deg(v)=p+q$, by Proposition~\ref{prop:tree}, we have
\[
 \LH_r(G,v)\cong\widetilde H_{r-1}(L_v)\cong
 \begin{cases}
 \Bbbk^{p+q-1},&r=1,\\
 0,&r\ne1.
 \end{cases}
\]
The cut link $L_v^{\mathrm{cut}}$ has one edge $a_i\to b_j$ for each pair $(i,j)$. Its underlying graph is $K_{p,q}$, so the edge boundary has rank $p+q-1$. Since no edge leaves a vertex $b_j$, there are no allowed $2$-paths. Thus, we have
\[
 \widetilde H_0(L_v^{\mathrm{cut}})=0,\quad
 H_1(L_v^{\mathrm{cut}})\cong\Bbbk^{(p-1)(q-1)}.
\]
For all $p,q\ge1$, $\widetilde H_0(L_v^{\mathrm{cut}})=0$ while $\LH_1(G,v)\ne0$. In degree two, $\LH_2(G,v)=0$ and $\widetilde H_1(L_v^{\mathrm{cut}})\cong\Bbbk^{(p-1)(q-1)}$, so the cut-link formula fails exactly when $p\ge2$ and $q\ge2$.

Subsection~\ref{sec:comparisons} uses the inclusion $P(L_v)\subseteq P(S_v)=S_v(P(G))$ from Proposition~\ref{prop:star-link-digraph}. Suppose $x\to v\to y$ with $x\ne y$ and $x\not\to y$ in $G$. The cut link contains the edge $x^-\to y^+$, but $xy$ is not an allowed path of $S_v$. Thus, the vertex map $x^-\mapsto x$, $y^+\mapsto y$ does not induce an inclusion $P(L_v^{\mathrm{cut}})\hookrightarrow P(S_v)$.

\end{remark}

\subsection{Comparison with local homology}\label{sec:comparisons}

Fix a regular path complex $P$ and $v\in V(P)$. In this subsection, write $S=S_v(P)$, $L=L_v(P)$, $Q_*=Q_*(P,v)$, and $I_*=\Omega_*(S)/\Omega_*(L)$. Proposition~\ref{prop:star-graph} gives $P=S\cup(P-v)$ and $S\cap(P-v)=L$. Although $S$ need not be induced in $P$, the link $L=S-v$ is induced in $S$.

For the comparison maps and the completed link below, naturality is taken with respect to pointed path-complex morphisms $f:(P,v)\to(P',v')$ that induce maps of both pairs
\[
 (P,P-v)\longrightarrow(P',P'-v'),\quad
 (S_v(P),L_v(P))\longrightarrow(S_{v'}(P'),L_{v'}(P')).
\]
An inclusion preserving the marked vertex induces the displayed maps of pairs. For a digraph map with $f^{-1}(\{v'\})=\{v\}$, preservation of directed reachability gives $f(R_v)\subseteq R_{v'}$. Since only $v$ maps to $v'$, the restriction to $R_v$ also sends $R_v\setminus\{v\}$ into $R_{v'}\setminus\{v'\}$. Proposition~\ref{prop:star-link-digraph} gives the map of star--link pairs.

\begin{definition}\label{def:defect}
The \textit{star defect complex} at $v$ is defined by
\begin{equation*}
 D_*(P,v):=\frac{\Omega_*(P)}{\Omega_*(P-v)+\Omega_*(S)}.
\end{equation*}
For a digraph $G$, we write $D_*(G,v)=D_*(P(G),v)$.
\end{definition}

\begin{theorem}\label{thm:j}
There is a short exact sequence
\begin{equation*}
 0\longrightarrow I_*\xrightarrow{\ j\ }Q_*
 \longrightarrow D_*\longrightarrow0.
\end{equation*}
It induces a long exact sequence
\begin{equation*}
 \cdots\longrightarrow H_{r+1}(D_*)\longrightarrow H_r(S,L)
 \xrightarrow{\ j_r\ }\LH_r(P,v)
 \longrightarrow H_r(D_*)\longrightarrow H_{r-1}(S,L)\longrightarrow\cdots .
\end{equation*}
\end{theorem}
\begin{proof}
By Lemma~\ref{lem:intersection}, one has
\[
 \Omega_*(S)\cap\Omega_*(P-v)
 =\Omega_*(S)\cap\mathcal R_*(V(P)\setminus\{v\})
 =\Omega_*(L).
\]
Thus, the inclusion of pairs $(S,L)\hookrightarrow(P,P-v)$ induces an injective chain map $I_*\to Q_*$, whose cokernel is $D_*(P,v)$. Taking homology, we have the desired long exact sequence.
\end{proof}

For the inclusion $i:L\hookrightarrow S$, we have $(i_\#)_{-1}=\mathrm{id}_{\Bbbk}$. Thus, one has $I_p=\Omega_p(S)/\Omega_p(L),p\ge0$ and $I_{-1}=0$. The short exact sequence
\[
 0\longrightarrow\widetilde\Omega_*(L)\xrightarrow{i_\#}
 \widetilde\Omega_*(S)\xrightarrow{\pi} I_*\longrightarrow0
\]
induces the long exact sequence
\[
 \cdots\longrightarrow\widetilde H_r(L)\xrightarrow{i_r}\widetilde H_r(S)
 \xrightarrow{\pi_*} H_r(S,L)\xrightarrow{\delta_r}\widetilde H_{r-1}(L)
 \xrightarrow{i_{r-1}}\widetilde H_{r-1}(S)\longrightarrow\cdots .
\]
For a relative cycle represented by $c\in\Omega_r(S)$, the connecting map is given by
\[
 \delta_r[c]=
 \begin{cases}
 [\partial c],&r\ge1,\\
 [\epsilon(c)],&r=0,
 \end{cases}
\]
where $\epsilon$ is the augmentation.

The two comparison maps have common domain $H_r(S,L)$ and take the form
\begin{equation*}
 \LH_r(P,v)\xleftarrow{\ j_r\ }H_r(S,L)
 \xrightarrow{\ \delta_r\ }\widetilde H_{r-1}(L).
\end{equation*}

\begin{corollary}\label{cor:iso-conditions}
Fix $r\ge0$. Each hypothesis below suffices for the indicated map to be an isomorphism:
\begin{align*}
 H_{r+1}(D_*)=H_r(D_*)=0 &\quad\Longrightarrow\quad j_r\text{ is an isomorphism},\\
 \widetilde H_r(S)=\widetilde H_{r-1}(S)=0
 &\quad\Longrightarrow\quad\delta_r\text{ is an isomorphism}.
\end{align*}
If both hypotheses hold, then $\delta_rj_r^{-1}$ gives a natural isomorphism $\LH_r(P,v)\cong\widetilde H_{r-1}(L)$.
\end{corollary}

\begin{proposition}\label{prop:convex}
Let $U=V(S)$. Suppose that $S=P|_U$ and that $U$ is convex with respect to the allowed paths of $P$. Then $D_*=0$, and $j_r$ is an isomorphism for each $r\ge0$.
\end{proposition}
\begin{proof}
Under the stated hypothesis, $S$ and $P-v$ are induced in $P$. Their cover satisfies $(\mathrm{P})$ by Proposition~\ref{prop:star-graph}. By Proposition~\ref{prop:mv-convex-cover}, one has $\Omega_*(P)=\Omega_*(S)+\Omega_*(P-v)$. It follows that $D_*=0$. Theorem~\ref{thm:j} then gives the conclusion.
\end{proof}

\begin{corollary}\label{cor:convex-star-digraph}
Let $G$ be a digraph and let $v\in V(G)$. Suppose that $R_v=V(S_v)$ is convex with respect to the allowed paths of $P(G)$. Then $D_*(G,v)=0$, and the inclusion of pairs $(S_v,L_v)\hookrightarrow(G,G-v)$ induces an isomorphism
\[
 j_r:H_r(S_v,L_v)\xrightarrow{\cong}\LH_r(G,v)
 \quad\text{for each }r\ge0.
\]
\end{corollary}
\begin{proof}
Put $P=P(G)$ and $U=R_v$. By Proposition~\ref{prop:star-link-digraph}, we have
\[
 S_v(P(G))=P(S_v)=P(G[U])=P|_U.
\]
The hypothesis states that $U$ is convex with respect to the allowed paths of $P$. Thus, Proposition~\ref{prop:convex} applies.
\end{proof}

\begin{corollary}\label{cor:perfect-star-local}
If $P$ is perfect, then $D_*(P,v)=0$ and $j_r$ is an isomorphism for each $r\ge0$. If $P$ is also monotone, then
\[
 \LH_r(P,v)\cong\widetilde H_{r-1}(L_v(P))
 \quad\text{for each }r\ge0.
\]
\end{corollary}
\begin{proof}
Every subsequence of an allowed path of a perfect complex is allowed. Hence, $S$ and $P-v$ are perfect, and $\Omega_*(Y)=\mathcal A_*(Y)$ for $Y=P,S,P-v$. The cover $P=S\cup(P-v)$ gives $D_*(P,v)=0$. If $P$ is also monotone, write $P=P(K)$ for a simplicial complex $K$. By Proposition~\ref{prop:simplicial-star-link}, $S$ corresponds to the closed star of $v$ in $K$, a cone with apex $v$. Thus, we have $\widetilde H_*(S)=0$, and Corollary~\ref{cor:iso-conditions} gives the formula.
\end{proof}

\begin{proposition}\label{prop:low-degree}
We have $D_0=D_1=0$. Thus, $j_0$ is always an isomorphism and $j_1$ is always surjective. Since $S$ is weakly connected, $\delta_1$ is also always surjective.
\end{proposition}
\begin{proof}
The spaces $\Omega_0(P)$ and $\Omega_1(P)$ are generated by allowed $0$-paths and $1$-paths. Each such path either avoids $v$ or belongs to $S$, proving $D_0=D_1=0$. In Theorem~\ref{thm:j}, the groups on either side of $j_0$ are $H_1(D_*)=H_0(D_*)=0$, and the group following $j_1$ is $H_1(D_*)=0$. Exactness gives the isomorphism in degree zero and surjectivity in degree one. By Proposition~\ref{prop:star-graph}, $S$ is weakly connected, so $\widetilde H_0(S)=0$. Exactness at $\widetilde H_0(L)$ in the sequence of the pair $(S,L)$ gives surjectivity of $\delta_1$.
\end{proof}

Since $j_1$ and $\delta_1$ are surjective, the assignment $j_1(a)\mapsto\delta_1(a)$ defines an isomorphism $\LH_1(P,v)\to\widetilde H_0(L)$ if and only if $\ker j_1=\ker\delta_1$.

\subsection{Failures of the comparison maps}\label{sec:examples}

In the following examples, we take $P=P(G)$ and use $S_v(P(G))=P(S_v)$ and $L_v(P(G))=P(L_v)$. Thus, the maps $j_r$ and $\delta_r$ are the comparison maps of Subsection~\ref{sec:comparisons} for the star and link path complexes. Each edge set consists precisely of the listed edges. The computations are valid over any field. The first two examples illustrate the failure of injectivity for $j_1$ and $\delta_1$. The final two demonstrate failure of surjectivity in degree two.

\begin{example}\label{ex:diamond}
Let $G$ be the digraph on four vertices with edges
\[
 a\to v,\quad v\to b,\quad a\to w,\quad w\to b.
\]
Neither $v$ nor $w$ is reachable from the other. Hence, $S_v$ is the path $a\to v\to b$, $L_v=\{a,b\}$, and $H_1(S_v,L_v)=\Bbbk$. However, we have
\[
 z=e_{avb}-e_{awb}\in\Omega_2(G),\quad
 \partial z=e_{vb}+e_{av}-e_{wb}-e_{aw}.
\]
In the local quotient, $\bar\partial[z]=[e_{av}+e_{vb}]$, so $\LH_1(G,v)=0$. Therefore, $j_1:\Bbbk\to0$ is not injective. Here $D_2=\Bbbk[z]$ and $D_p=0$ for $p\ne2$. The connecting homomorphism $H_2(D_*)\to H_1(S_v,L_v)$ sends $[z]$ to $[e_{av}+e_{vb}]$, which generates $\ker j_1$.

Since $B_1(v)=S_v$ and $B_2(v)=G$, the degree-one comparison map is $\Bbbk\to0$. Removing edges $a\to w$ and $w\to b$ leaves $S_v$ and $L_v$ unchanged, while the resulting digraph $G'$ satisfies $\LH_1(G',v)\cong\Bbbk$. This shows that the path homology of $L_v$ does not determine $\LH_*(G,v)$ for general digraphs.
\end{example}

\begin{figure}[htbp]
\centering
\begin{tikzpicture}[>=Stealth, every node/.style={circle,draw,inner sep=2.5pt},scale=1]
\node (a) at (0,0) {$a$};
\node (v) at (1.8,0.9) {$v$};
\node (w) at (1.8,-0.9) {$w$};
\node (b) at (3.6,0) {$b$};
\draw[->] (a)--(v); \draw[->] (v)--(b);
\draw[->] (a)--(w); \draw[->] (w)--(b);
\end{tikzpicture}
\caption{The diamond: $w$ lies outside the star digraph and supplies the term needed to cancel the non-allowed edge $ab$.}
\label{fig:diamond}
\end{figure}
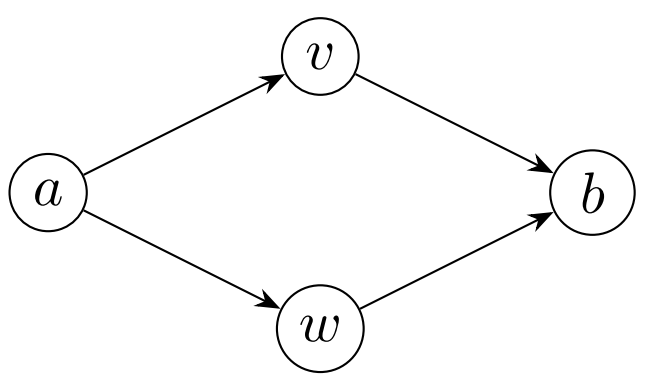

\begin{example}\label{ex:nonacyclic-star}
Consider the directed acyclic graph
\[
 v\to a,\quad a\to b,\quad b\to t,\quad v\to t.
\]
Every vertex is reachable from $v$, so $S_v=G$. The only allowed $2$-paths are $vab,abt$, with distinct non-allowed boundary edges $vb,at$, so $\Omega_2(G)=0$. Thus, the cycle
\[
 c=e_{va}+e_{ab}+e_{bt}-e_{vt}
\]
generates $H_1(S_v)=\Bbbk$. The link $L_v$ is $a\to b\to t$, so $\widetilde H_q(L_v)=0$ for each $q\ge-1$. Since $S_v=G$, the map $j_1$ is an isomorphism, whereas $\delta_1:\Bbbk\to\widetilde H_0(L_v)=0$ is not injective.
\end{example}

\begin{example}\label{ex:j-not-surjective}
Let $G$ be the following three-layer directed acyclic graph, with every arrow between consecutive layers included:
\[
 \{a,b\}\longrightarrow\{v,w\}\longrightarrow\{c,d\}.
\]
Then we have $w\notin S_v$. In $S_v$ each endpoint pair admits only one length-two path, so the non-allowed endpoint edge cannot cancel. Thus, we obtain $\Omega_2(S_v)=0$ and $H_2(S_v,L_v)=0$. Put
\begin{align*}
 u_v&=e_{avc}-e_{bvc}-e_{avd}+e_{bvd},\\
 u_w&=e_{awc}-e_{bwc}-e_{awd}+e_{bwd}.
\end{align*}
Then $z=u_v-u_w\in\Omega_2(G)$ and $\partial z=0$. The space $Q_2(G,v)$ has basis $q_{ij}=[e_{ivj}-e_{iwj}]$, where $i\in\{a,b\}$ and $j\in\{c,d\}$, and
\[
 \partial q_{ij}=[e_{iv}+e_{vj}].
\]
Thus, $\sum_{i,j}\lambda_{ij}q_{ij}$ is a cycle exactly when $\sum_j\lambda_{ij}=0$ for each $i$ and $\sum_i\lambda_{ij}=0$ for each $j$. The solution space is generated by the matrix $\left(\begin{smallmatrix}1&-1\\-1&1\end{smallmatrix}\right)$, corresponding to $[z]$. Since $\mathcal A_3(G)=0$, we have $Q_3(G,v)=0$. Thus, we obtain $\LH_2(G,v)=\Bbbk$ and $j_2:0\to\Bbbk$ is not surjective.
\end{example}

\begin{example}\label{ex:delta-not-surjective}
Let $L$ have the edges $u\to a\to t$ and $u\to b\to c\to t$. Add a vertex $v$ and the single edge $v\to u$, obtaining $G$. Then we have $S_v=G$ and $L_v=L$. The allowed $2$-paths are $vua,vub,uat,ubc,bct$, with distinct non-allowed boundary edges $va,vb,ut,uc,bt$. Thus, we obtain $\Omega_2(G)=0$ and $H_2(S_v,L_v)=0$. The underlying graph of $L$ has one cycle and $\Omega_2(L)=0$, so $H_1(L)=\Bbbk$. Consequently, $\delta_2:0\to\Bbbk$ is not surjective.
\end{example}

\subsection{The completed link}

Let $P$ be a regular path complex and let $v\in V(P)$. Write $S=S_v(P)$ and $L=L_v(P)$. We use $D_*(P,v)$ from Definition~\ref{def:defect}.

\begin{definition}\label{def:completed-link}
The \textit{completed link chain complex} is defined by
\[
 \widehat{\mathcal L}_*(P,v):=\Sigma_*(S)\cap\Sigma_*(P-v).
\]
For a digraph $G$, we write $\widehat{\mathcal L}_*(G,v)=\widehat{\mathcal L}_*(P(G),v)$. Since boundaries introduce no new vertices, we have $\widehat{\mathcal L}_p(P,v)\subseteq\mathcal R_p(V(L))$ for any $p\ge0$.
\end{definition}

\begin{theorem}\label{thm:completed-local}
For any regular path complex $P$ and vertex $v$, there exists a natural isomorphism
\[
 \LH_r(P,v)\cong
 H_r\!\left(\Sigma_*(S)/\widehat{\mathcal L}_*(P,v)\right)
 \quad\text{for each }r\ge0,
\]
and a natural chain isomorphism
\[
 \widehat{\mathcal L}_*(P,v)/\Sigma_*(L)\cong D_*(P,v)[-1].
\]
Under the chosen augmentation convention, we also have the long exact sequence
\[
 \cdots\to\widetilde H_r(S)\to\LH_r(P,v)
 \to\widetilde H_{r-1}(\widehat{\mathcal L}_*)
 \to\widetilde H_{r-1}(S)\to\cdots.
\]
If $\widetilde H_r(S)=\widetilde H_{r-1}(S)=0$, the connecting homomorphism yields a natural isomorphism $\LH_r(P,v)\cong\widetilde H_{r-1}(\widehat{\mathcal L}_*(P,v))$.
\end{theorem}
\begin{proof}
Set $B=P-v$. By Proposition~\ref{prop:star-graph}, every allowed path of $P$ lies in $S$ or $B$. Thus, we have
\[
 \mathcal A_*(P)=\mathcal A_*(S)+\mathcal A_*(B),\quad
 \Sigma_*(P)=\Sigma_*(S)+\Sigma_*(B).
\]
The isomorphism theorem gives
\[
 \Sigma_*(S)/\widehat{\mathcal L}_*(P,v)
 \cong\Sigma_*(P)/\Sigma_*(B).
\]
Taking homology and applying Theorem~\ref{thm:sigma-quasi} establishes the first claim.

Let $\pi_p$ denote the projection onto $\mathcal A_p(P)$ along the span of non-allowed regular paths, and define $\nu_p=1-\pi_p$. Let $\alpha\in S$ and $\beta\in P$ satisfy $\beta\preceq\alpha$. By Definition~\ref{def:star-path}, there exists $\gamma\in P$ containing $v$ such that $\alpha\preceq\gamma$. It follows that $\beta\preceq\alpha\preceq\gamma$, so $\beta\in S$. Any subsequence of a path avoiding $v$ also avoids $v$, so the same closure property holds for $B$ and $L$. Consequently, we have
\[
 \pi_p\partial\mathcal A_{p+1}(Y)\subseteq\mathcal A_p(Y),
 \quad \pi_p\Sigma_p(Y)=\mathcal A_p(Y),
 \quad Y\in\{S,B,L\}.
\]
Since $\mathcal A_*(S)\cap\mathcal A_*(B)=\mathcal A_*(L)$, the projection $\pi_p$ maps $\widehat{\mathcal L}_p$ into $\mathcal A_p(L)$. For $c\in\Omega_{p+1}(P)$, write $c=a+b$ with $a\in\mathcal A_{p+1}(S)$ and $b\in\mathcal A_{p+1}(B)$, and define
\[
 \Phi_p:D_{p+1}(P,v)\longrightarrow
 \widehat{\mathcal L}_p/\Sigma_p(L),\quad
 [c]\longmapsto[\nu_p\partial a].
\]
The identity $\nu_p\partial a=-\nu_p\partial b$ ensures this chain belongs to both completed complexes. A different decomposition alters the representative by $\nu_p\partial z\in\Sigma_p(L)$ for some $z\in\mathcal A_{p+1}(L)$. Any element of $\Omega_{p+1}(S)+\Omega_{p+1}(B)$ maps to zero. Hence, $\Phi_p$ is well-defined.

Take $k\in\widehat{\mathcal L}_p$, and write $k=u+\partial a=w+\partial b$, where $u\in\mathcal A_p(S)$, $a\in\mathcal A_{p+1}(S)$, $w\in\mathcal A_p(B)$, and $b\in\mathcal A_{p+1}(B)$. Then $a-b\in\Omega_{p+1}(P)$ and $\Phi_p([a-b])=[\nu_pk]=[k]$, which proves surjectivity. If $\Phi_p([c])=0$, choose $z\in\mathcal A_{p+1}(L)$ such that $\nu_p\partial z=\nu_p\partial a$. The projection identities imply $a-z\in\Omega_{p+1}(S)$ and $b+z\in\Omega_{p+1}(B)$, so $[c]=0$ in $D_{p+1}$. Next, for $p\ge1$, we have
\[
 \partial\Phi_p([c])
 =[-\partial(\pi_p\partial a)]
 =[-\nu_{p-1}\partial(\pi_p\partial a)]
 =-\Phi_{p-1}([\partial c]).
\]
The middle equality holds because $\pi_{p-1}(\widehat{\mathcal L}_{p-1})\subseteq\mathcal A_{p-1}(L)$.
Let $f$ be a morphism as in Subsection~\ref{sec:comparisons}, and put $L'=L_{v'}(P')$ and $\eta=\nu_p\partial a$. Since $f_\#$ preserves allowed chains and $f_\#\eta\in\widehat{\mathcal L}_p(P',v')$, we have
\[
 f_\#\eta-\nu'_p\partial(f_\#a)=\pi'_p f_\#\eta\in\mathcal A_p(L').
\]
Taking classes modulo $\Sigma_p(L')$ gives $f_\#\Phi_p([c])=\Phi'_p([f_\#c])$.

After augmentation, the inclusion $\widehat{\mathcal L}_*\hookrightarrow\Sigma_*(S)$ is the identity on $\Bbbk$ in degree $-1$, so the quotient is zero in that degree. The resulting short exact sequence with quotient $\Sigma_*(S)/\widehat{\mathcal L}_*$ gives the stated reduced homology sequence, even when $L=\varnothing$.
\end{proof}

\begin{proposition}\label{prop:completed-link-comparison}
There exists a short exact sequence of chain complexes
\[
 0\longrightarrow\Sigma_*(L)\longrightarrow\widehat{\mathcal L}_*(P,v)
 \longrightarrow D_*(P,v)[-1]\longrightarrow0.
\]
Under the augmentation convention, it induces the exact sequence
\[
 \begin{aligned}
 \cdots\longrightarrow H_{q+2}(D_*(P,v))
 &\longrightarrow\widetilde H_q(L)
 \longrightarrow\widetilde H_q(\widehat{\mathcal L}_*(P,v))\\
 &\longrightarrow H_{q+1}(D_*(P,v))
 \longrightarrow\widetilde H_{q-1}(L)\longrightarrow\cdots .
 \end{aligned}
\]
In particular, if $H_{q+2}(D_*(P,v))=H_{q+1}(D_*(P,v))=0$, the inclusion of augmented chain complexes $\widetilde\Omega_*(L)\hookrightarrow\widetilde{\widehat{\mathcal L}}_*(P,v)$ induces an isomorphism $\widetilde H_q(L)\cong\widetilde H_q(\widehat{\mathcal L}_*(P,v))$.
\end{proposition}
\begin{proof}
The short exact sequence follows from Theorem~\ref{thm:completed-local}. Augment $\Sigma_*(L)$ and $\widehat{\mathcal L}_*(P,v)$ by $\Bbbk$ in degree $-1$, using the identity map on these degree $-1$ terms. Using the identification $H_q(D_*[-1])\cong H_{q+1}(D_*)$ and Theorem~\ref{thm:sigma-quasi}, we have the desired result.
\end{proof}

\begin{corollary}\label{cor:perfect-completed-link}
If $P$ is a perfect path complex, then we have
\[
 \widehat{\mathcal L}_*(P,v)=\Sigma_*(L)=\Omega_*(L).
\]
\end{corollary}
\begin{proof}
By Corollary~\ref{cor:perfect-star-local}, $D_*(P,v)=0$. By Proposition~\ref{prop:completed-link-comparison}, we obtain $\widehat{\mathcal L}_*(P,v)=\Sigma_*(L)$. This completes the proof.
\end{proof}

\section{Products and radius bounds}\label{sec:kunneth-radius}

We compute the local homology and the radius filtration of a Cartesian product. Products of directed intervals attain the bounds of Section~\ref{sec:local}. We also give a radius bound for detecting local regularity when the chain groups vanish above a fixed degree.

\subsection{The local K\"unneth formula}

\begin{definition}\label{def:cartesian-product}
The \textit{Cartesian product} of digraphs $G$ and $H$ is the digraph $G\Box H$ on $V(G)\times V(H)$ with edges
\[
 (x,y)\to(x',y')\quad\Longleftrightarrow\quad
 \bigl(x=x'\text{ and }y\to y'\bigr)
 \text{ or }\bigl(y=y'\text{ and }x\to x'\bigr).
\]
\end{definition}

Let $\alpha=(x_0,\ldots,x_p)$ and $\beta=(y_0,\ldots,y_q)$ be regular paths. A \textit{shuffle} of $\alpha$ and $\beta$ is a path $\gamma=((x_{i_k},y_{j_k}))_{k=0}^{p+q}$, where $(i_k,j_k)$ runs from $(0,0)$ to $(p,q)$ by steps $(1,0)$ or $(0,1)$. We can recover $\alpha$ from $\gamma$ by projecting onto the first coordinate.

\begin{definition}\label{def:cartesian-product-path}
The \textit{Cartesian product} $P\Box Q$ of regular path complexes $P$ and $Q$ is the path complex on $V(P)\times V(Q)$ whose allowed paths are the shuffles of a path $\alpha\in P$ and a path $\beta\in Q$.
\end{definition}

For elementary paths, the signed shuffle product is
\[
 e_\alpha\times e_\beta=\sum_{\gamma}\operatorname{sgn}(\gamma)e_\gamma,
\]
where $\gamma$ runs over the shuffles of $\alpha$ and $\beta$. Here $\operatorname{sgn}(\gamma)=(-1)^{N(\gamma)}$, with $N(\gamma)$ the number of pairs of steps in which a step of $\beta$ precedes a step of $\alpha$. On $\partial$-invariant chains this product gives the natural chain isomorphism \cite[Theorem~6.6]{GLMYProduct}
\[
 \Phi:\Omega_*(P)\otimes_{\Bbbk}\Omega_*(Q)\xrightarrow{\cong}\Omega_*(P\Box Q),
 \quad u\otimes w\longmapsto u\times w.
\]
For $u\in\Omega_p(P)$ and $w\in\Omega_q(Q)$, we have
\[
 \partial(u\times w)=(\partial u)\times w+(-1)^p u\times\partial w.
\]

\begin{theorem}\label{thm:local-kunneth}
Let $P,Q$ be regular path complexes on finite vertex sets, and let $T\subseteq V(P)$ and $U\subseteq V(Q)$ be nonempty. The shuffle product induces a natural chain isomorphism
\[
 \overline\Phi:Q_*(P,T)\otimes_{\Bbbk}Q_*(Q,U)
 \xrightarrow{\cong}Q_*(P\Box Q,T\times U),
 \quad [u]\otimes[w]\longmapsto[u\times w].
\]
For each $n\ge0$, the induced product gives a homology isomorphism
\[
 \bigoplus_{i+j=n}\LH_i(P,T)\otimes_{\Bbbk}\LH_j(Q,U)
 \xrightarrow{\cong}\LH_n(P\Box Q,T\times U).
\]
\end{theorem}
\begin{proof}
Put $C=\Omega_*(P)$, $C'=\Omega_*(P-T)$, $D=\Omega_*(Q)$, and $D'=\Omega_*(Q-U)$. In the elementary tensor basis of $\mathcal A_*(P)\otimes\mathcal A_*(Q)$, write $z=\sum_{\alpha,\beta}z_{\alpha,\beta}e_\alpha\otimes e_\beta\in C\otimes D$. The coefficient of each shuffle $\gamma$ of $(\alpha,\beta)$ in $\Phi(z)$ is $\operatorname{sgn}(\gamma)z_{\alpha,\beta}$. A shuffle through $(\alpha_i,\beta_j)$ exists for each pair of positions $(i,j)$. Thus, $\Phi(z)$ is supported outside $T\times U$ if and only if
\[
 z_{\alpha,\beta}=0\quad\text{whenever $\alpha$ meets $T$ and $\beta$ meets $U$.}
\]

In the elementary path expansions of chains in $C$ and $D$, retain the terms meeting $T$ and $U$, respectively. Denote the projections by $\pi_T$ and $\pi_U$. The coefficient condition is $(\pi_T\otimes\pi_U)z=0$, and Lemma~\ref{lem:intersection} gives $\ker\pi_T=C'$ and $\ker\pi_U=D'$. Over $\Bbbk$ this gives
\[
 \Phi^{-1}\Omega_*((P\Box Q)-(T\times U))
 =\ker(\pi_T\otimes\pi_U)
 =C'\otimes D+C\otimes D'.
\]
Passing to the quotient gives $\overline\Phi$, since $(C\otimes D)/(C'\otimes D+C\otimes D')\cong(C/C')\otimes(D/D')$. The homology formula follows from the K\"unneth theorem over $\Bbbk$.
\end{proof}

\begin{corollary}\label{cor:point-local-kunneth}
For $v\in V(P)$, $w\in V(Q)$, and $n\ge0$, we have
\[
 \bigoplus_{i+j=n}\LH_i(P,v)\otimes_{\Bbbk}\LH_j(Q,w)
 \xrightarrow{\cong}\LH_n(P\Box Q,(v,w)).
\]
\end{corollary}
\begin{proof}
Take $T=\{v\}$ and $U=\{w\}$ in Theorem~\ref{thm:local-kunneth}.
\end{proof}

Let $\overline\Phi$ be the local shuffle map of Theorem~\ref{thm:local-kunneth} for $(P,T)$ and $(Q,U)$, and let $\overline\Phi'$ be the corresponding map for $(P',T')$ and $(Q',U')$.

\begin{proposition}\label{prop:local-kunneth-naturality}
Suppose the morphisms $f:P\to P'$ and $g:Q\to Q'$ satisfy $f^{-1}(T')\subseteq T$ and $g^{-1}(U')\subseteq U$. Then we have
\[
 (f\times g)_\#\circ\overline\Phi
 =\overline\Phi'\circ(f_\#\otimes g_\#).
\]
\end{proposition}
\begin{proof}
Since $(f\times g)^{-1}(T'\times U')\subseteq T\times U$, $f\times g$ induces a map from $Q_*(P\Box Q,T\times U)$ to $Q_*(P'\Box Q',T'\times U')$. The complement conditions on $f$ and $g$ give the corresponding maps on $Q_*(P,T)$ and $Q_*(Q,U)$. For any $a\in\Omega_i(P)$ and $b\in\Omega_j(Q)$, one has
\[
 \begin{aligned}
 (f\times g)_\#\overline\Phi([a]\otimes[b])
 &= [(f\times g)_\#(a\times b)]\\
 &= [f_\#a\times g_\#b]
 = \overline\Phi'(f_\#[a]\otimes g_\#[b]).
 \end{aligned}
\]
By linearity, the desired result follows.
\end{proof}

\begin{remark}
Set $P=Q=P(I)$ where $I=(0\to1)$, and take $T=\{0\}$, $U=\{1\}$. The allowed path
\[
 \gamma=((0,0),(1,0),(1,1))
\]
avoids $(0,1)$, even though both of its projections reduce to the edge $0\to1$.

The chain $e_\gamma\notin\Omega_2(P\Box Q)$, since its boundary contains $((0,0),(1,1))$, which is not an allowed edge. The second shuffle $\gamma'=((0,0),(0,1),(1,1))$ intersects $T\times U$, and
\[
 e_{01}\times e_{01}=e_\gamma-e_{\gamma'}\in\Omega_2(P\Box Q).
\]
The two shuffles carry opposite coefficients. In this example, imposing that the coefficient of the shuffle passing through $(0,1)$ be zero forces the coefficient of the other shuffle to vanish as well.
\end{remark}

\begin{corollary}\label{cor:strip-mv}
Let $X=(P\Box Q)-(T\times U)$, $A=(P-T)\Box Q$, and $B=P\Box(Q-U)$. Then we have
\[
 \Omega_*(X)=\Omega_*(A)+\Omega_*(B),\quad
 A\cap B=(P-T)\Box(Q-U).
\]
\end{corollary}
\begin{proof}
In the notation of Theorem~\ref{thm:local-kunneth}, the shuffle map gives
\[
 \Omega_*(X)=\Phi(C'\otimes D+C\otimes D')
 =\Omega_*(A)+\Omega_*(B).
\]
A shuffle belongs to both $A$ and $B$ exactly when its factor paths lie in $P-T$ and $Q-U$, respectively. Hence, $A\cap B=(P-T)\Box(Q-U)$.
\end{proof}

\begin{corollary}\label{cor:strip-mv-excision}
Let $X,A,B$ be as in Corollary~\ref{cor:strip-mv}, and put $Z=A\cap B$. There is a long exact sequence
\[
 \cdots\longrightarrow H_r(Z)\longrightarrow H_r(A)\oplus H_r(B)
 \longrightarrow H_r(X)\longrightarrow H_{r-1}(Z)\longrightarrow\cdots .
\]
The inclusion $(A,Z)\hookrightarrow(X,B)$ induces an isomorphism
\[
 H_r(A,Z)\xrightarrow{\cong}H_r(X,B),\quad r\ge0.
\]
\end{corollary}
\begin{proof}
By Corollary~\ref{cor:strip-mv} and Lemma~\ref{lem:intersection}, addition maps $\Omega_*(A)\oplus\Omega_*(B)$ onto $\Omega_*(X)$, and its kernel consists of the pairs $(z,-z)$ with $z\in\Omega_*(Z)$. We have the short exact sequence
\[
 0\longrightarrow\Omega_*(Z)\xrightarrow{z\mapsto(z,-z)}
 \Omega_*(A)\oplus\Omega_*(B)\xrightarrow{(a,b)\mapsto a+b}
 \Omega_*(X)\longrightarrow0,
\]
which leads to the long exact sequence. For the inclusion $(A,Z)\hookrightarrow(X,B)$, the chain isomorphism
\[
 \frac{\Omega_*(A)}{\Omega_*(Z)}
 =\frac{\Omega_*(A)}{\Omega_*(A)\cap\Omega_*(B)}
 \cong
 \frac{\Omega_*(A)+\Omega_*(B)}{\Omega_*(B)}
 =\frac{\Omega_*(X)}{\Omega_*(B)}
\]
induces the desired homology isomorphism.
\end{proof}

\subsection{Local regularity under products}\label{sec:regularity-products}

\begin{definition}\label{def:local-regular}
A vertex $v\in V(P)$ is \textit{locally homologically regular of dimension $n$} if $n\ge0$ and
\[
 \LH_r(P,v)\cong
 \begin{cases}\Bbbk,&r=n,\\0,&r\ne n.\end{cases}
\]
\end{definition}

The \textit{local Poincar\'e series} of $P$ at $v$ is defined as
\[
 \ell_{P,v}(t)=\sum_{r\ge0}\dim_{\Bbbk}\LH_r(P,v)\,t^r.
\]
By definition, $v$ is regular of dimension $n$ when $\ell_{P,v}(t)=t^n$.

\begin{corollary}\label{prop:product-regularity}
The vertex $(v,w)$ is regular of dimension $n$ if and only if $v$ and $w$ are regular and their dimensions sum to $n$.
\end{corollary}
\begin{proof}
By Theorem~\ref{thm:local-kunneth}, we have $\ell_{P\Box Q,(v,w)}(t)=\ell_{P,v}(t)\ell_{Q,w}(t)$. Since the two factors have nonnegative integer coefficients, their product equals $t^n$ if and only if they are $t^a$ and $t^b$ for some $a,b\ge0$ with $a+b=n$.
\end{proof}

\begin{corollary}\label{cor:tree-regular}
Let $T_1,\dots,T_m$ be directed trees. For each $i$, let $v_i$ be a vertex of $T_i$ of undirected degree $d_i\ge1$. Then
\[
 \LH_q\left(\mathop{\Box}_{i=1}^m T_i,(v_1,\ldots,v_m)\right)
 \cong\begin{cases}
 \Bbbk^{\prod_{i=1}^m(d_i-1)},&q=m,\\0,&q\ne m.
 \end{cases}
\]
The vertex $(v_1,\ldots,v_m)$ is regular of dimension $m$ if and only if each $d_i=2$.
\end{corollary}
\begin{proof}
Proposition~\ref{prop:tree} gives $\ell_{T_i,v_i}(t)=(d_i-1)t$. The product series is therefore $\bigl(\prod_{i=1}^m(d_i-1)\bigr)t^m$.
\end{proof}

At either endpoint $\varepsilon\in\{0,1\}$ of the directed interval $I$, the local chain complex is
\[
 0\longrightarrow\Bbbk[e_{01}]
 \xrightarrow{\;\partial\;}\Bbbk[e_\varepsilon]\longrightarrow0,
 \quad \partial[e_{01}]=(2\varepsilon-1)[e_\varepsilon].
\]
Since $\partial$ is an isomorphism, $\ell_{I,\varepsilon}(t)=0$. For each $v\in V(P)$, the local K\"unneth formula implies
\[
 \LH_q(P\Box P(I),(v,\varepsilon))=0,\quad q\ge0.
\]

\subsection{Products of radius filtrations}

For $s\ge0$, set
\[
 F_sQ_*(P,T)=\operatorname{im}\!\left(
 Q_*(P|_{N_s(T)},T)\longrightarrow Q_*(P,T)\right),\quad F_{-1}Q_*(P,T)=0.
\]
By Theorem~\ref{thm:chain-locality}, the inclusion map is injective and $F_sQ_p(P,T)=Q_p(P,T)$ for $s\ge p$.

\begin{theorem}\label{thm:filtered-kunneth}
Let $P,Q,T,U$ be as in Theorem~\ref{thm:local-kunneth}. For each integer $s\ge0$, the local shuffle map $\overline\Phi$ restricts to a chain isomorphism
\[
 \sum_{\substack{a,b\ge0\\a+b\le s}}
 F_aQ_*(P,T)\otimes_{\Bbbk}F_bQ_*(Q,U)
 \xrightarrow{\ \cong\ } F_sQ_*(P\Box Q,T\times U).
\]
The sum is taken inside $Q_*(P,T)\otimes_{\Bbbk}Q_*(Q,U)$ with its total grading.
\end{theorem}
\begin{proof}
For a path $\alpha=(x_0,\ldots,x_p)$ in $P$, let $\rho_T(\alpha)=\max_i d_P(x_i,T)$. Define $\rho_U$ and $\rho_{T\times U}$ in the same way. Write $\operatorname{Sh}(\alpha,\beta)$ for the shuffles of $\alpha$ and $\beta=(y_0,\ldots,y_q)$. Since each pair of positions $(i,j)$ occurs on some shuffle, we have
\[
 \begin{aligned}
 \max_{\gamma\in\operatorname{Sh}(\alpha,\beta)}\rho_{T\times U}(\gamma)
 &=\max_{i,j}\bigl(d_P(x_i,T)+d_Q(y_j,U)\bigr)\\
 &=\rho_T(\alpha)+\rho_U(\beta).
 \end{aligned}
\]
If $z$ has elementary tensor coefficients $z_{\alpha,\beta}$, each shuffle of $(\alpha,\beta)$ has coefficient $\pm z_{\alpha,\beta}$ in $\Phi(z)$. Thus, $\Phi(z)$ is supported in $N_s(T\times U)$ exactly when
\[
 z_{\alpha,\beta}\ne0\quad\Longrightarrow\quad
 \rho_T(\alpha)+\rho_U(\beta)\le s.
\]

For filtered spaces $A,B$ with filtrations $A_a,B_b$ and subspaces $V'\subseteq A$, $W'\subseteq B$, we have
\[
 (V'\otimes W')\cap\sum_{a+b\le s}A_a\otimes B_b
 =\sum_{a+b\le s}(V'\cap A_a)\otimes(W'\cap B_b).
\]
Choose bases of $A$ and $B$ such that every $A_a$ and $B_b$ is spanned by a subset of the corresponding basis. Note that $V'\cap A_a$ and $W'\cap B_b$ are spanned by the basis elements in $V'$ and $W'$, respectively. Both sides are spanned by the tensors $v\otimes w$ with
\[
 v\in V'\cap A_a,\quad w\in W'\cap B_b,\quad a+b\le s.
\]

Let $P_a=P|_{N_a(T)}$ and $Q_b=Q|_{N_b(U)}$. For bidegree $(p,q)$, we take $V'=\Omega_p(P)$, $W'=\Omega_q(Q)$, $A_a=\mathcal A_p(P_a)$, and $B_b=\mathcal A_q(Q_b)$. By Lemma~\ref{lem:intersection}, we have
\[
 V'\cap A_a=\Omega_p(P_a),\quad W'\cap B_b=\Omega_q(Q_b).
\]
It follows that
\[
 \Phi^{-1}\Omega_*\bigl((P\Box Q)|_{N_s(T\times U)}\bigr)
 =\sum_{a+b\le s}\Omega_*(P_a)\otimes\Omega_*(Q_b).
\]
The quotient map $\Omega_*(P)\to Q_*(P,T)$ sends $\Omega_*(P_a)$ to $F_aQ_*(P,T)$, and $\Omega_*(Q)\to Q_*(Q,U)$ maps $\Omega_*(Q_b)$ onto $F_bQ_*(Q,U)$. Taking tensor products and summing over $a+b\le s$ gives the product filtration.
\end{proof}

In particular, the shuffle product of chains in filtration degrees $a$ and $b$ lies in filtration degree $a+b$.

\subsection{Sharp radius bounds}\label{sec:cube-local-model}

The examples below achieve the bounds $s\ge r+1$ for isomorphism and $s\ge r$ for surjectivity.

Let $\Lambda V=\bigoplus_{p\ge0}\Lambda^p V$ be the exterior algebra of $V$ with $\Lambda^0 V=\Bbbk$. Given a basis $e_1,\ldots,e_n$ of $V$, the wedges $e_{i_1}\wedge\cdots\wedge e_{i_p}$ with $i_1<\cdots<i_p$ form a basis of $\Lambda^p V$. Thus, $\Lambda^p V$ has dimension $\binom np$ for $p\le n$ and vanishes for $p>n$.

Let $I$ denote the directed interval $0\to1$. For $n\ge 1$, let $C_n=I^{\Box n}$ be the directed $n$-cube. We denote the source of $C_n$ by $o=(0,\ldots,0)$. The undirected distance between two vertices of $C_n$ is the number of coordinates in which they differ.

\begin{proposition}\label{prop:cube-vertex-local}
For each vertex $v=(\varepsilon_1,\ldots,\varepsilon_n)\in\{0,1\}^n$, we have
\[
 \LH_q(C_n,v)=0\quad\text{for each }q\ge0.
\]
\end{proposition}
\begin{proof}
For $\varepsilon\in\{0,1\}$, the local complex of the directed interval is
\[
 Q_*(I,\varepsilon):\quad
 0\longrightarrow\Bbbk\xrightarrow{\,2\varepsilon-1\,}\Bbbk\longrightarrow0,
\]
in degrees $1,0$. Each of the two choices of endpoint gives an invertible differential. By Theorem~\ref{thm:local-kunneth}, we have
\[
 Q_*(C_n,v)\cong\bigotimes_{i=1}^n Q_*(I,\varepsilon_i),
\]
and the local K\"unneth formula yields the asserted vanishing.
\end{proof}

Multiplying the degree-one generator of each factor by $2\varepsilon_i-1$ changes its differential to $1$. The tensor product then identifies with
\[
 Q_p(C_n,v)\cong\Lambda^p\Bbbk^n,\quad
 \partial=\iota_\lambda,\quad \lambda(e_i)=1,
\]
where $e_1,\ldots,e_n$ is the standard basis and contraction by $\lambda$ is given by
\[
 \iota_\lambda(e_{i_1}\wedge\cdots\wedge e_{i_p})
 =\sum_{j=1}^p(-1)^{j-1}
 e_{i_1}\wedge\cdots\wedge\widehat{e_{i_j}}\wedge\cdots\wedge e_{i_p}.
\]
In degree zero, $\iota_\lambda=0$. Put $W=\ker\lambda$. Then $\Bbbk^n=\Bbbk e_1\oplus W$, and $\iota_\lambda(e_1\wedge\eta)=\eta$ and $\iota_\lambda\eta=0$ for $\eta\in\Lambda W$. Hence, we have $\ker(\iota_\lambda:\Lambda^p\Bbbk^n\to\Lambda^{p-1}\Bbbk^n)=\Lambda^pW$ for $p\ge1$, with the same kernel formula in degree zero.

\begin{theorem}\label{thm:cube-balls}
Let $v$ be a vertex of $C_n$. For each integer $0\le s<n$, we have
\begin{equation*}
 H_q(B_s(v),B_s(v)-v)\cong
 \begin{cases}\Bbbk^{\binom{n-1}{s}},&q=s,\\0,&q\ne s.\end{cases}
\end{equation*}
For $s\ge n$, all these relative homology groups vanish.
\end{theorem}
\begin{proof}
Write $v=(\varepsilon_1,\ldots,\varepsilon_n)$. In the $i$-th interval, the radius-zero ball contains only $\varepsilon_i$, and the radius-one ball contains the edge $0\to1$. Each elementary tensor in the standard basis of $Q_p(C_n,v)$ has $p$ edge factors, so it lies in $F_sQ_*(C_n,v)$ exactly when $p\le s$. By Theorem~\ref{thm:chain-locality}, one has $F_sQ_*(C_n,v)=Q_*(B_s(v),v)$. It follows that
\[
 Q_p(B_s(v),v)\cong
 \begin{cases}\Lambda^p\Bbbk^n,&0\le p\le s,\\0,&p>s.\end{cases}
\]
Suppose $s<n$ and write $C_*=Q_*(B_s(v),v)$. For $0\le q<s$, the maps $\partial_q$ and $\partial_{q+1}$ in $C_*$ agree with those in the exact complex $Q_*(C_n,v)$. Thus, we have $H_q(C_*)=0$. One has $\operatorname{im}\partial_{s+1}=0$ as $C_{s+1}=0$. It follows that
\[
 H_s(B_s(v),B_s(v)-v)
 \cong\ker\bigl(\iota_\lambda|_{\Lambda^s\Bbbk^n}\bigr)
 =\Lambda^s W\cong\Bbbk^{\binom{n-1}{s}}.
\]
If $s\ge n$, then $B_s(v)=C_n$, whose local homology at $v$ vanishes by Proposition~\ref{prop:cube-vertex-local}.
\end{proof}

\begin{corollary}\label{cor:iso-radius-optimal}
For $r\ge0$ and any vertex $v$ of $C_{r+1}$, we have
\[
 H_r(B_r(v),B_r(v)-v)\cong\Bbbk,\quad \LH_r(C_{r+1},v)=0.
\]
\end{corollary}
\begin{proof}
Apply Theorem~\ref{thm:cube-balls} with $n=r+1$ and $s=r$, and use Proposition~\ref{prop:cube-vertex-local}.
\end{proof}

Thus, the uniform isomorphism bound $s\ge r+1$ cannot be lowered to $s\ge r$, even for finite directed acyclic graphs. For fixed $r\ge1$, increasing $n>r$ makes the kernel dimension $\binom{n-1}{r}$ of the radius-$r$ comparison arbitrarily large, at every vertex of $C_n$.

\begin{example}\label{ex:square-transient}
Let $o=(0,0)$. In the radius-one ball $B_1(o)$ of $C_2$, $c=e_{00,10}-e_{00,01}$ represents a nonzero local $1$-class. The ball of radius two also contains $11$, and $z=e_{00,10,11}-e_{00,01,11}\in\Omega_2$ satisfies $[\partial z]=[c]$. Thus, the class becomes a boundary.
\end{example}

We next prove the sharpness of the surjectivity bound.

Let $J=(-1\to0\to1)$, $D_r=J^{\Box r}$, and $o=(0,\ldots,0)$, with $r\ge1$. Example~\ref{ex:interval} gives
\[
 Q_*(J,0):\quad0\longrightarrow\Bbbk^2
 \xrightarrow{(1,-1)}\Bbbk\longrightarrow0.
\]
Thus, $\LH_q(D_r,o)$ is nonzero only for $q=r$, where it equals $\Bbbk$. For one factor, the local $0$-chain appears at radius zero and the local $1$-chains at radius one. By Theorem~\ref{thm:filtered-kunneth}, $F_sQ_*(D_r,o)$ is the brutal truncation of the full local chain complex in degrees $\le s$.

\begin{corollary}\label{cor:surj-radius-optimal}
For any $r\ge1$, at the central vertex of $D_r$, we have
\[
 Q_r(B_{r-1}(o),o)=0,\quad
 H_r(B_{r-1}(o),B_{r-1}(o)-o)=0\longrightarrow\LH_r(D_r,o)=\Bbbk.
\]
The surjectivity bound $s\ge r$ is therefore optimal in degree $r$.
\end{corollary}

In $D_r$, the generator of $\LH_r(D_r,o)$ first has a representative at radius $r$. In $C_{r+1}$, the class at radius $r$ bounds at radius $r+1$.

\subsection{Detecting local regularity}\label{sec:regularity-radius}

Given integers $N\ge n\ge0$, the conditions
\[
 \LH_q(P,v)\cong
 \begin{cases}\Bbbk,&q=n,\\0,&q\ne n,\end{cases}
 \quad 0\le q\le N,
\]
can be checked at radius $N+1$ by Theorem~\ref{thm:chain-locality}. Local regularity of dimension $n$ additionally requires vanishing homology in degrees above $N$. A bound on chain degrees makes this a finite verification condition.

\begin{proposition}\label{prop:regularity-radius}
Suppose $\Omega_p(P)=0$ for each $p>D$, where $D$ is a nonnegative integer. For each $v\in V(P)$, the inclusion $P|_{N_D(v)}\hookrightarrow P$ induces a chain isomorphism
\[
 Q_*(P|_{N_D(v)},v)\xrightarrow{\cong}Q_*(P,v).
\]
\end{proposition}
\begin{proof}
For $p\le D$, the map is an isomorphism by Theorem~\ref{thm:chain-locality}. For $p>D$, both local chain groups vanish since $\Omega_p(P|_{N_D(v)})\subseteq\Omega_p(P)=0$.
\end{proof}

The induced $D$-neighborhood of $v$ determines whether $v$ is regular of a given dimension. Let $G$ be a finite directed acyclic graph. Take $D$ to be the maximum length of a directed path in $G$. Then $\Omega_p(G)=0$ for $p>D$.

Let $\Bbbk[m]$ denote the chain complex with $\Bbbk$ concentrated in degree $m$ and zero elsewhere, with trivial differential.

Let $G=A\vee_o B$ be the digraph obtained by identifying the source vertices of two disjoint digraphs $A$ and $B$ at $o$.

\begin{lemma}\label{lem:source-wedge}
There is a short exact sequence of chain complexes
\[
 0\longrightarrow\Bbbk[0]
 \xrightarrow{1\mapsto([e_o],-[e_o])}
 Q_*(A,o)\oplus Q_*(B,o)
 \xrightarrow{\,\psi\,} Q_*(G,o)\longrightarrow0.
\]
\end{lemma}
\begin{proof}
Let $i_A:A\hookrightarrow G$ and $i_B:B\hookrightarrow G$ be inclusions of digraphs. Set $\psi(u,w)=(i_A)_\#u+(i_B)_\#w$. Since $o$ is a source, each allowed path of positive length lies in one branch. For $p\ge2$, a regular $(p-1)$-path occurring in boundaries from both branches would be supported on $\{o\}$, which is impossible in positive degree. Thus, the boundary condition holds separately in the two branches. With $\Omega_1=\mathcal A_1$, this gives
\[
 \Omega_p(G)=\Omega_p(A)\oplus\Omega_p(B),\quad p\ge1.
\]
Since $G-o$ is the disjoint union of $A-o$ and $B-o$, we also have $\Omega_p(G-o)=\Omega_p(A-o)\oplus\Omega_p(B-o)$, so $Q_p(G,o)\cong Q_p(A,o)\oplus Q_p(B,o)$ for $p\ge1$. In degree zero, the inclusions of $A$ and $B$ into $G$ induce the addition map
\[
 \Bbbk[e_o]\oplus\Bbbk[e_o]\longrightarrow\Bbbk[e_o],\quad (a[e_o],b[e_o])\longmapsto(a+b)[e_o].
\]
The kernel in degree zero is $\Bbbk([e_o],-[e_o])$. Its generator has zero boundary, so the kernel complex is $\Bbbk[0]$.
\end{proof}

\begin{proposition}\label{prop:no-uniform-regularity-radius}
For each $s\ge0$ and $n\ge1$, there exist finite pointed directed acyclic graphs $(G_1,v_1)$ and $(G_2,v_2)$ with isomorphic rooted $s$-balls such that $v_1$ is regular of dimension $n$, while $\LH_n(G_2,v_2)\ne0$ and $\LH_q(G_2,v_2)\ne0$ for some $q>n$.
\end{proposition}
\begin{proof}
Pick an integer $d>\max\{s,1\}$. Let $C=I^{\Box(d+1)}$ with source $o$, and set $T=B_d^C(o)$. Attach the source of an extra directed interval to $o$ in both graphs, yielding $G=I\vee_o C$ and $G'=I\vee_o T$. By Proposition~\ref{prop:cube-vertex-local} and Theorem~\ref{thm:cube-balls}, $H_q(Q_*(I,0))=H_q(Q_*(C,o))=0$ for all $q\ge0$, whereas $H_q(Q_*(T,o))$ equals $\Bbbk$ when $q=d$ and vanishes otherwise.

For $Y\in\{C,T\}$, the direct sum $Q_*(I,0)\oplus Q_*(Y,o)$ has zero homology in degrees $0$ and $1$ since $d>1$. The connecting map in the long exact sequence of Lemma~\ref{lem:source-wedge} gives $\LH_1(I\vee_o Y,o)\cong H_0(\Bbbk[0])=\Bbbk$. The same sequence gives $\LH_0(I\vee_o Y,o)=0$ and $\LH_q(I\vee_o Y,o)\cong\LH_q(Y,o)$ for $q\ge2$. The local Poincar\'e series are
\[
 \ell_{G,o}(t)=t,\quad \ell_{G',o}(t)=t+t^d.
\]
The condition $s<d$ ensures the induced $s$-balls coincide. Now take the product of each graph with $n-1$ copies of $J=(-1\to0\to1)$, pointed at their central vertices. The local Poincar\'e series become $t^n$ and $t^n+t^{d+n-1}$. Additivity of Cartesian distance preserves the isomorphism of the $s$-balls. No product factor is needed when $n=1$.
\end{proof}

For a fixed $r$, Theorem~\ref{thm:chain-locality} allows us to compute $\LH_r(P,v)$ at radius $r+1$. Proposition~\ref{prop:no-uniform-regularity-radius} shows that a radius bound depending only on $n$ cannot force homology to vanish in all degrees greater than $n$, even for finite directed acyclic graphs.

\subsection{Radius persistence}\label{sec:radius-persistence}

Let $P$ be a regular path complex, and let $\varnothing\ne T\subseteq V(P)$. We employ the radius filtration $F_sQ_*(P,T)$ introduced above, indexed by $s\in\mathbb Z_{\ge0}$. By Theorem~\ref{thm:chain-locality}, the inclusion $P|_{N_s(T)}\hookrightarrow P$ induces the identification
\[
 F_sQ_*(P,T)\cong Q_*(P|_{N_s(T)},T).
\]
For each $s\ge0$, let $\iota_s:F_sQ_*(P,T)\hookrightarrow F_{s+1}Q_*(P,T)$ denote the inclusion map. Naturality in this subsection is considered with respect to morphisms of pairs $(P,P-T)$ that preserve the radius filtrations. One sufficient condition for a path-complex morphism $f:P\to P'$ is $f^{-1}(T')=T$. Such a map sends unmarked vertices to unmarked vertices, and satisfies $d_{P'}(f(x),T')\le d_P(x,T)$, because every allowed edge maps either to an edge or to a single vertex. Consequently, $f$ maps $N_s(T)$ into $N_s(T')$ and induces $F_sQ_*(P,T)\to F_sQ_*(P',T')$ for each $s$. This criterion applies to each factor in Cartesian products.

\begin{definition}\label{def:radius-rees}
The \textit{Rees chain complex} associated to the radius filtration is the graded $\Bbbk[t]$-chain complex defined by
\[
 \mathscr Q_p(P,T):=\bigoplus_{s\ge0}F_sQ_p(P,T)t^s,
 \quad \partial(ct^s):=(\partial c)t^s.
\]
An element $ct^s$ has homological degree $p$ and radius degree $s$. Multiplication by $t$ is given by $t\cdot(ct^s):=\iota_s(c)t^{s+1}$. The differential is $\Bbbk[t]$-linear, lowers homological degree by one, and preserves radius degree.
\end{definition}

\begin{definition}\label{def:radius-persistence}
For $r\ge0$, the \textit{radius-persistent local homology} in degree $r$ is the graded $\Bbbk[t]$-module
\[
 \mathsf L_r(P,T):=H_r(\mathscr Q_*(P,T)).
\]
Its radius grading and $t$-action are inherited from the Rees chain complex.
\end{definition}

Given a cycle $c\in F_sQ_r(P,T)$, write its homology class at radius $s$ as $[c]_s$. Since the differential preserves radius degree, the assignment $[ct^s]\mapsto[c]_s$ identifies the degree-$s$ component of $\mathsf L_r(P,T)$ with $H_r(F_sQ_*(P,T))$. Under this identification, multiplication by $t$ corresponds to the map induced by inclusion:
\[
 t\cdot[c]_s=(\iota_s)_*([c]_s)=[c]_{s+1}.
\]
For $c\in F_sQ_r(P,T)$ with $\partial c=0$, we have $t[c]_s=0$ if and only if $c\in\partial F_{s+1}Q_{r+1}(P,T)$.

\begin{proposition}\label{prop:radius-rees-product}
Each $\mathscr Q_p(P,T)$ is a finitely generated free graded $\Bbbk[t]$-module. For regular path complexes $P,Q$ on finite vertex sets and marked subsets $T\subseteq V(P)$, $U\subseteq V(Q)$, the shuffle product induces a natural isomorphism of graded chain complexes
\[
 \mathscr Q_*(P,T)\otimes_{\Bbbk[t]}\mathscr Q_*(Q,U)
 \xrightarrow{\cong}\mathscr Q_*(P\Box Q,T\times U).
\]
\end{proposition}
\begin{proof}
Choose a basis $c_1,\ldots,c_m$ of $Q_p(P,T)$ such that each $F_sQ_p(P,T)$ is spanned by a subset of this basis. Put $a_i=\min\{s:c_i\in F_sQ_p(P,T)\}$. Then
\[
 F_sQ_p(P,T)=\bigoplus_{a_i\le s}\Bbbk c_i,\quad
 \mathscr Q_p(P,T)=\bigoplus_{i=1}^m\Bbbk[t]\,c_it^{a_i}.
\]
Thus $\mathscr Q_p(P,T)$ is graded free on the $m$ generators $c_it^{a_i}$.

For $c\in F_aQ_p(P,T)$ and $d\in F_bQ_q(Q,U)$, the filtered shuffle product defines the $\Bbbk[t]$-balanced map
\[
 (ct^a,dt^b)\longmapsto(c\times d)t^{a+b}.
\]
The induced map on the tensor product preserves the radius grading and commutes with the differential by
\[
 \partial(c\times d)=(\partial c)\times d+(-1)^p c\times(\partial d).
\]
Choose a basis $\{d_j\}$ of $Q_q(Q,U)$ in the same way, with $b_j=\min\{s:d_j\in F_sQ_q(Q,U)\}$. By Theorem~\ref{thm:filtered-kunneth}, the elements
\[
 (c_i\times d_j)t^{a_i+b_j},\quad p+q=n,
\]
form a $\Bbbk[t]$-basis of $\mathscr Q_n(P\Box Q,T\times U)$ as $p,q,i,j$ vary. These are the images of the basis elements $c_it^{a_i}\otimes d_jt^{b_j}$ in degree $n$. The isomorphism commutes with induced maps because $(f\times g)_\#(c\times d)=f_\#c\times g_\#d$ for the morphisms in Proposition~\ref{prop:local-kunneth-naturality} that preserve the radius filtrations.
\end{proof}

\begin{theorem}\label{thm:persistent-kunneth}
For $r\ge0$, there exists a natural short exact sequence of graded $\Bbbk[t]$-modules
\[
 \begin{aligned}
 0\longrightarrow{}&
 \bigoplus_{i+j=r}\mathsf L_i(P,T)\otimes_{\Bbbk[t]}\mathsf L_j(Q,U)\\
 \longrightarrow{}&\mathsf L_r(P\Box Q,T\times U)\\
 \longrightarrow{}&
 \bigoplus_{i+j=r-1}\operatorname{Tor}^{\Bbbk[t]}_1
       \bigl(\mathsf L_i(P,T),\mathsf L_j(Q,U)\bigr)
 \longrightarrow0.
 \end{aligned}
\]
The sequence splits noncanonically as graded $\Bbbk[t]$-modules.
\end{theorem}
\begin{proof}
Proposition~\ref{prop:radius-rees-product} identifies the product complex with the tensor product of two complexes free in each degree over $\Bbbk[t]$. Apply the graded K\"unneth theorem over $\Bbbk[t]$ and Definition~\ref{def:radius-persistence} \cite[Section~5.2]{GakharPerea}.
\end{proof}

\begin{definition}\label{def:radius-barcode}
For integers $0\le a<b\le\infty$, the \textit{interval module} $\Bbbk_{[a,b)}:(\mathbb{R},\leq)\to \mathbf{Vec}_\Bbbk$ is a persistence module defined by
\[
\Bbbk_{[a,b)}(s)=\begin{cases}
                                                                                                              \Bbbk, & \mbox{if } a\le s<b \\
                                                                                                              0, & \mbox{otherwise}.
                                                                                                            \end{cases}
\]

The collection of intervals in a direct sum decomposition of $\mathsf L_r(P,T)$ forms its \textit{barcode}. The endpoints $a$ and $b$ of $\Bbbk_{[a,b)}$ are called the \textit{birth radius} and the \textit{death radius}, respectively.
\end{definition}

\begin{proposition}\label{prop:radial-barcode}
Each interval in the barcode of $\mathsf L_r(P,T)$ begins by radius $r$, and each finite interval ends by radius $r+1$. The barcode contains exactly $\dim_{\Bbbk}\LH_r(P,T)$ infinite intervals.
\end{proposition}
\begin{proof}
Write $Q_*=Q_*(P,T)$ and $Q_*^s=F_sQ_*$. For $s\ge r$, Theorem~\ref{thm:chain-locality} gives $Q_r^s=Q_r$. Since $Q_*^s$ is a subcomplex of $Q_*$, the two complexes have the same $r$-cycles. The inclusion $Q_*^s\hookrightarrow Q_*^{s+1}$ therefore induces the quotient map
\[
 Z_r(Q_*)/\partial Q_{r+1}^s
 \longrightarrow Z_r(Q_*)/\partial Q_{r+1}^{s+1},\quad [z]\longmapsto[z].
\]
This map is surjective. For $s\ge r+1$, we also have $Q_{r+1}^s=Q_{r+1}$, so $H_r(Q_*^s)\to\LH_r(P,T)$ is an isomorphism.

Each $H_r(Q_*^s)$ is finite-dimensional, and the structure maps are isomorphisms for $s\ge r+1$. Thus $\mathsf L_r(P,T)$ is finitely generated over $\Bbbk[t]$ and decomposes into a finite sum of interval modules. Surjectivity for $s\ge r$ gives the birth bound $a\le r$. Stabilization at $r+1$ gives the death bound $b\le r+1$ for each finite interval. At radius $r+1$, each infinite interval contributes one dimension to $\LH_r(P,T)$.
\end{proof}

\begin{example}\label{ex:cube-radius-barcodes}
For any vertex $v$ of $C_n$, Theorem~\ref{thm:cube-balls} gives
\[
 \mathsf L_q(C_n,v)\cong
 \begin{cases}
 \Bbbk_{[q,q+1)}^{\,\binom{n-1}{q}},&0\le q<n,\\
 0,&q\ge n.
 \end{cases}
\]
A nonzero class in degree $q$ is born at radius $q$ and dies at radius $q+1$. For $D_r=J^{\Box r}$ with center $o$, the filtered local complex is obtained by truncating above the radius degree. Its degree-$r$ homology is zero below radius $r$ and is $\Bbbk$ at every radius at least $r$. Thus, we obtain
\[
 \mathsf L_r(D_r,o)\cong\Bbbk_{[r,\infty)}.
\]
\end{example}

\begin{example}\label{ex:radial-tor}
Let $A=\Bbbk[t]$ and $M=A/(t)$, with its generator in radius degree zero. For the source of $I=(0\to1)$, we have $\mathsf L_0(I,0)=M$ and $\mathsf L_q(I,0)=0$ for $q>0$. The module $M$ corresponds to $\Bbbk_{[0,1)}$.

To compute the graded $\operatorname{Tor}$ term, use the free resolution
\[
 0\longrightarrow Ae\xrightarrow{\ e\mapsto t\ }A\longrightarrow M\longrightarrow0,
 \quad \deg e=1.
\]
Tensoring with $M$ makes the differential zero, since $tM=0$. Thus, $M\otimes_A M\cong M$ has its generator in radius degree zero, whereas $\operatorname{Tor}^A_1(M,M)\cong Me$ has its generator in radius degree one. Applying Theorem~\ref{thm:persistent-kunneth} to $C_2=I\Box I$, we obtain
\[
 \mathsf L_0(C_2,(0,0))\cong\Bbbk_{[0,1)},\quad
 \mathsf L_1(C_2,(0,0))\cong\Bbbk_{[1,2)}.
\]
The class of $e$ in $\operatorname{Tor}_1^A(M,M)$ has radius $1$ and is killed by $t$. It represents the interval $[1,2)$ in Example~\ref{ex:square-transient}.
\end{example}

\begin{definition}\label{def:stabilization-radius}
The \textit{stabilization radius} $\rho_r(P,T)$ is the least integer $s\ge0$ such that, for each integer $a\ge s$, the inclusion $F_aQ_*(P,T)\hookrightarrow Q_*(P,T)$ induces an isomorphism
\[
 H_r(F_aQ_*(P,T))\xrightarrow{\cong}\LH_r(P,T).
\]
\end{definition}

Theorem~\ref{thm:chain-locality} gives $\rho_r(P,T)\le r+1$. For each vertex $v$ of $C_{r+1}$, Corollary~\ref{cor:iso-radius-optimal} identifies the map induced by $B_r(v)\hookrightarrow C_{r+1}$ in degree $r$ as
\[
 H_r(B_r(v),B_r(v)-v)\cong\Bbbk\longrightarrow\LH_r(C_{r+1},v)=0.
\]
Thus, $\rho_r(C_{r+1},v)>r$, and the upper bound gives $\rho_r(C_{r+1},v)=r+1$.

\begin{example}\label{ex:early-radius-isomorphism}
Fix $r\ge1$ and $v\in V(C_{r+1})$. In degree $r$, the comparison map is $0\to0$ at each radius $s<r$ and $\Bbbk\to0$ at radius $r$. Thus, the first radius giving an isomorphism can be strictly smaller than the stabilization radius.
\end{example}

\section{Geometric models for local path homology}\label{sec:simplicial}\label{sec:geometry}

In this section, we use simplicial and cubical models to compute local path homology and compare it with the classical local homology of the corresponding geometric realization.

\subsection{The simplicial link formula}

Let $K$ be a finite simplicial complex, and assign distinct labels to its vertices to define $P(K)$. By Proposition~\ref{prop:simplicial-star-link} and Corollary~\ref{cor:perfect-star-local}, we have $L_v(P(K))=P(\Lk_K(v))$ and
\[
 \LH_r(P(K),v)\cong\widetilde H_{r-1}(L_v(P(K)))
 \cong\widetilde H_{r-1}(\Lk_K(v);\Bbbk).
\]
The labels on the simplicial link are inherited from $K$. The simplicial chain identification further yields
\[
 \LH_r(P(K),v)\cong H_r(K,K-v;\Bbbk)
 \cong H_r(|K|,|K|\setminus\{v\};\Bbbk).
\]
The second isomorphism follows from simplicial–singular comparison. Normalizing barycentric coordinates away from $v$ retracts $|K|\setminus\{v\}$ onto $|K-v|$ \cite{Hatcher}. This relates local path homology to classical local homology.

\begin{definition}
A digraph $G$ is called a \textit{transitive directed acyclic graph} if it contains no directed cycles and satisfies the condition that $x\to y$ and $y\to z$ imply $x\to z$.
\end{definition}

A Hasse diagram has an edge only for a covering relation. For the chain $a<v<b$, its edges are $a\to v$ and $v\to b$, while the transitive digraph also has $a\to b$.

\begin{corollary}\label{cor:transitive-link}
Let $G$ be a transitive directed acyclic graph, and let $\Delta(G)$ denote the order complex of its vertex poset. For any $v\in V(G)$, we have
\[
 \Delta(S_v)=\overline{\st}_{\Delta(G)}(v),\quad
 \Delta(L_v)=\Lk_{\Delta(G)}(v).
\]
The composite $\delta_r j_r^{-1}$ gives $\LH_r(G,v)\cong\widetilde H_{r-1}(\Lk_{\Delta(G)}(v);\Bbbk)$ for all $r\ge0$.
\end{corollary}
\begin{proof}
The allowed paths are precisely the increasing chains, so $P(G)=P(\Delta(G))$. Sending $e_{x_0\cdots x_p}$ to $[x_0,\ldots,x_p]$ identifies $\Omega_*(G)=\Sigma_*(G)=\mathcal A_*(G)$ with $C_*(\Delta(G);\Bbbk)$ \cite{GMY}. The star and link identities follow from Propositions~\ref{prop:star-link-digraph} and~\ref{prop:simplicial-star-link}. By Corollaries~\ref{cor:perfect-star-local} and~\ref{cor:perfect-completed-link}, we have $D_*(G,v)=0$ and $\widehat{\mathcal L}_*(G,v)=\Omega_*(L_v)$. The closed star is a cone, so the composite $\delta_rj_r^{-1}:\LH_r(G,v)\to\widetilde H_{r-1}(L_v)$ is an isomorphism.
\end{proof}

More generally, let $X=A\cup B$ be a subcomplex cover of a simplicial complex, with $A,B$ equipped with vertex labels from $X$. Under the simplicial chain identifications, one has
\[
 \Sigma_*(X)=\Omega_*(X)=C_*(X),\quad K_*(A,B)=C_*(A\cap B),\quad \mathcal T_*=0.
\]
These identifications lead to the classical Mayer--Vietoris sequence and excision isomorphism within path homology.

\begin{proposition}\label{prop:perfect-radius-one}
Let $P$ be a perfect regular path complex, and let $v\in V(P)$. The inclusion $P|_{N_1(v)}\hookrightarrow P$ induces a chain isomorphism
\[
 Q_*(P|_{N_1(v)},v)\xrightarrow{\cong}Q_*(P,v).
\]
\end{proposition}
\begin{proof}
Since $P$ is perfect, we have $\Omega_p(P)=\mathcal A_p(P)$. The same holds for each induced restriction of $P$. It follows that $Q_p(P,v)$ has basis $[e_\alpha]$ indexed by the allowed $p$-paths $\alpha$ containing $v$.

Let $x\ne v$ be a vertex of such a path $\alpha$. Deleting all vertices except $x$ and $v$ gives an allowed edge, either $xv$ or $vx$. It follows that $x\in N_1(v)$ and $\alpha\in P|_{N_1(v)}$. The two local chain groups therefore have the same basis in each degree, and the inclusion sends each $[e_\alpha]$ to itself.
\end{proof}

Let $P$ be a perfect monotone path complex. A vertex $v$ of $P$ is regular of dimension $n$ if and only if
\[
 \widetilde H_q(L_v(P))\cong
 \begin{cases}\Bbbk,&q=n-1,\\0,&q\ne n-1.\end{cases}
\]
The regularity criterion also holds for any regular path complex $P$ if the comparison maps $j_r$ and $\delta_r$ of Corollary~\ref{cor:iso-conditions} are isomorphisms for each $r\ge0$.

\subsection{The cubical chain model}\label{sec:positive-grid}

\begin{definition}\label{def:oriented-lattice}
Let $n\ge1$. A \textit{lattice digraph in $\mathbb Z^n$} is a finite digraph $G=(V,E)$ with $V\subseteq\mathbb Z^n$ and edges being one-step directed edges.
\end{definition}

To form a lattice digraph $G$ on $V$, choose any subset of the nearest-neighbor lattice edges with both endpoints in $V$ and orient each chosen edge independently. The resulting digraph may contain directed cycles.

\begin{definition}\label{def:grid-cubical-complex}
Let $\sigma=a+\sum_{j=1}^p[0,1]e_{i_j}$ be an axis-parallel unit cube, with $a\in\mathbb Z^n$, $p\ge1$, and $i_1<\cdots<i_p$. We say that $\sigma$ is a \textit{cube with a product orientation in $G$} if all its vertices and edges belong to $G$ and, for each $j$, all its edges parallel to $e_{i_j}$ are directed along the same vector $\varepsilon_je_{i_j}$, where $\varepsilon_j\in\{\pm1\}$. Let $X(G)$ consist of these cubes, together with the vertices of $G$ as zero-dimensional cubes. The product orientation restricts to every face, so $X(G)$ is a finite cubical complex. We write $X_p(G)$ for its set of $p$-dimensional cubes.
\end{definition}

For a cube $\sigma$, we denote its vertex set by $\operatorname{Vert}(\sigma)$. Explicitly, for $\sigma=a+\sum_{j=1}^p[0,1]e_{i_j}$, we have
\[
 \operatorname{Vert}(\sigma)
 =\left\{a+\sum_{j=1}^p\delta_je_{i_j}:\delta_j\in\{0,1\}\right\}.
\]
For a zero-dimensional cube $\sigma=\{a\}$, we have $\operatorname{Vert}(\sigma)=\{a\}$.

The complex $X(G)$ includes only cubes with product orientations. If every edge is $x\to x+e_i$, a cube belongs to it precisely when all its vertices and edges are present.

Orient $\sigma=a+\sum_{j=1}^p[0,1]e_{i_j}$ by the ordered basis $(e_{i_1},\ldots,e_{i_p})$, where $i_1<\cdots<i_p$. Let $\sigma_j^-$ and $\sigma_j^+$ be the faces given by $x_{i_j}=a_{i_j}$ and $x_{i_j}=a_{i_j}+1$, respectively. With each face oriented by its remaining coordinate directions, the cubical boundary is
\[
 \partial_{\mathrm{cub}}\sigma
 =\sum_{j=1}^p(-1)^{j-1}(\sigma_j^+-\sigma_j^-).
\]
For $\sigma\in X_p(G)$, write $u_j=\varepsilon_j e_{i_j}$ for the directed edge step along the $i_j$-axis. The source of $\sigma$ is
\[
 b=a+\sum_{\varepsilon_j=-1}e_{i_j}.
\]
For $\pi\in S_p$, let $\gamma_\pi$ be the path starting at $b$ with successive steps $u_{\pi(1)},\ldots,u_{\pi(p)}$. Define
\[
 \Psi(\sigma)=\sum_{\pi\in S_p}\operatorname{sgn}(\pi)e_{\gamma_\pi},
 \quad \Phi(\sigma)=\varepsilon_\sigma\Psi(\sigma).
\]
Here, $\varepsilon_\sigma=\prod_{j=1}^p\varepsilon_j$ is the sign of the ordered basis $(u_1,\ldots,u_p)$ relative to $(e_{i_1},\ldots,e_{i_p})$. In degree zero, we set $\Phi(a)=\Psi(a)=e_a$.

\begin{theorem}\label{thm:grid-cubical-chains}
For a lattice digraph $G$ in $\mathbb Z^n$, we have a chain isomorphism
\[
 \Phi:C_*^{\mathrm{cub}}(X(G);\Bbbk)\xrightarrow{\cong}\Omega_*(G),
\]
which is natural with respect to subgraph inclusions and compatible with induced vertex deletions. In particular, $\Omega_p(G)=0$ for all $p>n$.
\end{theorem}
\begin{proof}
In degree zero, $\Phi$ maps each vertex $x$ to $e_x$. For $p\ge1$, write $c=\sum_\gamma c_\gamma e_\gamma\in\Omega_p(G)$. Fix an allowed path $\gamma=(x_0,\ldots,x_p)$ and write $w_k=x_k-x_{k-1}$. Two consecutive steps cannot be opposites, since this would lead to both orientations of the same edge.

For $1\le k<p$, write $w_k=\pm e_i$ and $w_{k+1}=\pm e_j$. Deleting $x_k$ from $\gamma$ gives a path $\eta_k$ with a step of lattice length two from $x_{k-1}$ to $x_{k+1}$. Since $\partial c$ is allowed, the coefficient of $e_{\eta_k}$ in $\partial c$ is zero.

If $i=j$, the two steps follow the same axis. Since they cannot be opposite, $w_k=w_{k+1}$. The only intermediate lattice vertex is $x_k$, so the coefficient of $e_{\eta_k}$ is $(-1)^k c_\gamma$. Hence $c_\gamma=0$.

If $i\ne j$, the two coordinate steps can be taken in either order in the lattice. Taking $w_k$ first gives the intermediate vertex $x_k$. Taking $w_{k+1}$ first gives $x_{k-1}+w_{k+1}$ and defines the reordered path $\tau_k\gamma$. Set $c_{\tau_k\gamma}=0$ if this path is not allowed in $G$. Deleting the intermediate vertex gives $\eta_k$ with sign $(-1)^k$ in both orders, so
\[
 c_\gamma+c_{\tau_k\gamma}=0.
\]

Suppose $c_\gamma\ne0$. If an axis occurs more than once, choose successive occurrences $w_i,w_j$ on that axis, with $i<j$. Move $w_j$ to position $i+1$, and denote the resulting path by $\gamma'$. Each of the $j-i-1$ interchanges involves distinct axes, so
\[
 c_{\gamma'}=(-1)^{j-i-1}c_\gamma\ne0.
\]
Thus $\gamma'$ is allowed. If $w_j=w_i$, the equal adjacent steps give $c_{\gamma'}=0$. If $w_j=-w_i$, then $\gamma'$ traverses the same edge in both directions, contrary to the definition of $G$. Hence, $w_1,\ldots,w_p$ lie on distinct axes.

Write $b=x_0$ and index the steps as $u_j=\varepsilon_je_{i_j}$ with $i_1<\cdots<i_p$. Since adjacent transpositions generate $S_p$, the coefficients of the permutation paths satisfy
\[
 c_{\gamma_\pi}=\operatorname{sgn}(\pi)c_{\gamma_{\mathrm{id}}}\ne0,
 \quad \pi\in S_p.
\]
In particular, every $\gamma_\pi$ is allowed. To prove that $\sigma=b+\sum_{j=1}^p[0,1]u_j$ belongs to $X_p(G)$, write $b_J=b+\sum_{j\in J}u_j$ for $J\subseteq\{1,\ldots,p\}$. Each edge of the product-oriented cube has the form
\[
 b_J\longrightarrow b_{J\cup\{k\}},\quad k\notin J.
\]
This edge lies on any permutation path that takes the steps indexed by $J$ first and $u_k$ next. Hence all the directed edges of $\sigma$ belong to $G$. Summing the terms of $c$ on the permutation paths of $\sigma$ gives $c_{\gamma_{\mathrm{id}}}\Psi(\sigma)$.

Conversely, fix a cube $\sigma\in X_p(G)$. Two permutations that exchange the steps at positions $k$ and $k+1$ give the same path after deleting the $k$-th vertex. The permutations have opposite signs, so these terms cancel in $\partial\Psi(\sigma)$. Deleting an endpoint gives coefficient $(-1)^{j-1}$ at the target of $u_j$ and $(-1)^j$ at its source. The target face is $\sigma_j^+$ for $\varepsilon_j=1$ and $\sigma_j^-$ for $\varepsilon_j=-1$. Thus, we have
\[
 \partial\Psi(\sigma)=\sum_{j=1}^p(-1)^{j-1}\varepsilon_j
 \bigl(\Psi(\sigma_j^+)-\Psi(\sigma_j^-)\bigr).
\]
Since $\varepsilon_{\sigma_j^\pm}=\varepsilon_\sigma\varepsilon_j$, it follows that
\[
 \partial\Phi(\sigma)=\sum_{j=1}^p(-1)^{j-1}
 \bigl(\Phi(\sigma_j^+)-\Phi(\sigma_j^-)\bigr)
 =\Phi(\partial_{\mathrm{cub}}\sigma).
\]
So one has $\Phi(\sigma)\in\Omega_p(G)$. Distinct cubes have disjoint supports as each permutation path uniquely determines its underlying cube. The coefficient relations now imply
\[
 \Omega_p(G)=\bigoplus_{\sigma\in X_p(G)}\Bbbk\Phi(\sigma).
\]
Thus, $\Phi$ is a chain isomorphism. Since $X(G)$ has dimension at most $n$, $\Omega_p(G)=0$ for $p>n$.

For $G'\subseteq G$ and $\sigma\in X(G')$, the directed edges of $\sigma$ give $\Phi_{G'}(\sigma)=\Phi_G(\sigma)$ in $\Omega_*(G)$. Thus, $\Phi_G|_{C_*^{\mathrm{cub}}(X(G'))}=\Phi_{G'}$ under the natural inclusions. For $T\subseteq V(G)$, we have
\[
 X(G-T)=\{\sigma\in X(G):\operatorname{Vert}(\sigma)\cap T=\varnothing\}.
\]
Hence, $\Phi_G$ restricts to $\Phi_{G-T}:C_*^{\mathrm{cub}}(X(G-T);\Bbbk)\xrightarrow{\cong}\Omega_*(G-T)$.
\end{proof}

\begin{corollary}\label{cor:grid-high-vanishing}
For a lattice digraph $G$ in $\mathbb Z^n$, we have $H_q(G)=0$ for all $q\ge n$.
\end{corollary}
\begin{proof}
Choose integers $a_i<b_i$ such that $X(G)$ is contained inside $\prod_{i=1}^n[a_i,b_i]$, and let $B$ denote the cubical subdivision of this box into unit cubes. Then $X(G)$ is a subcomplex of $B$. Since $|B|$ is contractible and $n\ge1$,
\[
 H_n^{\mathrm{cub}}(B;\Bbbk)=0,\quad C_{n+1}^{\mathrm{cub}}(B;\Bbbk)=0.
\]
Accordingly,
\[
 \begin{aligned}
 \ker\bigl(\partial_n:C_n^{\mathrm{cub}}(B)\to C_{n-1}^{\mathrm{cub}}(B)\bigr)
 &=\operatorname{im}\bigl(\partial_{n+1}:C_{n+1}^{\mathrm{cub}}(B)\to C_n^{\mathrm{cub}}(B)\bigr)\\
 &=0.
 \end{aligned}
\]
The chain inclusion $C_*^{\mathrm{cub}}(X(G))\hookrightarrow C_*^{\mathrm{cub}}(B)$ therefore implies that $X(G)$ admits no nonzero $n$-dimensional cycles. Applying Theorem~\ref{thm:grid-cubical-chains} yields $H_n(G)=0$. For $q>n$, the same theorem gives $\Omega_q(G)=0$, so $H_q(G)=0$.
\end{proof}

This proves the vanishing in \cite[Theorem~3.3]{GrigorLiu} for the lattice digraphs considered here. On such a digraph, deleting a vertex from an allowed path gives a regular face, since the digraph has no pair of opposite edges. Thus, the boundary on allowed chains agrees with the boundary used in that theorem.

\begin{example}\label{ex:grid-cyclic-square}
Take $a=(0,0)$, $b=(1,0)$, $c=(1,1)$, $d=(0,1)$ and edges $a\to b\to c\to d\to a$. These orientations are not a product orientation. Thus, $X(G)$ is the boundary of the square, and Theorem~\ref{thm:grid-cubical-chains} gives $\Omega_2(G)=0$ and $H_1(G)\cong\Bbbk$. The presence of every edge does not suffice to include a cube in $X(G)$.
\end{example}

\begin{example}\label{ex:grid-opposite-edges}
Let $G$ have vertices $0,1\in\mathbb Z$ and edges $0\to1$, $1\to0$. For each $p\ge0$, its two allowed $p$-paths alternate between $0$ and $1$. An interior deletion is nonregular, and an endpoint deletion is allowed. Thus, $\Omega_p(G)=\mathcal A_p(G)\cong\Bbbk^2$ for each $p\ge0$. With opposite edges, neither the one-dimensional cubical chain model nor the vanishing $\Omega_p(G)=0$ for $p>1$ can hold.
\end{example}

\subsection{The cubical link formula}\label{sec:grid-link}

Let $G$ be a lattice digraph in $\mathbb Z^n$, and let $X(G)$ be the complex of cubes with product orientations.

\begin{definition}\label{def:grid-geometric-link}
Let $v\in V(G)$. The \textit{cubical link} $\Lambda_v$ has a vertex $\eta e_i$, with $\eta\in\{\pm1\}$, precisely when either $v\to v+\eta e_i$ or $v+\eta e_i\to v$ is an edge of $G$. A set $\{\eta_1e_{i_1},\ldots,\eta_pe_{i_p}\}$ is a simplex precisely when the axes $i_1,\ldots,i_p$ are distinct and
\[
 v+\sum_{j=1}^p[0,1]\eta_je_{i_j}\in X(G).
\]
The sign $\eta_j$ gives the displacement from $v$, independently of the edge orientation.
\end{definition}

The cubical vertex link $\Lambda_v$ of $X(G)$ at $v$ has at most $2n$ vertices. It forms a subcomplex of the boundary of the $n$-dimensional cross-polytope, and its vertices correspond to geometric edge directions. By contrast, the vertices of $L_v(P(G))=P(L_v)$ are the vertices reachable to or from $v$, excluding $v$ itself. The homology of $P(\Lambda_v)$ need not coincide with that of $L_v(P(G))$ (Example~\ref{ex:cubical-path-links}).

\begin{theorem}\label{thm:grid-geometric-link}
For each $p\ge0$, there is an isomorphism
\begin{equation*}
 Q_p(G,v)\cong\widetilde C_{p-1}(\Lambda_v;\Bbbk),
\end{equation*}
which commutes with the cubical and augmented simplicial boundaries. Taking homology gives
\begin{equation*}
 \LH_q(G,v)\cong\widetilde H_{q-1}(\Lambda_v;\Bbbk),\quad q\ge0.
\end{equation*}
\end{theorem}
\begin{proof}
By Theorem~\ref{thm:grid-cubical-chains}, the classes $[\Phi(\sigma)]$ of the $p$-cubes containing $v$ form a basis of $Q_p(G,v)$. In the relative cubical boundary, each facet disjoint from $v$ vanishes.

For $p\ge1$, write the geometric edge directions of $\sigma$ at $v$ as $d_j=\eta_j e_{i_j}$, where $\eta_j\in\{\pm1\}$ and $i_1<\cdots<i_p$. The correspondence $\sigma\mapsto[d_1,\ldots,d_p]$ identifies the $p$-cubes containing $v$ with the $(p-1)$-simplices of $\Lambda_v$, ordered by coordinate index. Define
\[
 F_p([\Phi(\sigma)])=(-1)^p\Bigl(\prod_{j=1}^p\eta_j\Bigr)[d_1,\ldots,d_p],
 \quad F_0([e_v])=1\in\widetilde C_{-1}(\Lambda_v;\Bbbk).
\]
By definition, each $F_p$ is an isomorphism.

Removing $d_j$ from $[d_1,\dots,d_p]$ yields a face with boundary coefficient $(-1)^{j-1}$. The corresponding facet of $\sigma$ is $\sigma_j^-$ if $\eta_j=1$ and $\sigma_j^+$ if $\eta_j=-1$. The coefficient of this facet in $\partial_{\mathrm{cub}}\sigma$ is $(-1)^j\eta_j$. Note that
\[
 (-1)^j\eta_j\,(-1)^{p-1}\prod_{k\ne j}\eta_k
 =(-1)^{j-1}\,(-1)^p\prod_{k=1}^p\eta_k.
\]
We have $F_{p-1}\partial=\partial F_p$ for $p\ge1$. Passing to homology gives the desired isomorphism.
\end{proof}

\begin{corollary}\label{cor:grid-geometric-realization}
Let $G$ be a lattice digraph, and let $v$ be a vertex of $G$. Then there is an isomorphism
\begin{equation*}
 \LH_\ast(G,v)\cong H_\ast\bigl(|X(G)|,|X(G)|\setminus\{v\};\Bbbk\bigr).
\end{equation*}
\end{corollary}
\begin{proof}
Let $X=|X(G)|$ and $Y=|X(G-v)|$. Suppose $\sigma$ is a $p$-cube with $v$ as its origin. For $x=(t_1,\dots,t_p)\ne 0$, the ray starting at $v$ and passing through $x$ intersects the face opposite $v$ at
\[
 r_\sigma(x)=\frac{x}{\max_{1\le i\le p}t_i}.
\]
Define the homotopy that moves $x$ along the segment from $x$ to $r_\sigma(x)$ by
\[
 H_\tau(x)=(1-\tau)x+\tau r_\sigma(x),\quad 0\le\tau\le1.
\]
The map $r_\sigma$ retracts $\sigma\setminus\{v\}$ onto the union of the faces opposite $v$. These homotopies agree on shared faces containing $v$, since the omitted coordinates vanish. Extending by the identity on cubes disjoint from $v$ yields a deformation retraction of $X\setminus\{v\}$ onto $Y$.

By Theorem~\ref{thm:grid-cubical-chains}, one obtains $\LH_\ast(G,v)\cong H_\ast(X,Y;\Bbbk)$. The deformation retraction gives $H_\ast(X,Y;\Bbbk)\cong H_\ast(X,X\setminus\{v\};\Bbbk)$, which completes the proof.
\end{proof}

\begin{remark}\label{rem:grid-multiple-punctures}
Let $C_2=I\Box I$ be a directed square, and let $T$ be the vertex set of two opposite vertices of $C_2$. Write $X=X(C_2)$ for convenience. The complex $X(C_2-T)$ consists of two isolated vertices, while $|X|\setminus T$ is path connected. Since $|X|$ is contractible, the long exact sequences of the two pairs give
\[
 H_1(X,X(C_2-T);\Bbbk)\cong\Bbbk,\quad H_1(|X|,|X|\setminus T;\Bbbk)=0.
\]
\end{remark}
Listing the vertices of each nonempty simplex in increasing order gives the allowed paths of $P(\Lambda_v)$. By Theorem~\ref{thm:grid-geometric-link}, one obtains
\[
 \LH_q(G,v)\cong\widetilde H_{q-1}\bigl(P(\Lambda_v)\bigr),\quad q\ge0.
\]

\begin{example}\label{ex:cubical-path-links}
Realize the digraph $G$ of Example~\ref{ex:diamond} as the directed unit square with $a=(0,0)$, $v=(1,0)$, $w=(0,1)$, and $b=(1,1)$. The cubical link $\Lambda_v$ is the edge with endpoints $-e_1$ and $e_2$. The link path complex $L_v(P(G))$ consists of the vertices $a,b$. Their reduced homology groups in degree zero are
\[
 \widetilde H_0(P(\Lambda_v))=0,\quad
 \widetilde H_0(L_v(P(G)))\cong\Bbbk.
\]
The cubical link formula gives $\LH_1(G,v)=0$. Example~\ref{ex:diamond} gives the same result on path chains. The cycle $e_{av}+e_{vb}$ represents a generator of $H_1(S_v,L_v)$, but in $Q_*(G,v)$ it satisfies $\partial[e_{avb}-e_{awb}]=[e_{av}+e_{vb}]$.
\end{example}

By Theorem~\ref{thm:grid-geometric-link}, $v$ is regular of dimension $d$ if and only if $\widetilde H_{d-1}(\Lambda_v;\Bbbk)\cong\Bbbk$ and $\widetilde H_j(\Lambda_v;\Bbbk)=0$ for each $j\ge-1$ with $j\ne d-1$.

\section{Matching Gluing}\label{sec:matching-gluing}

We follow the matching construction of \cite{GrigorLiu}. The interface consists of the edges present in both layers between matched vertices. The prism chains give a mapping cone for the glued graph and, after vertex deletion, for its local chain complex.

\subsection{The mapping-cone model}

Let $A$ and $B$ be digraphs with disjoint vertex sets, and let $M$ be a finite set. Suppose that $j_A:x\mapsto a_x$ and $j_B:x\mapsto b_x$ are injections from $M$ into $V(A)$ and $V(B)$, respectively. Let $G$ consist of $A,B$ and the interlayer edges
\[
 a_x\longrightarrow b_x,\quad x\in M,
\]
and no other edges between the two layers. The endpoints of these edges are pairwise disjoint. Directed cycles within either layer are permitted.

\begin{definition}\label{def:matching-interface}
The \textit{common interface} is the digraph $K$ with vertex set $M$ and edge relation defined by
\begin{equation*}
 x\longrightarrow y\text{ in }K
 \quad\Longleftrightarrow\quad
 a_x\longrightarrow a_y\text{ in }A\text{ and }
 b_x\longrightarrow b_y\text{ in }B.
\end{equation*}
It has natural embeddings $j_A,j_B$ into the two sides, whose images need not be induced subgraphs.
\end{definition}

The image $j_A(K)$ is induced on the matched vertices of $A$ precisely when each edge $a_x\to a_y$ between them has a corresponding edge $b_x\to b_y$ in $B$. The same criterion holds with $A$ and $B$ interchanged.

For a digraph $Y$, let $\Lambda_p(V(Y))$ denote the vector space with basis the elementary $p$-paths on $V(Y)$. Its alternating boundary $\partial^\Lambda$ retains all faces, and $\partial^\Lambda_0=0$. We write $\mathcal A_p(Y)$ for the space generated by the allowed $p$-paths of $Y$.

\begin{definition}\label{def:matching-nonregular}
Following \cite[Section~2]{GrigorLiu}, define
\[
 \Omega_0^\Lambda(Y)=\mathcal A_0(Y),\qquad
 \Omega_p^\Lambda(Y)=\{c\in\mathcal A_p(Y):\partial^\Lambda c\in\mathcal A_{p-1}(Y)\},\quad p\ge1.
\]
The homology of this complex is denoted by $H_q^\Lambda(Y)$. Throughout this subsection, $\Omega_*(Y)$ denotes the complex with the regular boundary of Definition~\ref{def:regular-boundary}.
\end{definition}

For an elementary path $c=e_{x_0\cdots x_q}$ on $M$, define
\[
 P_i(c)=e_{a_{x_0}\cdots a_{x_i}b_{x_i}\cdots b_{x_q}},\quad 0\le i\le q.
\]

\begin{definition}\label{def:matching-prism}
The \textit{prism operator} is the linear map $h:\Lambda_q(M)\to\Lambda_{q+1}(V(G))$ determined by
\[
 h(c)=\sum_{i=0}^q(-1)^iP_i(c).
\]
\end{definition}

\begin{proposition}\label{prop:matching-prism-boundary}
For either boundary $d=\partial^\Lambda$ or $d=\partial$, the prism operator satisfies
\[
 d h(c)+h(d c)=j_Bc-j_Ac.
\]
In particular, $h$ restricts to maps $\Omega_q^\Lambda(K)\to\Omega_{q+1}^\Lambda(G)$ and $\Omega_q(K)\to\Omega_{q+1}(G)$ for each $q\ge0$.
\end{proposition}
\begin{proof}
First take $d=\partial^\Lambda$. Let $d_j$ delete the vertex at position $j$, counting from zero. In $d h(c)$, the faces $d_iP_i(c)=d_iP_{i-1}(c)$ occur with opposite signs for $1\le i\le q$. The endpoint terms $d_0P_0(c)$ and $d_{q+1}P_q(c)$ contribute $j_Bc$ and $-j_Ac$. For $q\ge1$, the remaining faces satisfy
\[
 d_jP_i(c)=
 \begin{cases}
 P_{i-1}(d_jc),&j<i,\\
 P_i(d_{j-1}c),&j>i+1.
 \end{cases}
\]
With their boundary signs, these terms sum to $-h(dc)$. For $q=0$, the identity is $d e_{a_{x_0}b_{x_0}}=e_{b_{x_0}}-e_{a_{x_0}}$.

If $c$ has a consecutive repetition, each $P_i(c)$ also has one. Thus the identity descends to regular chains. If $c$ is an allowed path in $K$, each $P_i(c)$ is allowed in $G$ by the definition of the interface and the matching edges $a_x\to b_x$. Hence, whenever $c$ and $dc$ are allowed chains, both $h(c)$ and $d h(c)=j_Bc-j_Ac-h(dc)$ are allowed.
\end{proof}

\begin{definition}
Let $f:C_*\to D_*$ be a chain map. The \textit{mapping cone} of $f$ is the chain complex specified by
\[
 \operatorname{Cone}(f)_p=D_p\oplus C_{p-1},\quad
 d(u,c)=(\partial u+f(c),-\partial c).
\]
For ordinary chain complexes, chain groups in negative degrees vanish.
\end{definition}

\begin{theorem}\label{thm:matching-cone}
Let $f=(-j_A,j_B):\Omega_*(K)\to\Omega_*(A)\oplus\Omega_*(B)$. There is a chain isomorphism
\begin{equation*}
 \Phi: \operatorname{Cone}(f) \to \Omega_*(G)
\end{equation*}
natural under induced restrictions of the matching structure. Explicitly, the graded vector space decomposes as
\begin{equation*}
 \Omega_p(G)=\Omega_p(A)\oplus\Omega_p(B)\oplus h\Omega_{p-1}(K),
\end{equation*}
and the chain isomorphism sends $(u_A,u_B,c)$ to $u_A+u_B+h(c)$. The same formula gives a chain isomorphism for $\Omega_*^\Lambda$.
\end{theorem}
\begin{proof}
For $p=0$, the decomposition follows from $V(G)=V(A)\sqcup V(B)$. Fix $p\ge1$ and let $d$ be either $\partial$ or $\partial^\Lambda$. Take $z\in\mathcal A_p(G)$ with $dz\in\mathcal A_{p-1}(G)$. Write $[w](z)$ for the coefficient of the elementary path $w$ in $z$.

For each $x\in M$, give both $a_x$ and $b_x$ the label $x$. The unmatched vertex outside $M$ is assigned its own label. Then a path can be labelled as $\xi=x_0\cdots x_{p-1}$ whose consecutive labels are distinct. For $0\le i<p$, let $\lambda_i(\xi)$ be the coefficient in $z$ of the elementary path
\[
 P_i(\xi)=e_{a_{x_0}\cdots a_{x_i} b_{x_i}\cdots b_{x_{p-1}}}.
\]
Here, $a_x$ and $b_x$ are the vertices labelled $x$ in the respective layers. Define $\lambda_i(\xi)=0$ whenever there is no required path.

Suppose at least one of $P_i(\xi)$ and $P_{i+1}(\xi)$ is allowed for some $0\le i<p-1$. The face
\[
 F_i(\xi)=e_{a_{x_0}\cdots a_{x_i} b_{x_{i+1}}\cdots b_{x_{p-1}}}
\]
is regular and contains the forbidden step $a_{x_i}\to b_{x_{i+1}}$. Any allowed $p$-path with $F_i(\xi)$ as a face must have an additional vertex between $a_{x_i}$ and $b_{x_{i+1}}$. The matching edges give two possible insertions. Inserting $b_{x_i}$ gives $P_i(\xi)$, and inserting $a_{x_{i+1}}$ gives $P_{i+1}(\xi)$. In both cases the inserted vertex has position $i+1$, so
\[
 0=[F_i(\xi)](dz)=(-1)^{i+1}\bigl(\lambda_i(\xi)+\lambda_{i+1}(\xi)\bigr).
\]
If neither $P_i(\xi)$ nor $P_{i+1}(\xi)$ is allowed, both coefficients are zero, so the same relation holds. Iterating this relation from index $0$ gives
\[
 \lambda_i(\xi)=(-1)^i\lambda_0(\xi),\qquad 0\le i<p.
\]
When $\lambda_0(\xi)\ne0$, both $P_0(\xi)$ and $P_{p-1}(\xi)$ are allowed. Their horizontal steps give the entire path $\xi$ in $B$ and $A$, respectively, so $\xi$ is allowed in $K$. Hence, the mixed part of $z$ is $h(c)$ for the unique chain
\[
 c=\sum_{\xi}\lambda_0(\xi)e_\xi\in\mathcal A_{p-1}(K),
\]
where the sum ranges over the allowed $(p-1)$-paths of $K$.

For $p\ge2$, put $r=p-2$. Write $z=u_A+u_B+h(c)$, where $u_A$ and $u_B$ consist of the paths lying entirely in $A$ and $B$, respectively. By Proposition~\ref{prop:matching-prism-boundary}, we have
\[
 dz=(du_A-j_Ac)+(du_B+j_Bc)-h(dc).
\]
Thus, the paths in $dz$ that meet the layers $A$ and $B$ come from $-h(dc)$.

For $\eta=y_0\cdots y_r$ in the support of $dc$, consider
\[
 P_i(\eta)=e_{a_{y_0}\cdots a_{y_i}b_{y_i}\cdots b_{y_r}}.
\]
If $P_j(\eta')=P_i(\eta)$, we have $j=i$ and $\eta'=\eta$. So the coefficient of this path in $h(dc)$ is $(-1)^i[\eta](dc)$. Since the mixed part of $dz$ is $-h(dc)$, we obtain
\[
 [P_i(\eta)](dz)=(-1)^{i+1}[\eta](dc),\quad 0\le i\le r.
\]
Suppose $\eta$ is not allowed in $K$. If one of its steps is missing in $A$, take $i=r$; then the corresponding prism path $P_r(\eta)$ contains this missing $A$-step. If one of its steps is missing in $B$, take $i=0$; then $P_0(\eta)$ contains this missing $B$-step. For $d=\partial^\Lambda$, a consecutive repetition in $\eta$ also makes $P_0(\eta)$ nonregular. Thus, in every case, the specific path $P_i(\eta)$ is non-allowed. Because $dz$ is an allowed chain, we get $[P_i(\eta)](dz)=0$. The coefficient formula
\[
 [P_i(\eta)](dz)=(-1)^{i+1}[\eta](dc)
\]
implies $[\eta](dc)=0$. Hence, every non-allowed path has zero coefficient in $dc$, and $dc\in\mathcal A_{p-2}(K)$. For $p=1$, $c$ has degree zero and $dc=0$.

By Proposition~\ref{prop:matching-prism-boundary}, $d h(c)\in\mathcal A_{p-1}(G)$. Write
\[
 z-h(c)=u_A+u_B,
\]
where $u_A$ and $u_B$ are supported entirely in $A$ and $B$, respectively. Then
\[
 du_A+du_B=dz-dh(c)\in\mathcal A_{p-1}(G).
\]
The supports of $du_A$ and $du_B$ lie in disjoint layers, so each is allowed. Thus, one has $u_A\in\Omega_p(A)$ and $u_B\in\Omega_p(B)$. The cone differential is
\[
 (u_A,u_B,c)\longmapsto (du_A-j_Ac,\ du_B+j_Bc,\ -dc).
\]
For an induced restriction, the inclusions preserve each term $P_i(\xi)$ of every retained prism, so they commute with $\Phi$.
\end{proof}

\begin{corollary}\label{cor:matching-les}
There exists a long exact sequence
\begin{equation*}
 \begin{aligned}
 \cdots\longrightarrow H_q(K)
 &\xrightarrow{\alpha_q}H_q(A)\oplus H_q(B)
 \longrightarrow H_q(G)\\
 &\longrightarrow H_{q-1}(K)
 \xrightarrow{\alpha_{q-1}}H_{q-1}(A)\oplus H_{q-1}(B)
 \longrightarrow\cdots ,
 \end{aligned}
\end{equation*}
where $\alpha_q=(-j_{A*},j_{B*})$. The sequence ends with $H_0(G)\to0$.
\end{corollary}
\begin{proof}
Set $D_*=\Omega_*(A)\oplus\Omega_*(B)$ and $C_*=\Omega_*(K)$. The cone fits into a short exact sequence $0\to D_*\to\operatorname{Cone}(f)\to C_{*-1}\to0$, where the last complex carries differential $-\partial$. For a cycle $c\in C_q$, the lift $(0,c)\in\operatorname{Cone}(f)_{q+1}$ has boundary $(f(c),0)$. The connecting homomorphism therefore maps $[c]$ to $[f(c)]=(-j_{A*}[c],j_{B*}[c])$. Combining the associated homology long exact sequence with Theorem~\ref{thm:matching-cone} completes the proof. The same argument establishes the corresponding sequence for $H_*^\Lambda$.
\end{proof}

\begin{example}\label{ex:matching-interface-choice}
Let $Z=Z_4$ be the four-vertex directed acyclic graph
\[
 0\longrightarrow1\longrightarrow2\longrightarrow3,
 \quad 0\longrightarrow3.
\]
Every allowed path of length at least two lies on the long branch and has a non-allowed interior face with no other allowed filling. Hence, we have
\[
 \Omega_p(Z)=0,\quad p\ge2,\quad
 H_0(Z)\cong H_1(Z)\cong\Bbbk,
\]
where $H_1(Z)$ is generated by $[z]$, with $z=e_{01}+e_{12}+e_{23}-e_{03}$.

Let $A=Z\Box Z$, and let $B$ consist of $16$ isolated vertices. Join each vertex of $A$ to its copy in $B$ by a matching edge. The common interface $K$ has $16$ vertices and no edges, so $H_q(K)=H_q(B)=0$ for $q\ge1$. Corollary~\ref{cor:matching-les} yields an isomorphism induced by the inclusion $A\hookrightarrow G$:
\[
 H_2(A)\xrightarrow{\cong}H_2(G),\quad
 H_2(A)\cong H_1(Z)\otimes_{\Bbbk}H_1(Z)\cong\Bbbk.
\]
The Cartesian product K\"unneth formula gives the last isomorphism. The matching bijection $A\to B$ sends edges to pairs of distinct isolated vertices, so it induces no interface chain map from $A$ to $B$. Replacing $K$ by $A$ would also give the identity as the $A$-component in degree two. Since $H_2(B)=0$, exactness would then force $H_2(A)\to H_2(G)$ to be zero, although inclusion gives the isomorphism above.
\end{example}

Theorem~\ref{thm:matching-cone} extends to the setting where matching edges may carry either orientation.

Let $A$ and $B$ be digraphs with disjoint vertex sets, let $M=M_+\sqcup M_-$ be a finite set, and let $j_A:x\mapsto a_x$ and $j_B:x\mapsto b_x$ be injections into $V(A)$ and $V(B)$, respectively. Define $G$ by
\[
 \begin{aligned}
 V(G)&=V(A)\sqcup V(B),\\
 E(G)&=E(A)\cup E(B)
 \cup\{(a_x,b_x):x\in M_+\}
 \cup\{(b_x,a_x):x\in M_-\}.
 \end{aligned}
\]
Thus, the matching edges have pairwise disjoint endpoints; $M_+$ indexes edges directed from $A$ to $B$, and $M_-$ indexes edges directed from $B$ to $A$. For $\varepsilon\in\{+,-\}$, let $K_\varepsilon$ denote the digraph on vertex set $M_\varepsilon$ with
\[
 x\to y\text{ in }K_\varepsilon
 \quad\Longleftrightarrow\quad
 a_x\to a_y\text{ in }A\text{ and }b_x\to b_y\text{ in }B.
\]
Set $K=K_+\sqcup K_-$, with no edges between $M_+$ and $M_-$. The embeddings into $A$ and $B$ induce chain maps, also denoted $j_A$ and $j_B$. We write $c=(c_+,c_-)\in\Omega_p(K)$, where $c_\varepsilon\in\Omega_p(K_\varepsilon)$.

\begin{theorem}\label{thm:matching-mixed-directions}
With the notation above, define $f:\Omega_*(K)\to\Omega_*(A)\oplus\Omega_*(B)$ by
\[
 f(c_+,c_-)=\bigl(-j_Ac_++j_Ac_-,\ j_Bc_+-j_Bc_-\bigr).
\]
Then there is a chain isomorphism
\[
 \Phi:\operatorname{Cone}(f)\xrightarrow{\cong}\Omega_*(G)
\]
natural under induced restrictions of the matching.
\end{theorem}
\begin{proof}
Let $G_\varepsilon$ consist of $A$, $B$, and the matching edges indexed by $M_\varepsilon$. We first prove that every path in the support of $c\in\Omega_p(G)$ crosses between the layers at most once. Suppose otherwise, and choose a path $w$ with nonzero coefficient for which the number of edges between its first two crossings is minimal. After interchanging $A$ and $B$ if necessary, the relevant segment of $w$ has the form
\[
 a_x\to b_x=y_0\to y_1\to\cdots\to y_m=b_z\to a_z,
 \qquad x\in M_+,\ z\in M_-.
\]
Since the matching endpoints are disjoint, $x\ne z$, and hence $m\ge1$.

Let $F$ be the face obtained from $w$ by deleting $b_x$, which occurs at position $k$. The resulting face contains the step $a_x\to y_1$, which is absent from $G$. Therefore
\[
 [F](\partial c)=0.
\]
Any allowed path having $F$ as a face must insert a vertex $u$ at position $k$ with
\[
 a_x\to u\to y_1.
\]
If $u\in B$, the only possible matching edge from $a_x$ gives $u=b_x$, so the path is $w$. If $u\in A$, the matching structure gives $u=a_y$ and $y_1=b_y$ for a unique $y\in M_+$. Let $w'$ denote this second path, and set $[w'](c)=0$ when it does not exist. Both $w$ and $w'$ delete their vertex at position $k$, so
\[
 0=[F](\partial c)=(-1)^k\bigl([w](c)+[w'](c)\bigr).
\]
Since $[w](c)\ne0$, the path $w'$ exists and
\[
 [w'](c)=-[w](c)\ne0.
\]

If $m=1$, then $b_y=y_1=b_z$, so $y=z$, contradicting $M_+\cap M_-=\varnothing$. If $m>1$, the first two crossings of $w'$ are $a_y\to b_y$ and $b_z\to a_z$, and there are only $m-1$ edges between them. This contradicts the choice of $w$. Thus every path with nonzero coefficient in $c$ crosses at most once.

Separate the terms of $c$ into three classes: paths contained in $A$, paths contained in $B$, and paths crossing from $A$ to $B$ or from $B$ to $A$. A face of a mixed path that remains in one layer is obtained by deleting an endpoint and is allowed. Every non-allowed face of a mixed path therefore still contains vertices from both layers in the same order as the original path. Consequently, its coefficient can only receive contributions from paths in the same class. The boundary condition for $c$ thus holds separately for the pure $A$-terms, the pure $B$-terms, and the two classes of mixed terms. Therefore
\[
 \Omega_*(G)=\Omega_*(G_+)+\Omega_*(G_-),\qquad
 \Omega_*(G_+)\cap\Omega_*(G_-)=D_*,
 \quad D_*=\Omega_*(A)\oplus\Omega_*(B).
\]

Let $h_+$ be the prism map from $A$ to $B$ over $K_+$, and let $h_-$ be the prism map from $B$ to $A$ over $K_-$. Applying Theorem~\ref{thm:matching-cone} in the two directions gives
\[
 \Omega_p(G)=D_p\oplus h_+\Omega_{p-1}(K_+)\oplus h_-\Omega_{p-1}(K_-).
\]
The identities
\[
 \partial h_++h_+\partial=j_B-j_A,\qquad
 \partial h_-+h_-\partial=j_A-j_B
\]
give the two signs in $f$. Hence
\[
 \Phi(u,c_+,c_-)=u+h_+(c_+)+h_-(c_-)
\]
is the required chain isomorphism. It is natural under induced restrictions because every retained prism has the same elementary-path terms before and after restriction.

For $c\in\Omega_p^\Lambda(G)$, passing to regular chains does not change the coefficients of the allowed paths. The argument above therefore shows again that every path with nonzero coefficient crosses at most once. A non-allowed face of a mixed path still records the same order of the two layers, even when it contains a consecutive repetition. The same separation by path class applies to $\partial^\Lambda$, and Theorem~\ref{thm:matching-cone} gives the same cone formula.

\end{proof}

\begin{remark}
For $\varepsilon\in\{+,-\}$, let $S_\varepsilon=A[\{a_x:x\in M_\varepsilon\}]$. Suppose that for each edge $a_x\to a_y$ in $S_\varepsilon$, the corresponding edge $b_x\to b_y$ belongs to $B$. Then $j_A$ identifies $K_\varepsilon$ with $S_\varepsilon$. With these identifications and the signs in $f$, the long exact cone sequence for $H_*^\Lambda$ agrees with the sequence in \cite[Theorem~3.1]{GrigorLiu} through the map
\[
 H_2^\Lambda(G)\longrightarrow H_1^\Lambda(S_+)\oplus H_1^\Lambda(S_-).
\]
The cone sequence continues to $H_0^\Lambda(G)\to0$. If $G$ has no pair of opposite edges, deleting an interior vertex of an allowed path cannot produce a consecutive repetition. Hence the regular and nonregular boundaries agree on allowed chains of $G$, its sublayers, and its interfaces.
\end{remark}

\subsection{Local gluing}

Assume that all matching edges share the same direction. For mixed orientations, adopt the sign conventions from Theorem~\ref{thm:matching-mixed-directions}. Given $T\subseteq V(G)$, define
\begin{equation*}
 T_A=T\cap V(A),\quad T_B=T\cap V(B),\quad
 U_T=\{x\in M:a_x\in T_A\text{ or }b_x\in T_B\}.
\end{equation*}
We define the local quotient corresponding to the empty marked set to be zero. Removing $T$ yields the common interface $K-U_T$.

\begin{theorem}\label{thm:matching-local-cone}
The interface embeddings induce a natural chain map
\[
 f_T:Q_*(K,U_T)\longrightarrow Q_*(A,T_A)\oplus Q_*(B,T_B),
\]
and there exists a natural chain isomorphism
\begin{equation*}
 Q_*(G,T)\cong\operatorname{Cone}(f_T).
\end{equation*}
As a consequence, we obtain the long exact sequence
\begin{equation*}
 \begin{aligned}
 \cdots\longrightarrow\LH_q(K,U_T)
 &\xrightarrow{\alpha_q^T}\LH_q(A,T_A)\oplus\LH_q(B,T_B)\\
 &\longrightarrow\LH_q(G,T)
 \longrightarrow\LH_{q-1}(K,U_T)\longrightarrow\cdots
 \end{aligned}
\end{equation*}
where $\alpha_q^T=H_q(f_T)$, and all terms in negative degrees vanish.
\end{theorem}
\begin{proof}
The embeddings send $\Omega_*(K-U_T)$ into $\Omega_*(A-T_A)$ and $\Omega_*(B-T_B)$, so $f$ descends to
\[
 f_T([c])=(-[j_Ac],[j_Bc]).
\]
Under Theorem~\ref{thm:matching-cone}, the inclusion $G-T\subseteq G$ becomes the inclusion of the cone of
\[
 \Omega_*(K-U_T)\longrightarrow\Omega_*(A-T_A)\oplus\Omega_*(B-T_B)
\]
into the cone for $f$. Since the cone is a direct sum in each degree, its quotient has terms
\[
 Q_p(A,T_A)\oplus Q_p(B,T_B)\oplus Q_{p-1}(K,U_T)
\]
and differential $(u,c)\mapsto(\partial u+f_T(c),-\partial c)$. These are the groups and differential of $\operatorname{Cone}(f_T)$. In the cone sequence, the connecting map sends a cycle class $[c]$ to $[f_T(c)]$, giving the displayed map on local homology.
\end{proof}

A matching index $x$ is removed from the interface as soon as $a_x\in T$ or $b_x\in T$.

\begin{corollary}\label{cor:matching-point-local}
If $v=a_x$ is a matched vertex on side $A$, then
\begin{equation*}
 Q_*(G,v)\cong\operatorname{Cone}\bigl(Q_*(K,x)\xrightarrow{-j_A}Q_*(A,v)\bigr).
\end{equation*}
Thus, the inclusion $G-v\hookrightarrow G$ induces isomorphisms $H_q(G-v)\to H_q(G)$ for all $q\ge0$ if and only if $-j_A:Q_*(K,x)\to Q_*(A,v)$ is a quasi-isomorphism. If $v\in V(A)$ is unmatched, then $Q_*(G,v)\cong Q_*(A,v)$. The analogous statements hold on side $B$.
\end{corollary}
\begin{proof}
In Theorem~\ref{thm:matching-local-cone}, take $T=\{a_x\}$. Here $T_B=\varnothing$ and $U_T=\{x\}$, so $Q_*(B,T_B)=0$ and the interface term is $Q_*(K,x)$, so the local complex is $\operatorname{Cone}(-j_A)$. By exactness for $(G,G-a_x)$, deletion preserves homology exactly when this cone is acyclic. By the long exact sequence of the cone, acyclicity is equivalent to $-j_A:Q_*(K,x)\to Q_*(A,a_x)$ inducing a homology isomorphism in each degree. For an unmatched vertex $v$ in $A$, the interface term also vanishes, leaving $Q_*(A,v)$. For $v=b_x$, the surviving component of $f_T$ is $+j_B$.
\end{proof}

\begin{corollary}\label{cor:matching-layer-excision}
Identify $K$ with its image under $j_A$ inside $A$. There exist natural isomorphisms
\begin{equation*}
 H_q(G,B)\cong H_q(A,K),\quad q\ge0.
\end{equation*}
Therefore, the inclusion $B\hookrightarrow G$ induces isomorphisms $H_q(B)\to H_q(G)$ for all $q\ge0$ if and only if the inclusion $K\hookrightarrow A$ induces isomorphisms $H_q(K)\to H_q(A)$ for all $q\ge0$.
\end{corollary}
\begin{proof}
Take $T=V(A)$ in Theorem~\ref{thm:matching-local-cone}. The map $\operatorname{Cone}(-j_A)\to\Omega_*(A,K)$, $(u,c)\mapsto[u]$, is surjective. On the kernel, $(j_Aa,c)\mapsto(a,c)$ carries the differential to $(a,c)\mapsto(\partial a-c,-\partial c)$. This identifies the kernel with the contractible complex $\operatorname{Cone}(-\mathrm{id}_{\Omega_*(K)})$. The projection therefore induces $H_q(G,B)\cong H_q(A,K)$. The long exact sequences of $(G,B)$ and $(A,K)$ give the equivalence of the two inclusion criteria.
\end{proof}
The relative group $H_q(A,K)$ uses the subdigraph $j_A(K)$, which need not be induced. The formula applies to the whole layer $A$ without a separation condition.

\subsection{Criteria for short exact sequences}\label{sec:matching-short-exact}

The long exact sequence from Corollary~\ref{cor:matching-les} yields
\begin{equation*}
 0\longrightarrow\operatorname{coker}\alpha_q
 \longrightarrow H_q(G)
 \longrightarrow\ker\alpha_{q-1}\longrightarrow0.
\end{equation*}
The natural maps assemble into a short exact sequence
\begin{equation*}
 0\longrightarrow H_q(K)\xrightarrow{\alpha_q}
 H_q(A)\oplus H_q(B)\longrightarrow H_q(G)\longrightarrow0,
\end{equation*}
if and only if both $\alpha_q$ and $\alpha_{q-1}$ are injective. Even when all matching edges point from $A$ to $B$, $\alpha_{q-1}$ can have nonzero kernel.

The next example illustrates why injectivity of $\alpha_{q-1}$ is necessary, even when the common interface is induced in both layers.

\begin{example}\label{ex:matching-ten-vertex}
Let $Z=Z_4$ denote the digraph from Example~\ref{ex:matching-interface-choice}. Take disjoint copies $Z^-,Z^+$ with vertices $i^-,i^+$, and adjoin new vertices $s,t$. Construct $A$ by adding edges $s\to i^-$ to $Z^-$, and $B$ by adding edges $i^+\to t$ to $Z^+$, for $0\le i\le3$. Introduce matching edges $i^-\to i^+$ and no further edges. The resulting digraph $G$ is acyclic and has ten vertices. Its common interface is $Z$, which is induced on the matched vertices within each layer.

The layers $A$ and $B$ form a source cone and a sink cone, respectively, so their reduced homology vanishes in every degree. Consequently,
\[
 \alpha_1:H_1(Z)\longrightarrow H_1(A)\oplus H_1(B)
 \quad\text{is the map}\quad\Bbbk\longrightarrow0.
\]
Corollary~\ref{cor:matching-les} gives
\[
 H_2(G)\cong\ker\alpha_1\cong\Bbbk,\quad
 H_1(G)=0,\quad H_0(G)\cong\Bbbk.
\]
Here $\alpha_2:0\to0$ is injective. The cokernel of $H_2(A)\oplus H_2(B)\to H_2(G)$ equals $\ker\alpha_1\cong\Bbbk$.

The nonzero degree-two homology class admits an explicit representative. Let $z^-,z^+$ denote the layerwise copies of the cycle $z=e_{01}+e_{12}+e_{23}-e_{03}$, and let $sz^-$ and $z^+t$ stand for termwise prefixing and suffixing by the cone vertices. Then we have
\[
 c=sz^-+h(z)-z^+t\in\Omega_2(G),\quad
 \partial c=z^-+(z^+-z^-)-z^+=0.
\]
Each path in $sz^-$ starts at $s$. No path in $h(z)$ or $z^+t$ contains $s$, so every path in $sz^-$ retains its coefficient in $c$. Since $sz^-\ne0$, we have $c\ne0$. Since $\Omega_p(Z)=0$ for $p\ge2$, the decomposition in Theorem~\ref{thm:matching-cone} gives $\Omega_3(G)=0$. Thus, $0\ne[c]\in H_2(G)$.
\end{example}

\begin{remark}\label{rem:matching-source-version}
In \cite[Theorem~3.2]{GrigorLiu}, the interface is the subdigraph $S$ of $A$ induced on the sources of the matching edges. The stated sequence for $q\ge2$ reads
\[
 0\longrightarrow H_q^\Lambda(S)\longrightarrow
 H_q^\Lambda(A)\oplus H_q^\Lambda(B)\longrightarrow H_q^\Lambda(G)\longrightarrow0.
\]
In Example~\ref{ex:matching-ten-vertex}, the edges $i^-\to i^+$ match $Z^-$ to $Z^+$, and $S$ and $K$ both identify with $Z_4$. Since $G$ is acyclic, $\partial^\Lambda$ and $\partial$ agree on its allowed paths. In $A$ we have $\partial^\Lambda(sz^-)=z^-$, and in $B$ we have $\partial^\Lambda(z^+t)=z^+$. Thus, $\alpha_1([z])=(-[z^-],[z^+])=0$. Since $[z]$ generates $H_1^\Lambda(K)$ and $H_2^\Lambda(A)=H_2^\Lambda(B)=0$, the connecting homomorphism gives $H_2^\Lambda(G)\cong\ker\alpha_1\cong\Bbbk$. In degree two the sequence of \cite[Theorem~3.2]{GrigorLiu} would therefore require a surjection $0\to\Bbbk$.
\end{remark}

The proof of \cite[Theorem~1.1]{GrigorLiu} uses the short exact sequence of \cite[Theorem~3.2]{GrigorLiu}. For the induction on vertical layers, we need the higher homology of a layer $A$ and the kernel of $H_1(K)\to H_1(A)$ for a subdigraph $K\subseteq A$.

\begin{lemma}\label{lem:outdegree-one}
Let $A$ be a digraph where every vertex has outdegree at most one, and let $K\subseteq A$ be a subdigraph. Then
\[
 \Omega_q^\Lambda(A)=\Omega_q^\Lambda(K)=0,\quad H_q(A)=H_q(K)=0,
 \quad q\ge2.
\]
The inclusion $K\hookrightarrow A$ induces an injection $H_1(K)\hookrightarrow H_1(A)$.
\end{lemma}
\begin{proof}
The claim for $\Omega_*^\Lambda$ follows from \cite[Lemma~3.5]{GrigorLiu} applied to $A$ and $K$. For the regular boundary, let $\gamma=x_0\cdots x_p$ have coefficient $c_\gamma\ne0$ in $c\in\Omega_p(A)$, with $p\ge2$. Suppose $x_{i-1}\ne x_{i+1}$ for some $1\le i<p$. Deleting $x_i$ gives a regular face with a missing edge $x_{i-1}\to x_{i+1}$, since the only outgoing edge at $x_{i-1}$ leads to $x_i$. To obtain this face from an allowed $p$-path, one must insert a vertex between $x_{i-1}$ and $x_{i+1}$, and that vertex must be $x_i$. The face therefore has coefficient $(-1)^i c_\gamma\ne0$ in $\partial c$, contradicting $\partial c\in\mathcal A_{p-1}(A)$. Hence, $x_{i-1}=x_{i+1}$ at each interior index, and $\gamma$ alternates along a pair $a\rightleftarrows b$.

Distinct pairs $a\rightleftarrows b$ have disjoint vertex sets by the outdegree bound. For each pair, let $u_p$ and $v_p$ denote the alternating $p$-paths starting at $a$ and $b$, respectively. Every interior deletion gives consecutive equal vertices, so the regular boundary is
\[
 \partial u_p=v_{p-1}+(-1)^p u_{p-1},\quad
 \partial v_p=u_{p-1}+(-1)^p v_{p-1}.
\]
Thus, the alternating paths form a basis of $\Omega_p(A)$ for $p\ge2$. For each pair, $\partial(su_p+tv_p)=0$ exactly when $t=-(-1)^p s$, and $\partial v_{p+1}=u_p-(-1)^p v_p$. Hence $\ker\partial_p=\Bbbk(u_p-(-1)^p v_p)=\im\partial_{p+1}$ on its alternating paths. Hence, $H_p(A)=0$. Applying the same calculation to the bidirectional pairs in $K$ gives $H_p(K)=0$.

Write $B_1(Y)=\partial\Omega_2(Y)$. The preceding computation gives
\[
 B_1(A)=\spann_{\Bbbk}\{e_{ab}+e_{ba}:a\rightleftarrows b\text{ in }A\},
 \quad B_1(A)\cap\Omega_1(K)=B_1(K).
\]
In a linear combination of the generators $e_{ab}+e_{ba}$, the edges $a\to b$ and $b\to a$ have the same coefficient. If either edge is absent from $K$, that coefficient must vanish for the chain to lie in $\Omega_1(K)$. This proves the intersection formula and hence injectivity of $H_1(K)\to H_1(A)$.
\end{proof}

The regular and nonregular boundaries can give different first homology groups even for $Y=(a\rightleftarrows b)$. Indeed, $\partial^\Lambda e_{aba}=e_{ba}-e_{aa}+e_{ab}$ and $\partial e_{aba}=e_{ba}+e_{ab}$. The class of $e_{ab}+e_{ba}$ generates $H_1^\Lambda(Y)\cong\Bbbk$, while $H_1(Y)=0$. We now prove the vanishing result from \cite[Theorem~1.1]{GrigorLiu} for the boundaries $\partial^\Lambda$ and $\partial$.

Let $G=(V,E)$ be a finite digraph with $V\subseteq\mathbb Z^2$ and only horizontal and vertical edges. On each row, list the vertices from left to right as $v_1,\ldots,v_k$. The horizontal edges on that row are exactly $v_i\to v_{i+1}$ for $1\le i<k$.

\begin{theorem}\label{thm:temporal-vanishing}
If each vertex of $G$ has at most one outgoing vertical edge, then
\[
 H_q^\Lambda(G)=H_q(G)=0,\quad q\ge2.
\]
\end{theorem}
\begin{proof}
We induct on the number $m$ of distinct first coordinates appearing in $V$. When $m\le1$, apply Lemma~\ref{lem:outdegree-one}. For $m\ge2$, let $x_0$ be the minimal first coordinate and set
\[
 A=G[V\cap(\{x_0\}\times\mathbb Z)],\quad
 B=G[V\setminus V(A)].
\]
Each row meeting $A$ and $B$ contributes one edge from its vertex in $A$ to its leftmost vertex in $B$. These edges form a matching since distinct rows have disjoint endpoints. The map $j_A:K\hookrightarrow A$ identifies the interface with a subdigraph of $A$, whose edges are all vertical. Thus, Lemma~\ref{lem:outdegree-one} applies to $K\subseteq A$.

By Lemma~\ref{lem:outdegree-one} and the prism decomposition from Theorem~\ref{thm:matching-cone}, one has
\[
 \Omega_2^\Lambda(G)=\Omega_2^\Lambda(B)\oplus h\Omega_1^\Lambda(K),\quad
 \Omega_p^\Lambda(G)=\Omega_p^\Lambda(B),\quad p\ge3.
\]
Let $u+h(c)$ be a $\partial^\Lambda$-cycle, with $u\in\Omega_2^\Lambda(B)$ and $c\in\Omega_1^\Lambda(K)$. The component of $\partial^\Lambda(u+h(c))$ in $\Omega_1^\Lambda(A)$ is $-j_Ac$. Since $j_A$ is injective on chains, we have $c=0$. Thus, $G$ and $B$ have the same $\partial^\Lambda$-cycles in degree two.
For $q\ge2$, the $(q+1)$-chain spaces agree, so the $q$-boundaries satisfy
\[
 \partial^\Lambda\Omega_{q+1}^\Lambda(B)=\partial^\Lambda\Omega_{q+1}^\Lambda(G).
\]
For $q\ge3$, the $q$-cycles also coincide because $\Omega_q^\Lambda(B)=\Omega_q^\Lambda(G)$ and $\partial^\Lambda$ on $B$ is the restriction of $\partial^\Lambda$ on $G$. Thus, the inclusion $B\hookrightarrow G$ induces $H_q^\Lambda(B)\cong H_q^\Lambda(G)$ for each $q\ge2$.

For the regular boundary, Lemma~\ref{lem:outdegree-one} yields $H_q(A)=H_q(K)=0$ for $q\ge2$. Corollary~\ref{cor:matching-les} gives $H_q(B)\cong H_q(G)$ for $q\ge3$, and the exact segment
\[
 0\longrightarrow H_2(B)\longrightarrow H_2(G)\longrightarrow H_1(K)
 \xrightarrow{(-j_{A*},\,j_{B*})}H_1(A)\oplus H_1(B).
\]
By Lemma~\ref{lem:outdegree-one}, $j_{A*}:H_1(K)\to H_1(A)$ is injective. Hence, $(-j_{A*},j_{B*})$ is injective. Exactness makes $H_2(G)\to H_1(K)$ zero, so the inclusion $B\hookrightarrow G$ induces an isomorphism on $H_2$. Removing the leftmost layer preserves the horizontal successor edges and the vertical outdegree bound. The induction hypothesis gives $H_q^\Lambda(B)=H_q(B)=0$ for $q\ge2$. The maps induced by $B\hookrightarrow G$ now give $H_q^\Lambda(G)=H_q(G)=0$ in these degrees.
\end{proof}

\subsection{Gluing and suspension}\label{sec:matching-suspension}

Let $A$ and $B$ be the two layers of $G$, with matching edges $a_x\to b_x$ indexed by the vertices of the interface $K$.

\begin{proposition}\label{prop:matching-regularity}
If $\LH_*(A,a_x)=0$, then the connecting homomorphism yields
\begin{equation*}
 \LH_q(G,a_x)\cong\LH_{q-1}(K,x),\quad q\ge0.
\end{equation*}
Similarly, if $\LH_*(B,b_x)=0$, then $\LH_q(G,b_x)\cong\LH_{q-1}(K,x)$ for $q\ge0$.
\end{proposition}
\begin{proof}
For $T=\{a_x\}$, we have $T\cap V(B)=\varnothing$, so the quotient for layer $B$ is $Q_*(B,\varnothing)=0$. The exact sequence of Theorem~\ref{thm:matching-local-cone} gives
\[
 \LH_q(A,a_x)\longrightarrow\LH_q(G,a_x)\longrightarrow
 \LH_{q-1}(K,x)\longrightarrow\LH_{q-1}(A,a_x).
\]
Since $\LH_*(A,a_x)=0$, we have the isomorphism $\LH_q(G,a_x)\to\LH_{q-1}(K,x)$. For $T=\{b_x\}$, the isomorphism $\LH_q(G,b_x)\to\LH_{q-1}(K,x)$ follows from $\LH_*(B,b_x)=0$.
\end{proof}

\begin{corollary}\label{cor:matching-regularity}
Let $n\ge1$. Suppose each unmatched vertex is regular of dimension $n$ within its layer, each vertex of $K$ is regular of dimension $n-1$, and
\[
 \LH_q(A,a_x)=\LH_q(B,b_x)=0,
 \quad x\in V(K),\ q\ge0.
\]
Then every vertex of $G$ is regular of dimension $n$.
\end{corollary}
\begin{proof}
At $a_x$ and $b_x$, Proposition~\ref{prop:matching-regularity} identifies the local homology in degree $q$ with $\LH_{q-1}(K,x)$. This group equals $\Bbbk$ when $q=n$ and vanishes otherwise. At an unmatched vertex, the local homology groups coincide with those of its layer by Corollary~\ref{cor:matching-point-local}.
\end{proof}

\begin{definition}\label{def:matching-suspension}
For a nonempty digraph $H$, the \textit{matching suspension} $\mathsf S_M H$ consists of two disjoint copies $H^-,H^+$, two new vertices $a,b$, and the edges
\[
 a\longrightarrow x^-,\quad x^-\longrightarrow x^+,\quad
 x^+\longrightarrow b,\quad x\in V(H).
\]
\end{definition}

\begin{theorem}\label{thm:matching-suspension}
For each $x\in V(H)$ and $q\ge0$, we have isomorphisms
\begin{equation*}
 \begin{aligned}
 \widetilde H_q(\mathsf S_M H)&\cong\widetilde H_{q-1}(H),\\
 \LH_q(\mathsf S_M H,a)\cong\LH_q(\mathsf S_M H,b)
 &\cong\widetilde H_{q-1}(H),\\
 \LH_q(\mathsf S_M H,x^-)\cong\LH_q(\mathsf S_M H,x^+)
 &\cong\LH_{q-1}(H,x).
 \end{aligned}
\end{equation*}
If $H$ has no directed cycles, then neither does $\mathsf S_M H$.
\end{theorem}
\begin{proof}
Put $X=\mathsf S_M H$, $A=a*H^-$, and $B=H^+*b$. Since $\widetilde H_*(A)=\widetilde H_*(B)=0$, the connecting map in Corollary~\ref{cor:matching-les} gives
\[
 H_q(X)\xrightarrow{\cong}H_{q-1}(H),\quad q\ge2.
\]
In low degrees the same sequence reads
\[
 0\longrightarrow H_1(X)\longrightarrow H_0(H)
 \xrightarrow{(-\epsilon,\epsilon)}\Bbbk^2
 \longrightarrow H_0(X)\longrightarrow0,
\]
where $\epsilon:H_0(H)\to\Bbbk$ is the augmentation. Since $H$ is nonempty, $\epsilon$ is surjective. The map $(-\epsilon,\epsilon)$ has kernel $\ker\epsilon$ and image $\Bbbk(-1,1)$. Summing the two coordinates identifies its cokernel with $\Bbbk$. Exactness now gives $H_1(X)\cong\ker\epsilon=\widetilde H_0(H)$ and $H_0(X)\cong\Bbbk$.

At the cone vertices, Corollary~\ref{cor:matching-point-local} and the exact sequences of the cone pairs give
\[
 \begin{aligned}
 \LH_q(X,a)&\cong H_q(A,H^-)\cong\widetilde H_{q-1}(H),\\
 \LH_q(X,b)&\cong H_q(B,H^+)\cong\widetilde H_{q-1}(H).
 \end{aligned}
\]
For $x\in V(H)$, the deletions $A-x^-$ and $B-x^+$ are nonempty cones. Thus, $\LH_*(A,x^-)=\LH_*(B,x^+)=0$, and Proposition~\ref{prop:matching-regularity} gives
\[
 \LH_q(X,x^-)\cong\LH_{q-1}(H,x)\cong\LH_q(X,x^+),\quad q\ge0.
\]
The vertex $a$ is a source, $b$ is a sink, and every matching edge points from $H^-$ to $H^+$. Hence, every directed cycle in $X$ lies in $H^-$ or $H^+$.
\end{proof}

Suppose $d\ge0$, each vertex of $H$ is regular of dimension $d$, and $\widetilde H_q(H)$ is $\Bbbk$ for $q=d$ and zero otherwise. By Theorem~\ref{thm:matching-suspension}, $\mathsf S_M H$ has reduced homology $\Bbbk$ in degree $d+1$ and vanishes elsewhere, and every vertex of $\mathsf S_M H$ is regular of dimension $d+1$.

\begin{example}\label{ex:matching-ten-vertex-local}
The digraph $G$ in Example~\ref{ex:matching-ten-vertex} is $\mathsf S_M Z_4$. For each $x\in V(Z_4)$, $Z_4-x$ is a directed tree. The long exact sequence of the pair $(Z_4,Z_4-x)$ gives
\[
 \LH_q(Z_4,x)\cong
 \begin{cases}\Bbbk,&q=1,\\0,&q\ne1.\end{cases}
\]
Theorem~\ref{thm:matching-suspension} then gives
\[
 \LH_q(G,v)\cong
 \begin{cases}\Bbbk,&q=2,\\0,&q\ne2,\end{cases}
 \quad v\in V(G).
\]
For the full graph, Example~\ref{ex:matching-ten-vertex} gives $H_0(G)\cong H_2(G)\cong\Bbbk$ and $H_1(G)=0$.
\end{example}

\section{Additivity and interactions}\label{sec:separation}

Fix a regular path complex $P$ on a finite vertex set, together with pairwise disjoint nonempty subsets $T_1,\ldots,T_n$. Let
\[
 T=\bigcup_{i=1}^nT_i,\quad
 d=\min_{i\ne j}d_P(T_i,T_j).
\]
We use distance in the underlying undirected graph, with distance $\infty$ between components and $d=\infty$ for $n=1$. All chain complexes are unaugmented. We use the local quotient $Q_*(P,T)$ of Subsection~\ref{sec:local-definition}.

\subsection{Additivity and the interaction complex}\label{sec:additivity-interaction}

Let $K$ be a finite simplicial complex and $T\subseteq V(K)$. Write $K[T]$ and $K-T$ for the induced subcomplexes on $T$ and its complement, respectively. Classical results show
\[
 \LH_r(P(K),T)\cong H_r(K,K-T;\Bbbk)
 \cong H_r(|K|,|K|\setminus|K[T]|;\Bbbk).
\]
For $p=\sum_{x\in V(K)}\lambda_x x\in |K|\setminus|K[T]|$, define
\[
 r(p)=\frac{\sum_{x\notin T}\lambda_x x}{\sum_{x\notin T}\lambda_x}\in |K-T|.
\]
The denominator is positive since $p\notin|K[T]|$. The vertices of $r(p)$ lie in the simplex supporting $p$, and $p$ and $r(p)$ both have positive total weight outside $T$. Thus the segment from $p$ to $r(p)$ stays in that simplex and avoids $|K[T]|$. Since $r(p)=p$ on $|K-T|$, the straight-line homotopy $H(p,t)=(1-t)p+t\,r(p)$ is a strong deformation retraction of $|K|\setminus|K[T]|$ onto $|K-T|$.

If the vertices of $T$ are pairwise nonadjacent, each simplex meeting $T$ contains a unique vertex $v\in T$. The natural projections therefore give a chain isomorphism
\[
 Q_*(P(K),T)\xrightarrow{\cong}\bigoplus_{v\in T}Q_*(P(K),v),
\]
sending the class of each such simplex to its class in the $v$-summand. If two vertices of $T$ span an edge, then $|K[T]|$ contains the entire edge, including its interior.

The inclusions $P-T\subseteq P-T_i$ induce a natural chain map
\begin{equation*}
 \rho:Q_*(P,T)\longrightarrow\bigoplus_{i=1}^nQ_*(P,T_i),
 \quad[c]\longmapsto([c],\ldots,[c]).
\end{equation*}
Write $\rho_r=H_r(\rho)$. If a path-complex morphism $f:P\to P'$ satisfies $f^{-1}(T_i')=T_i$ for each $i$, then its maps on local quotients commute with $\rho$.

\begin{theorem}\label{thm:separated-additivity}
The chain map $\rho$ is injective in each degree and bijective for $0\le p<d$. The induced map
\begin{equation*}
 \rho_r:\LH_r(P,T)\longrightarrow\bigoplus_{i=1}^n\LH_r(P,T_i)
\end{equation*}
is an isomorphism whenever $d>r+1$ and is surjective whenever $d>r$.
\end{theorem}
\begin{proof}
The intersection formula for induced restrictions gives
\[
 \bigcap_i\Omega_*(P-T_i)=\Omega_*(P-T),
\]
so $\ker\rho=0$. Let $p<d$ and choose a class in each $Q_p(P,T_i)$. By Theorem~\ref{thm:chain-locality}, the $i$-th class has a representative $c_i$ supported in $N_p(T_i)$. For $j\ne i$, this neighborhood misses $T_j$, so $c_i\in\Omega_p(P-T_j)$ by Lemma~\ref{lem:intersection}. Thus, $[\sum_i c_i]$ maps to the chosen tuple. The neighborhoods $N_p(T_i)$ may overlap.

Put $C_*=Q_*(P,T)$ and $D_*=\bigoplus_i Q_*(P,T_i)$. If $d>r$, each $z\in Z_r(D_*)$ has a lift $c\in C_r$. Since $\rho(\partial c)=\partial z=0$ and $\rho$ is injective, $c$ is a cycle. Thus, $\rho_r$ is surjective. If $d>r+1$, then $\rho(C_{r+1})=D_{r+1}$ as well, so $\rho(\partial C_{r+1})=\partial D_{r+1}$ and $\rho_r$ is injective.
\end{proof}

\begin{definition}\label{def:interaction}
The \textit{multiple-mark interaction complex} is defined by
\begin{equation*}
 \mathcal I_*(P;T_1,\ldots,T_n)
 :=\operatorname{coker}\!\left(
 Q_*(P,T)\xrightarrow{\rho}\bigoplus_iQ_*(P,T_i)\right).
\end{equation*}
\end{definition}

Write $\mathcal I_*=\mathcal I_*(P;T_1,\ldots,T_n)$. By Theorem~\ref{thm:separated-additivity}, the sequence
\[
 0\longrightarrow Q_*(P,T)\xrightarrow{\rho}
 \bigoplus_i Q_*(P,T_i)\longrightarrow\mathcal I_*\longrightarrow0
\]
is exact.

\begin{proposition}\label{prop:interaction-mv}
For two marked sets $T_1,T_2$, the following conditions are equivalent.
\begin{enumerate}[label=\textup{(\roman*)},leftmargin=2.4em]
 \item The natural map $\rho_r:\LH_r(P,T_1\cup T_2)\to\LH_r(P,T_1)\oplus\LH_r(P,T_2)$ is an isomorphism for all $r\ge0$.
 \item $H_q(\mathcal I_*(P;T_1,T_2))=0$ for all $q\ge0$.
 \item For each $r\ge0$, the inclusion $P-T_1\hookrightarrow P$ induces an isomorphism $\LH_r(P-T_1,T_2)\to\LH_r(P,T_2)$.
\end{enumerate}
Condition \textup{(iii)} may equivalently be stated with $T_1$ and $T_2$ interchanged.
\end{proposition}

\begin{proof}
Write $C=\Omega_*(P)$, $U=\Omega_*(P-T_1)$, and $V=\Omega_*(P-T_2)$. The difference map
\[
 (C/U)\oplus(C/V)\longrightarrow C/(U+V),
 \quad([a],[b])\longmapsto[b-a]
\]
is a surjective chain map. If $b-a=u+v$ with $u\in U$ and $v\in V$, then $a+u=b-v$ represents both classes. Thus, the kernel is the diagonal image of $C/(U\cap V)$. Since $U\cap V=\Omega_*(P-(T_1\cup T_2))$ by Lemma~\ref{lem:intersection}, we obtain
\begin{equation*}
 \mathcal I_*(P;T_1,T_2)
 \cong\frac{C}{U+V}
 =\frac{\Omega_*(P)}{\Omega_*(P-T_1)+\Omega_*(P-T_2)}.
\end{equation*}
This is the defect complex of Section~\ref{sec:mv} for $(P-T_1,P-T_2)$. Since $T_1\cap T_2=\varnothing$, we have $U/(U\cap V)=Q_*(P-T_1,T_2)$. The inclusion $U\hookrightarrow C$ gives the short exact sequence
\begin{equation*}
 0\longrightarrow Q_*(P-T_1,T_2)
 \longrightarrow Q_*(P,T_2)
 \longrightarrow\mathcal I_*(P;T_1,T_2)\longrightarrow0.
\end{equation*}
Taking homology of the short exact sequence following Definition~\ref{def:interaction} gives \textup{(i)}$\Longleftrightarrow$\textup{(ii)}. The homology sequence of the displayed short exact sequence gives \textup{(ii)}$\Longleftrightarrow$\textup{(iii)}. Interchanging $T_1$ and $T_2$ gives the other assertion.
\end{proof}

In a fixed degree $r$, the weaker condition $H_{r+1}(\mathcal I_*)=H_r(\mathcal I_*)=0$ suffices for both $\rho_r$ and the map $\LH_r(P-T_1,T_2)\to\LH_r(P,T_2)$ induced by $P-T_1\hookrightarrow P$ to be isomorphisms.

\subsection{Additivity criteria}\label{sec:additivity-criteria}

Two paths with the same endpoints can contribute to the same $\partial$-invariant chain even when they meet different marked sets. The following diamond gives such a pair.

\begin{example}\label{ex:interaction-diamond}
Let $G$ have edges $a\to x$, $x\to b$, $a\to y$, and $y\to b$, and set $P=P(G)$. There is no directed path from $x$ to $y$ or from $y$ to $x$, and no allowed path contains both vertices. Nevertheless,
\[
 \Omega_2(P)=\Bbbk(e_{axb}-e_{ayb}),\quad
 \LH_*(P,x)=\LH_*(P,y)=0,\quad
 \LH_1(P,\{x,y\})=\Bbbk.
\]
In $Q_*(P,x)$, the class $[e_{axb}-e_{ayb}]$ has boundary $[e_{ax}+e_{xb}]$, which spans the degree-one cycle space. The degree-one boundary is onto $Q_0(P,x)$, so $\LH_*(P,x)=0$. Symmetry gives $\LH_*(P,y)=0$. Deleting both marks leaves the isolated vertices $a,b$. Since $\widetilde H_*(G)=0$, the connecting map gives $\LH_1(P,\{x,y\})\cong\widetilde H_0(G-\{x,y\})\cong\Bbbk$.
\end{example}

For $p\ge1$, let $W_p(a,b)$ denote the $(a,b)$ endpoint component of $\Omega_p(P)$. Let $\pi_{T_i}$ be the projection operator that retains elementary path terms containing a vertex of $T_i$ and discards all remaining terms, and define
\[
 J_p(a,b)=\{i:\pi_{T_i}W_p(a,b)\ne0\}.
\]

\begin{proposition}\label{prop:endpoint-independence}
Let $r\ge0$. Assume that
\begin{equation*}
 |J_p(a,b)|\le1,\quad 1\le p\le r+1,\ a,b\in V(P).
\end{equation*}
Then $\rho$ is an isomorphism in chain degrees $0,\dots,r+1$ and induces an isomorphism on homology in degree $r$. If the condition holds for every $p\ge1$, then
\[
 \rho:Q_*(P,T)\xrightarrow{\cong}\bigoplus_{i=1}^n Q_*(P,T_i)
\]
is a chain isomorphism.
\end{proposition}
\begin{proof}
In degree zero, the chain map $\rho$ sends the basis element $[e_v]$ of $Q_0(P,T)$ (for $v\in T_i$) to $[e_v]$ in the $i$-th summand. This assignment gives a bijection between the bases. For $1\le p\le r+1$, take a representative $c\in\Omega_p(P)$ for a class in $Q_p(P,T_i)$ and decompose $c=\sum_{a,b}c_{a,b}$. Define
\[
 c_i=\sum_{J_p(a,b)=\{i\}}c_{a,b}.
\]
By Lemma~\ref{lem:endpoints}, $c_i\in\Omega_p(P)$. The hypothesis implies $\pi_{T_i}(c-c_i)=0$ and $\pi_{T_j}(c_i)=0$ for $j\ne i$. Lemma~\ref{lem:intersection} then gives $c-c_i\in\Omega_p(P-T_i)$ and $c_i\in\Omega_p(P-T_j)$. Thus, $\rho([c_i]_T)$ has $i$-th component $[c]_{T_i}$ and $j$-th component zero for each $j\ne i$. Summing these lifts establishes surjectivity. Injectivity follows from Theorem~\ref{thm:separated-additivity}, and the isomorphisms in degrees $r-1,r,r+1$ imply the homology isomorphism.
\end{proof}

To see how distance enters the endpoint criterion, let two allowed $p$-paths have the same initial and terminal vertices. Suppose they meet $x\in T_i$ and $y\in T_\ell$ at positions $j$ and $k$, respectively. The two routes from $x$ to $y$ through the common endpoints give
\[
 d_P(x,y)\le\min\{j+k,\ 2p-j-k\}\le p.
\]
Thus, $d>p$ implies $|J_p(a,b)|\le1$ for each endpoint pair $(a,b)$. For $d\le p$, the same inequality can still be checked directly on $W_p(a,b)$.

For $S\subseteq V(P)$, define $\operatorname{Pred}_P(S)$ and $\operatorname{Succ}_P(S)$ by
\[
 \operatorname{Pred}_P(S)=\{x_0:(x_0,\ldots,x_p)\in P,\ x_p\in S\},\quad
 \operatorname{Succ}_P(S)=\{x_p:(x_0,\ldots,x_p)\in P,\ x_0\in S\}.
\]
Here $p\ge0$, so both sets contain $S$.

\begin{corollary}\label{cor:structured-separation}
The natural map
\[
 \rho:Q_*(P,T)\longrightarrow\bigoplus_{i=1}^n Q_*(P,T_i)
\]
is a chain isomorphism under any of the following conditions:
\begin{enumerate}[label=\textup{(\roman*)},leftmargin=2.4em]
 \item For each $i\ne j$, we have
 \[
  \operatorname{Pred}_P(T_i)\cap\operatorname{Pred}_P(T_j)=\varnothing
  \quad\text{or}\quad
  \operatorname{Succ}_P(T_i)\cap\operatorname{Succ}_P(T_j)=\varnothing.
 \]
 \item $P=P(G)$, where $G$ is a directed tree and no edge joins distinct marked sets.
 \item $P=P(K)$ for a simplicial complex $K$, and no edge joins distinct marked sets.
\end{enumerate}
\end{corollary}
\begin{proof}
In \textup{(i)}, if paths with endpoints $a,b$ meet distinct marked sets, then $a$ is a common predecessor and $b$ a common successor. The hypothesis excludes this pair. In \textup{(ii)}, $\Omega_p=0$ for $p\ge2$, and each vertex or edge meets at most one marked set. In \textup{(iii)}, the simplices form a chain basis. A simplex meeting two marked sets would contain an edge joining them.
\end{proof}

The endpoint criterion can hold at small separation distances, but the uniform bounds of Theorem~\ref{thm:separated-additivity} are sharp. For the examples below, write $I=(0\to1)$ and $C_m=I^{\Box m}$.

Recall that $\Bbbk[m]$ denotes the chain complex equal to $\Bbbk$ in degree $m$ and zero in every other degree, with zero differential.

\begin{example}\label{ex:separation-cube}
Let $o=(0,\ldots,0)$ and $t=(1,\ldots,1)$, with $m\ge1$. Then we have $d(o,t)=m$ and
\[
 \LH_*(C_m,o)=\LH_*(C_m,t)=0,
 \quad
 \mathcal I_*(C_m;\{o\},\{t\})\cong\Bbbk[m].
\]
Consequently, we obtain
\begin{equation*}
 \LH_q(C_m,\{o,t\})\cong
 \begin{cases}\Bbbk,&q=m-1,\\0,&q\ne m-1.\end{cases}
\end{equation*}
\end{example}
\begin{proof}
Each corner-local complex is a tensor product of factors $Q_*(I,0)$ or $Q_*(I,1)$, whose homology vanishes in each degree. The K\"unneth formula therefore gives $\LH_*(C_m,o)=\LH_*(C_m,t)=0$. In degrees $p<m$, Theorem~\ref{thm:separated-additivity} identifies $Q_p(C_m,\{o,t\})$ with $Q_p(C_m,o)\oplus Q_p(C_m,t)$. In degree $m$, the full cube has a single top shuffle generator, and every length-$m$ path runs from $o$ to $t$. Thus, the joint local chain group in dimension $m$ is $\Bbbk$, while the direct sum of the two corner-local groups is $\Bbbk^2$. Under these identifications, the comparison map is the diagonal embedding. All degrees above $m$ vanish. Hence, the cokernel is $\Bbbk[m]$, and the interaction sequence gives the asserted joint local homology groups.
\end{proof}

\begin{example}\label{ex:separation-truncated-cube}
For $r\ge1$, set $G_r=B_r^{C_{r+1}}(o)$, the full cube with its all-one vertex removed. Choose a vertex $z$ on layer $r$. Then
\begin{equation*}
 \begin{gathered}
 d(o,z)=r,\quad
 \LH_r(G_r,o)=\Bbbk,\quad \LH_*(G_r,z)=0,\\
 \LH_r(G_r,\{o,z\})=0.
 \end{gathered}
\end{equation*}
Accordingly, the comparison map
\[
 \rho_r:\LH_r(G_r,\{o,z\})
 \longrightarrow\LH_r(G_r,o)\oplus\LH_r(G_r,z)
\]
is the homomorphism $0\to\Bbbk$ and fails to be surjective.
\end{example}
\begin{proof}
The local homology at the source follows from Theorem~\ref{thm:cube-balls}. The vertex $z$ is maximal, and its star digraph $S_z$ is the lower $r$-cube formed by all subsets of $z$. By Proposition~\ref{prop:star-link-digraph}, $S_z(P(G_r))=P(S_z)$ is an induced subcomplex inside $P(G_r)$. Every directed path whose endpoints lie in $V(S_z)$ remains within $V(S_z)$, so Proposition~\ref{prop:convex} yields $\LH_q(G_r,z)\cong\LH_q(S_z,z)$. These groups vanish by the local K\"unneth formula applied at a corner of the cube.

Every length-$r$ allowed path in $G_r$ begins at $o$. The top shuffle chains $c_1,\dots,c_{r+1}$ associated to the faces $x_i=0$ generate $\Omega_r(G_r)$. Therefore $Q_r(G_r,o)=Q_r(G_r,\{o,z\})=\Omega_r(G_r)$, and the projection $Q_r(G_r,\{o,z\})\to Q_r(G_r,o)$ is the identity on $\Omega_r(G_r)$. By the exterior-algebra model in Subsection~\ref{sec:cube-local-model}, the top-dimensional local cycle space at the source is one-dimensional, spanned by an alternating sum of the $c_i$ with all coefficients nonzero.

Let $c_j$ be the chain of the unique lower $r$-face through $z$. It spans $Q_r(G_r,z)$, and $\partial[c_j]_z\ne0$ since $\LH_r(G_r,z)=Q_{r+1}(G_r,z)=0$. Let $c$ be a cycle in $Q_r(G_r,\{o,z\})$. Its image at $o$ is a multiple of the alternating sum. Normalize this sum so that the coefficient of $c_j$ is $1$, and denote the scalar by $\lambda$. Since $[c_i]_z=0$ for $i\ne j$, the image of $c$ at $z$ is $\lambda[c_j]_z$. This image is a cycle, so $0=\lambda\,\partial[c_j]_z$. As $\partial[c_j]_z\ne0$, we have $\lambda=0$ and hence $c=0$.
\end{proof}

\begin{corollary}\label{cor:optimal-separation}
The bound $d\ge r+2$ for $\rho_r$ to be an isomorphism is sharp for each $r\ge0$. The bound $d\ge r+1$ for surjectivity is sharp for each $r\ge1$. Both bounds are sharp even for finite directed acyclic graphs.
\end{corollary}
\begin{proof}
In Example~\ref{ex:separation-cube}, take $m=r+1$. At distance $r+1$, the comparison is $\rho_r:\Bbbk\to0$. For $r\ge1$, Example~\ref{ex:separation-truncated-cube} gives $\rho_r:0\to\Bbbk$ at distance $r$. Both examples are finite directed acyclic graphs.
\end{proof}

\subsection{Products and higher-order interactions}\label{sec:interaction-products}

Let $P$ and $P'$ be regular path complexes on finite vertex sets. Fix $n\ge1$, pairwise disjoint nonempty sets $T_1,\ldots,T_n\subseteq V(P)$, and a nonempty set $U\subseteq V(P')$. Put $T=\bigcup_{i=1}^n T_i$. In $P\Box P'$, the marked sets are $T_i\times U$ for $1\le i\le n$, with union $T\times U$. Write $\rho:Q_*(P,T)\to\bigoplus_{i=1}^n Q_*(P,T_i)$ for the comparison map of Definition~\ref{def:interaction}.

\begin{theorem}\label{thm:interaction-product}
There is a natural chain isomorphism
\begin{equation*}
 \mathcal I_*(P\Box P';T_1\times U,\dots,T_n\times U)\cong\mathcal I_*(P;T_1,\dots,T_n)\otimes_{\Bbbk}Q_*(P',U).
\end{equation*}
For each $q\ge0$, this induces an isomorphism
\[
 H_q\bigl(\mathcal I_*(P\Box P';T_1\times U,\ldots,T_n\times U)\bigr)\cong\bigoplus_{a+b=q}
 H_a\bigl(\mathcal I_*(P;T_1,\dots,T_n)\bigr)\otimes_{\Bbbk}\LH_b(P',U).
\]
\end{theorem}
\begin{proof}
Put $D_*=Q_*(P',U)$. By Theorem~\ref{thm:local-kunneth}, the product comparison sends $[c]_T\otimes v$ to $([c]_{T_i}\otimes v)_{i=1}^n$. Thus it equals $\rho\otimes\mathrm{id}_{D_*}$. The cokernel of $\rho\otimes\mathrm{id}_{D_*}$ is $(\operatorname{coker}\rho)\otimes_{\Bbbk}D_*$ by right exactness of the tensor product. The K\"unneth formula over $\Bbbk$ gives the homology isomorphism.
\end{proof}

\begin{example}\label{ex:higher-interaction}
Let $m\ge3$, and let $F_m$ be the directed acyclic graph formed by the paths
\[
 a\to x_i\to b,\quad 1\le i\le m.
\]
Set $T=\{x_1,\dots,x_m\}$. Every nonempty proper subset $J\subsetneq\{1,\dots,m\}$ satisfies
\[
 \LH_*(F_m,\{x_i:i\in J\})=0,
\]
and the natural comparison map
\[
 Q_*(F_m,\{x_i:i\in J\})\longrightarrow\bigoplus_{i\in J}Q_*(F_m,x_i)
\]
is a chain isomorphism. For the full marked set $T$, we have
\begin{equation*}
 \LH_q(F_m,T)\cong
 \begin{cases}\Bbbk,&q=1,\\0,&q\ne1,\end{cases}
 \quad
 \mathcal I_*(F_m;\{x_1\},\dots,\{x_m\})\cong\Bbbk[2].
\end{equation*}
\end{example}
\begin{proof}
For $c=\sum_i c_i e_{ax_ib}$, the coefficient of $e_{ab}$ in $\partial c$ is $-\sum_i c_i$. All other faces are allowed, and there are no allowed paths of length three. Hence, we have
\[
 \Omega_2(F_m)=
 \left\{\sum_{i=1}^m c_i e_{ax_ib}:\sum_i c_i=0\right\},
 \quad \Omega_p(F_m)=0,\quad p\ge3.
\]
For a $1$-chain $z=\sum_i(\lambda_i e_{ax_i}+\mu_i e_{x_i b})$, we have
\[
 \partial z=-\Bigl(\sum_i\lambda_i\Bigr)e_a
 +\sum_i(\lambda_i-\mu_i)e_{x_i}
 +\Bigl(\sum_i\mu_i\Bigr)e_b.
\]
Thus, $\partial z=0$ if and only if $\lambda_i=\mu_i$ for each $i$ and $\sum_i\lambda_i=0$. In this case, $\sum_i\lambda_i e_{ax_i b}$ is the unique $2$-chain with boundary $z$. The graph is connected, so $\widetilde H_*(F_m)=0$. For $S\subsetneq T$, at least one branch remains and the calculation gives $\widetilde H_*(F_m-S)=0$. For $S=T$, only the isolated vertices $a,b$ remain. Since $\widetilde H_*(F_m)=0$, the connecting map identifies $\LH_q(F_m,S)$ with $\widetilde H_{q-1}(F_m-S)$ for every $S\subseteq T$. This group vanishes for $S\subsetneq T$ and equals $\widetilde H_{q-1}(\{a,b\})$ for $S=T$.

Identify $\Omega_2(F_m)$ with $W=\{c\in\Bbbk^m:\sum_i c_i=0\}$ by sending $\sum_i c_i e_{ax_i b}$ to $(c_i)$. For $J\subseteq\{1,\ldots,m\}$, let $S=\{x_i:i\in J\}$ and let $\pi_J:W\to\Bbbk^J$ be the coordinate projection. A chain lies in $\Omega_2(F_m-S)$ exactly when $c_i=0$ for each $i\in J$. Hence,
\[
 Q_2(F_m,S)=W/\ker\pi_J\cong\pi_J(W).
\]
For $J=\{i\}$, the coordinate $c_i$ takes every value in $\Bbbk$, since another coordinate can be chosen to equal $-c_i$. Under $Q_2(F_m,x_i)\cong\Bbbk$, the degree-two comparison for $S$ is the inclusion $\pi_J(W)\hookrightarrow\Bbbk^J$.

For $S=T$, this inclusion is $W\hookrightarrow\Bbbk^m$, whose cokernel is $\Bbbk$ via the sum of coordinates. For a proper subset $J$, choose $k\notin J$. Any $(d_i)_{i\in J}\in\Bbbk^J$ extends to a vector in $W$ by taking $c_i=d_i$ for $i\in J$, $c_k=-\sum_{i\in J}d_i$, and all remaining coordinates zero. Thus, $\pi_J(W)=\Bbbk^J$.

In degrees zero and one, each vertex or edge meeting $S$ meets exactly one $x_i\in S$, so the comparison sends each basis element to the same vertex or edge in its unique summand. All chain groups vanish above degree two. The comparison is therefore a chain isomorphism for $S\subsetneq T$, and its cokernel for $S=T$ is $\Bbbk[2]$.
\end{proof}

\begin{corollary}\label{cor:higher-interaction-degree}
Let $m\ge3$ and $r\ge1$. With $F_m$ and $T_0=\{x_1,\ldots,x_m\}$ as in Example~\ref{ex:higher-interaction}, put
\[
 J=(-1\to0\to1),\quad G=F_m\Box J^{\Box(r-1)},\quad T=T_0\times\{0\}^{r-1},
\]
where $G=F_m$ and $T=T_0$ for $r=1$. The finite directed acyclic graph $G$ satisfies $\LH_*(G,S)=0$ for each $S\subsetneq T$, and its only nonzero local homology group at $T$ is $\LH_r(G,T)\cong\Bbbk$.
\end{corollary}
\begin{proof}
Each $S\subseteq T$ has the form $S=S_0\times\{0\}^{r-1}$. By Example~\ref{ex:interval}, the only nonzero local group of $J$ at $0$ is $\LH_1(J,0)\cong\Bbbk$. For $S_0\ne\varnothing$, Theorem~\ref{thm:local-kunneth} gives
\[
 \LH_q(G,S)\cong\LH_{q-r+1}(F_m,S_0),\quad q\ge0,
\]
where groups in negative degrees are zero. Example~\ref{ex:higher-interaction} gives zero for $S_0\subsetneq T_0$ and $\Bbbk$ in degree $q=r$ for $S_0=T_0$. For $S=\varnothing$, $Q_*(G,S)=\Omega_*(G)/\Omega_*(G)=0$. The full interaction complex of $F_m$ is $\Bbbk[2]$. By Theorem~\ref{thm:interaction-product}, taking the product with each copy of $(J,\{0\})$ shifts the interaction homology up one degree, giving $\Bbbk$ in degree $r+1$ and zero elsewhere.
\end{proof}

\subsection{Interactions under matching gluing}\label{sec:interaction-gluing}

Let $G$ be a matching gluing with layers $A,B$ and interface $K$ as in Section~\ref{sec:matching-gluing}.

Take disjoint marked sets $T_1,\ldots,T_n$ among the matched endpoints in $A$. Identify them with subsets of $V(K)$ through $j_A$, and put $T=\bigcup_iT_i$. Let
\[
 \eta:\mathcal I_*(K;T_1,\ldots,T_n)
 \longrightarrow\mathcal I_*(A;T_1,\ldots,T_n)
\]
be the chain map induced by the $A$-component of the matching-cone map.

\begin{theorem}\label{thm:matching-interaction-cone}
There is a natural chain isomorphism
\begin{equation*}
 \mathcal I_*(G;T_1,\ldots,T_n)
 \cong\operatorname{Cone}(\eta).
\end{equation*}
\end{theorem}
\begin{proof}
Apply Theorem~\ref{thm:matching-local-cone} to $T$ and to each $T_i$. The $B$-term is $Q_*(B,\varnothing)=0$. If $\rho_A$ and $\rho_K$ denote the comparison maps on the layer and the interface, the degree-$p$ map on cones is $(u,c)\mapsto(\rho_{A,p}u,\rho_{K,p-1}c)$. The cokernel of this degree-$p$ map on cones is
\[
 \mathcal I_p(G;T_1,\ldots,T_n)
 \cong\mathcal I_p(A;T_1,\ldots,T_n)
 \oplus\mathcal I_{p-1}(K;T_1,\ldots,T_n),
\]
with differential $(u,c)\mapsto(\partial u+\eta c,-\partial c)$, as in $\operatorname{Cone}(\eta)$.
\end{proof}

If all matching edges point from $A$ to $B$, then $\eta=-j_{A*}$, where $j_{A*}$ is induced by the interface inclusion $j_A:K\hookrightarrow A$. The map $(u,c)\mapsto(u,-c)$ gives a chain isomorphism $\operatorname{Cone}(-j_{A*})\cong\operatorname{Cone}(j_{A*})$. For mixed orientations, $\eta$ has the signs specified in Theorem~\ref{thm:matching-mixed-directions}.

\begin{example}\label{ex:matching-interaction-cancellation}
Let $G=F_m\Box I$ with $m\ge3$ as in Example~\ref{ex:higher-interaction}. Direct the matching from $A$ to $B$ and mark the vertices $T=\{x_1,\ldots,x_m\}$ in $A$. Distinct marks are at distance two. Since $K=A=F_m$ under the matching identification, $\eta=-\mathrm{id}_{\Bbbk[2]}$. By Theorem~\ref{thm:matching-interaction-cone}, the interaction complex has a single nonzero differential, $-\mathrm{id}:\Bbbk\to\Bbbk$ from degree three to degree two. Its homology is zero.

For each $S\subseteq T$, Theorem~\ref{thm:matching-local-cone} gives $Q_*(G,S)\cong\operatorname{Cone}(-\mathrm{id}_{Q_*(F_m,S)})$. This cone is contractible, so $G-S\hookrightarrow G$ induces a homology isomorphism. For $S\ne\varnothing$, Theorem~\ref{thm:local-kunneth} also identifies $Q_*(G,S)$ with $Q_*(F_m,S)\otimes Q_*(I,0)$, where $Q_*(I,0)$ is contractible.
\end{example}

\begin{example}\label{ex:matching-obstruction-shift}
Let $m\ge2$. Take $A=F_{m+1}$ with branches $a\to x_i\to b$, $0\le i\le m$. Let $K=F_m$ retain only branches $i=1,\ldots,m$, and let $B$ be a copy of $K$. Match every vertex of $K$ from $A$ to $B$ and mark $T=\{x_1,\ldots,x_m\}\subseteq A$. The resulting graph has $2m+5$ vertices, and we have
\[
 \widetilde H_*(G)=0,\quad
 \widetilde H_q(G-T)\cong
 \begin{cases}\Bbbk,&q=1,\\0,&q\ne1.\end{cases}
\]
For each $S\subsetneq T$, the inclusion $G-S\hookrightarrow G$ induces isomorphisms $H_q(G-S)\to H_q(G)$ for each $q\ge0$.
\end{example}
\begin{proof}
For $S\subsetneq T$, the graphs $A-S$, $B$, and $K-S$ each contain at least one branch and have zero reduced homology by the branch calculation in Example~\ref{ex:higher-interaction}. On $H_0$ the interface map is $\Bbbk\to\Bbbk^2$, $c\mapsto(-c,c)$. Corollary~\ref{cor:matching-les} gives $\widetilde H_*(G-S)=0$, including $S=\varnothing$. Since $G-S$ and $G$ are connected, $G-S\hookrightarrow G$ induces an isomorphism on $H_0$. For $S=T$, the layer $A-T$ is the path $a\to x_0\to b$, and the interface consists of $a,b$. The map on $H_0$ is
\[
 \Bbbk^2\longrightarrow\Bbbk^2,\quad(\lambda,\mu)\longmapsto(-(\lambda+\mu),\lambda+\mu).
\]
Its kernel is $\Bbbk(1,-1)$. The matching sequence gives $\widetilde H_1(G-T)\cong\Bbbk$ and $\widetilde H_q(G-T)=0$ for $q\ne1$. In $K$, the equation $\sum_{i=1}^m c_i=0$ gives interaction complex $\Bbbk[2]$. In $A$, any marked coefficients $c_1,\ldots,c_m$ extend to a $2$-chain by taking $c_0=-\sum_{i=1}^m c_i$. Thus the degree-two comparison on $A$ surjects onto $\Bbbk^m$. In degrees zero and one each basis path meets at most one mark, so the comparison is also surjective there. Its cokernel is zero in every degree. Theorem~\ref{thm:matching-interaction-cone} then gives $\mathcal I_*(G;\{x\}_{x\in T})\cong\operatorname{Cone}(\Bbbk[2]\to0)\cong\Bbbk[3]$.
\end{proof}

Corollary~\ref{cor:matching-parallel-criterion} gives the deletion criterion in terms of $j_{A*}$.

\subsection{Interactions in lattice digraphs}\label{sec:interaction-cubical}

Let $T$ be a vertex set of a lattice digraph $G$. In the cubical basis of Theorem~\ref{thm:grid-cubical-chains}, a cube $\sigma$ contributes to $Q_*(G,v)$ exactly when $v\in\operatorname{Vert}(\sigma)$. Set
\[
 C_\sigma=T\cap\operatorname{Vert}(\sigma).
\]
Let $\Bbbk^{C_\sigma}$ be the space of functions $C_\sigma\to\Bbbk$, and let $\Bbbk\mathbf1_{C_\sigma}$ be its subspace of constant functions. Both spaces are zero when $C_\sigma=\varnothing$. For a codimension-one face $\tau$ of $\sigma$, denote its coefficient in $\partial_{\mathrm{cub}}\sigma$ by $[\sigma:\tau]\in\{\pm1\}$, with the orientations fixed in Subsection~\ref{sec:positive-grid}.

\begin{theorem}\label{thm:grid-interaction}
For each $p\ge0$, there exists an isomorphism
\begin{equation*}
 \mathcal I_p(G;\{v\}_{v\in T})\cong
 \bigoplus_{\sigma\in X_p(G)}
 \frac{\Bbbk^{C_\sigma}}{\Bbbk\mathbf1_{C_\sigma}}.
\end{equation*}
Under these identifications, the boundary component from $\sigma$ to a codimension-$1$ face $\tau$ is
\[
 \partial_{\tau\sigma}:
 \frac{\Bbbk^{C_\sigma}}{\Bbbk\mathbf1_{C_\sigma}}
 \longrightarrow
 \frac{\Bbbk^{C_\tau}}{\Bbbk\mathbf1_{C_\tau}},
 \quad
 [f]\longmapsto[\sigma:\tau]\,[f|_{C_\tau}].
\]
All remaining components vanish.
\end{theorem}
\begin{proof}
Write $e_{\sigma,S}=[\Phi(\sigma)]_S$. By Theorem~\ref{thm:grid-cubical-chains}, the comparison map in degree $p$ is
\[
 \rho:\bigoplus_{\substack{\sigma\in X_p(G)\\ C_\sigma\ne\varnothing}}
 \Bbbk e_{\sigma,T}
 \longrightarrow \bigoplus_{\sigma\in X_p(G)}\Bbbk^{C_\sigma},
 \quad e_{\sigma,T}\longmapsto\mathbf1_{C_\sigma},
\]
where the $v$-coordinate in the $\sigma$-summand is the coefficient of $e_{\sigma,\{v\}}$. Taking the cokernel gives the stated direct sum. Each summand has dimension $\max\{|C_\sigma|-1,0\}$.

For $v\in C_\sigma$, the relative cubical boundary is
\[
 \partial e_{\sigma,\{v\}}
 =\sum_{\substack{\tau\text{ facet of }\sigma\\ v\in C_\tau}}
 [\sigma:\tau]\,e_{\tau,\{v\}}.
\]
Thus, for a facet $\tau$ of $\sigma$, the $\tau$-component of the boundary of $f\in\Bbbk^{C_\sigma}$ is $[\sigma:\tau]\,f|_{C_\tau}$. Since $\mathbf1_{C_\sigma}|_{C_\tau}=\mathbf1_{C_\tau}$, this formula descends to the quotient by constant functions.
\end{proof}

For $C_\sigma\ne\varnothing$, the quotient by constant functions is $\widetilde H^0(C_\sigma;\Bbbk)$, with $C_\sigma$ regarded as a discrete space.

For the opposite marks $o,t$ of $C_m$ in Example~\ref{ex:separation-cube}, the full $m$-cube is the only cube containing both. Hence, the theorem gives
\[
 \mathcal I_*(C_m;\{o\},\{t\})\cong\Bbbk[m].
\]
Subsection~\ref{sec:grid-deletion} uses the cubical description of $\mathcal I_*(G;\{v\}_{v\in T})$ to characterize when $G-T\hookrightarrow G$ induces a homology isomorphism.

\section{Parallel vertex deletion}\label{sec:parallel-deletion}

Let $G$ be a finite digraph, $T\subseteq V(G)$, and $r\ge0$. Suppose distinct vertices of $T$ are at distance at least $r+3$. If $\LH_q(G,v)=0$ for each $v\in T$ and $0\le q\le r+1$, the inclusion $G-T\hookrightarrow G$ induces a homology isomorphism through degree $r$. For lattice digraphs in $\mathbb Z^n$, pairwise distance at least $n+1$ suffices in all degrees if each marked vertex has an acyclic cubical link.

\subsection{Local deletion condition}\label{sec:deletion-local-stability}

For $T\subseteq V(G)$, we have a morphism $\iota_q:H_q(G-T)\to H_q(G)$ induced by the inclusion $G-T\hookrightarrow G$. The long exact sequence of the pair $(G,G-T)$ gives
\[
 \cdots\longrightarrow\LH_{q+1}(G,T)
 \longrightarrow H_q(G-T)\xrightarrow{\iota_q}H_q(G)
 \longrightarrow\LH_q(G,T)\longrightarrow\cdots .
\]

\begin{definition}\label{def:deletion-certificate}
A nonempty set $T\subseteq V(G)$ satisfies the \textit{local deletion condition of degree $r$} if
\[
 \LH_q(G,T)=0,\quad 0\le q\le r+1.
\]
A vertex $v$ satisfies the local deletion condition of degree $r$ if $\{v\}$ does.
\end{definition}

\begin{proposition}\label{prop:deletion-certificate}
A nonempty set $T\subseteq V(G)$ satisfies the local deletion condition of degree $r$ if and only if $\iota_q$ is an isomorphism for $0\le q\le r$ and $\iota_{r+1}$ is surjective.
\end{proposition}
\begin{proof}
For $q\ge1$, the relative long exact sequence gives
\[
 0\longrightarrow\operatorname{coker}\iota_q
 \longrightarrow\LH_q(G,T)
 \longrightarrow\ker\iota_{q-1}\longrightarrow0,
\]
and $\LH_0(G,T)\cong\operatorname{coker}\iota_0$. Hence, the vanishing of $\LH_q(G,T)$ for $0\le q\le r+1$ is equivalent to the surjectivity of $\iota_q$ for $0\le q\le r+1$ and the injectivity of $\iota_q$ for $0\le q\le r$. The same argument gives the equivalence in all degrees.

For $0\le q\le r+1$, Theorem~\ref{thm:chain-locality} gives
\[
 H_q\bigl(G[N_{r+2}^G(T)],G[N_{r+2}^G(T)]-T\bigr)
 \xrightarrow{\cong}\LH_q(G,T).
\]
The largest boundary needed for these groups is $\partial_{r+2}:Q_{r+2}\to Q_{r+1}$. Thus, allowed paths of length at most $r+2$ suffice.
\end{proof}

Let $S,T\subseteq V(G)$ be disjoint, with $T\ne\varnothing$. Put $d=d_G(S,T)$, with $d=\infty$ if $S=\varnothing$. The inclusion $G-S\hookrightarrow G$ induces the chain map
\[
 j:Q_*(G-S,T)\longrightarrow Q_*(G,T).
\]

\begin{theorem}\label{thm:deletion-local-stability}
The map $j_p$ is injective in each degree and is an isomorphism for $p<d$. The induced map $H_q(j)$ is surjective for $q<d$ and is an isomorphism for $q+1<d$.
\end{theorem}
\begin{proof}
If $S=\varnothing$, then $j$ is the identity map. For $S\ne\varnothing$, set $\mathcal I_*=\mathcal I_*(G;S,T)$. The short exact sequence constructed in the proof of Proposition~\ref{prop:interaction-mv} yields
\[
 0\longrightarrow Q_*(G-S,T)\xrightarrow{j}Q_*(G,T)
 \longrightarrow\mathcal I_*\longrightarrow0.
\]
By Theorem~\ref{thm:separated-additivity}, the comparison map for $S,T$ is an isomorphism in chain degrees $p<d$, so its cokernel $\mathcal I_p$ vanishes in those degrees. Consequently, $j$ is injective in every degree, and $j_p$ is surjective for $p<d$. In the corresponding long exact homology sequence, $H_q(\mathcal I_*)=0$ for $q<d$, which implies surjectivity of $H_q(j)$. If additionally $q+1<d$, then $H_{q+1}(\mathcal I_*)=0$, which gives injectivity.
\end{proof}

\begin{corollary}\label{cor:deletion-update-radius}
Let $S\subseteq V(G)$ and $q\ge0$. For every vertex $v\in V(G)\setminus N_{q+1}^G(S)$, the inclusion $G-S\hookrightarrow G$ induces an isomorphism
\[
 \LH_q(G-S,v)\xrightarrow{\cong}\LH_q(G,v).
\]
Therefore, after deleting $S$, it is sufficient to recheck the local deletion condition of degree $r$ only at vertices in $N_{r+2}^G(S)\setminus S$. The radii $q+1$ and $r+2$ are sharp, even within the class of finite directed acyclic graphs.
\end{corollary}
\begin{proof}
Apply Theorem~\ref{thm:deletion-local-stability} with $T=\{v\}$. Taking $0\le q\le r+1$ recovers the stated region for rechecking the local deletion condition of degree $r$. All neighborhoods are measured in the original graph $G$, prior to deletion.

To verify sharpness, let $G=C_{q+1}=I^{\Box(q+1)}$, where $I=(0\to1)$, and let $o,t$ denote its source and sink. Proposition~\ref{prop:cube-vertex-local} and Theorem~\ref{thm:cube-balls} give
\[
 d_G(o,t)=q+1,\quad
 \LH_q(C_{q+1},o)=0,\quad
 \LH_q(C_{q+1}-t,o)\cong\Bbbk.
\]
Here $C_{q+1}-t=B_q^{C_{q+1}}(o)$. Hence, deletion at distance $q+1$ can alter the degree-$q$ local homology. Setting $q=r+1$, we see that $o$ satisfies the local deletion condition of degree $r$ in $C_{r+2}$ but fails it after deleting $t$, which establishes sharpness of the radius $r+2$.
\end{proof}

\subsection{Parallel deletion}\label{sec:deletion-parallel}

Theorem~\ref{thm:deletion-local-stability} allows us to preserve the local deletion condition of degree $r$ at the remaining marked vertices while deleting sufficiently separated marked vertices.

\begin{theorem}\label{thm:robust-parallel-deletion}
Let $T\subseteq V(G)$. Assume every $v\in T$ satisfies the local deletion condition of degree $r$ and
\[
 d_G(v,w)\ge r+3,\quad v,w\in T,\ v\ne w.
\]
Then for any $A\subseteq B\subseteq T$, the homomorphism induced by the inclusion $G-B\hookrightarrow G-A$
\[
 H_q(G-B)\longrightarrow H_q(G-A)
\]
is an isomorphism for $0\le q\le r$ and surjective for $q=r+1$.
\end{theorem}
\begin{proof}
For $A\subseteq T$ and $v\in T\setminus A$, we have $d_G(A,v)\ge r+3$. By Theorem~\ref{thm:deletion-local-stability}, the inclusion $G-A\hookrightarrow G$ yields isomorphisms
\[
 \LH_q(G-A,v)\xrightarrow{\cong}\LH_q(G,v)=0,
 \quad 0\le q\le r+1.
\]
Thus, $v$ satisfies the local deletion condition of degree $r$ in $G-A$. Write $B\setminus A=\{v_1,\ldots,v_k\}$ and put $A_j=A\cup\{v_1,\ldots,v_j\}$, with $A_0=A$. Proposition~\ref{prop:deletion-certificate} applies to each inclusion $G-A_j\hookrightarrow G-A_{j-1}$. The inclusions for $j=k,k-1,\ldots,1$ compose to $G-B\hookrightarrow G-A$.
\end{proof}

In particular, the marked vertices can be deleted in any order. The distance condition is equivalent to
\[
 N_{r+2}^G(v)\cap(T\setminus\{v\})=\varnothing,\quad v\in T.
\]

All distances are measured in the original graph $G$. The balls $B_{r+2}^G(v)$ are not required to be pairwise disjoint.

\begin{corollary}\label{cor:parallel-block-deletion}
Let $T_1,\ldots,T_n$ be disjoint nonempty vertex sets, each satisfying the local deletion condition of degree $r$. Suppose that
\[
 d_G(T_i,T_j)\ge r+3,\quad i\ne j.
\]
For $J\subseteq\{1,\ldots,n\}$, put $T_J=\bigcup_{i\in J}T_i$. Then, for any $J\subseteq J'\subseteq\{1,\ldots,n\}$, the map induced by the inclusion $G-T_{J'}\hookrightarrow G-T_J$
\[
 H_q(G-T_{J'})\longrightarrow H_q(G-T_J)
\]
is an isomorphism for $0\le q\le r$ and is surjective for $q=r+1$.
\end{corollary}
\begin{proof}
For $i\notin J$, each block in $T_J$ has distance at least $r+3$ from $T_i$. Thus, Theorem~\ref{thm:deletion-local-stability} gives
\[
 \LH_q(G-T_J,T_i)\xrightarrow{\cong}\LH_q(G,T_i)=0,
 \quad 0\le q\le r+1.
\]
By Proposition~\ref{prop:deletion-certificate}, deleting $T_i$ from $G-T_J$ induces isomorphisms up through degree $r$ and a surjection in degree $r+1$. Adding the indices of $J'\setminus J$ sequentially expresses the claimed map as a composition of maps with these properties.
\end{proof}

This corollary applies to unions of the designated blocks. To delete a proper subset of a single block using this reasoning, one must verify the local deletion condition of degree $r$ for that subset.

\begin{proposition}\label{prop:parallel-fixed-degree}
Let $T=\{v_1,\dots,v_n\}$, and set $d=\min_{i\ne j}d_G(v_i,v_j)$, with $d=\infty$ whenever $n\le1$. For the homomorphism induced by the inclusion $G-T\hookrightarrow G$
\[
 \iota_r:H_r(G-T)\longrightarrow H_r(G),
\]
the following statements hold:
\begin{enumerate}[label=\textup{(\roman*)},leftmargin=2.4em]
 \item If $d\ge r+3$ and $\LH_r(G,v_i)=\LH_{r+1}(G,v_i)=0$ for every $i$, then $\iota_r$ is an isomorphism.
 \item If $d\ge r+2$ and $\LH_r(G,v_i)=0$ for every $i$, then $\iota_r$ is surjective.
\end{enumerate}
\end{proposition}
\begin{proof}
For case \textup{(i)}, we have $d>q+1$ for $q=r,r+1$. Theorem~\ref{thm:separated-additivity} then identifies $\LH_q(G,T)$ with $\bigoplus_i\LH_q(G,v_i)=0$ at both degrees. For case \textup{(ii)}, $d>r+1$ gives $\LH_r(G,T)\cong\bigoplus_i\LH_r(G,v_i)=0$. Consider the long exact sequence
\[
 \LH_{r+1}(G,T)\longrightarrow H_r(G-T)
 \xrightarrow{\iota_r}H_r(G)\longrightarrow\LH_r(G,T).
\]
Since $\LH_r(G,T)=0$, exactness gives surjectivity of $\iota_r$ in \textup{(i)} and \textup{(ii)}. In \textup{(i)}, $\LH_{r+1}(G,T)=0$ also gives injectivity.
\end{proof}

The inclusions $G_t-T\hookrightarrow G_t$ also commute with the structure maps of an increasing filtration $(G_t)_t$.

\begin{corollary}\label{cor:parallel-persistence}
Let $(G_t)_t$ be a finite increasing filtration of digraphs on a fixed vertex set $V$, and let $T\subseteq V$. Suppose that for each $t$, every $v\in T$ satisfies the local deletion condition of degree $r$ in $G_t$ and
\[
 d_{G_t}(v,w)\ge r+3,\quad v,w\in T,\ v\ne w.
\]
Then the inclusions $G_t-T\hookrightarrow G_t$ induce isomorphisms of persistence modules
\[
 \{H_q(G_t-T)\}_t\xrightarrow{\cong}\{H_q(G_t)\}_t,
 \quad 0\le q\le r.
\]
Accordingly, their barcodes coincide in these degrees.
\end{corollary}
\begin{proof}
Write $i_t:G_t-T\hookrightarrow G_t$ for the inclusion. For $a\le t$, write $j_{a,t}:G_a\hookrightarrow G_t$ and $j^T_{a,t}:G_a-T\hookrightarrow G_t-T$. Then $j_{a,t}i_a=i_tj^T_{a,t}$. Applying $H_q$ gives a morphism of persistence modules, and Theorem~\ref{thm:robust-parallel-deletion} makes each $H_q(i_t)$ an isomorphism for $q\le r$.
\end{proof}

If the filtration has a final parameter $t_{\max}$, then
\[
 d_{G_t}(v,w)\ge d_{G_{t_{\max}}}(v,w).
\]
Thus, the distance condition need only be checked in $G_{t_{\max}}$, whereas the local deletion condition of degree $r$ must hold at every parameter. For varying marked sets $T_t$, the same argument applies if each inclusion $G_t-T_t\hookrightarrow G_t$ induces homology isomorphisms in the required degrees and $G_a-T_a\subseteq G_t-T_t$ whenever $a\le t$.

\subsection{Interaction criteria}\label{sec:deletion-interactions}

Let $T=\{v_1,\ldots,v_n\}$, and put
\[
 \mathcal I_*=\mathcal I_*(G;\{v_1\},\ldots,\{v_n\}).
\]
Write $\iota_q:H_q(G-T)\to H_q(G)$ for the homomorphism induced by the inclusion $G-T\hookrightarrow G$, and $b_q:H_q(G)\to\LH_q(G,T)$ for the homomorphism induced by the quotient $\Omega_*(G)\to Q_*(G,T)$.

\begin{theorem}[Interaction criterion for parallel deletion]\label{thm:parallel-interaction-criterion}
Suppose each inclusion $G-v_i\hookrightarrow G$ induces a homology isomorphism in every degree. Then the connecting homomorphism yields an isomorphism
\[
 \theta_q:H_{q+1}(\mathcal I_*)\xrightarrow{\cong}\LH_q(G,T),
 \quad q\ge0.
\]
Set $a_q=\theta_q^{-1}\circ b_q$. There exists a natural long exact sequence
\[
 \begin{aligned}
 \cdots\longrightarrow H_{r+1}(G)
 &\xrightarrow{a_{r+1}}H_{r+2}(\mathcal I_*)
 \longrightarrow H_r(G-T)\\
 &\xrightarrow{\iota_r}H_r(G)
 \xrightarrow{a_r}H_{r+1}(\mathcal I_*)\longrightarrow\cdots .
 \end{aligned}
\]
In particular, the inclusion $G-T\hookrightarrow G$ induces a homology isomorphism in every degree if and only if $H_q(\mathcal I_*)=0$ for all $q\ge0$.
\end{theorem}
\begin{proof}
The long exact sequence for the pair $(G,G-v_i)$ implies $\LH_j(G,v_i)=0$ for all $i$ and $j\ge0$. Taking homology of the short exact sequence from Definition~\ref{def:interaction} produces the connecting map $\theta_q$, whose neighboring terms are $\bigoplus_i\LH_{q+1}(G,v_i)$ and $\bigoplus_i\LH_q(G,v_i)$. Both groups vanish, so $\theta_q$ is an isomorphism. In the exact sequence of $(G,G-T)$, identify $\LH_q(G,T)$ with $H_{q+1}(\mathcal I_*)$ via $\theta_q^{-1}$. The quotient map $b_q$ then becomes $a_q=\theta_q^{-1}b_q$. As $\mathcal I_0=0$, these isomorphisms give $H_*(\mathcal I_*)=0$ exactly when $\LH_*(G,T)=0$. Exactness for $(G,G-T)$ identifies $\LH_*(G,T)=0$ with $H_q(G-T)\to H_q(G)$ being an isomorphism in each degree.
\end{proof}

For a fixed degree $r$, exactness gives
\[
 \ker\iota_r\cong\operatorname{coker}a_{r+1},\quad
 \operatorname{coker}\iota_r\cong\im a_r.
\]
Consequently,
\[
 \iota_r\text{ is an isomorphism}
 \quad\Longleftrightarrow\quad
 a_{r+1}\text{ is surjective and }a_r=0.
\]
In particular, $H_{r+2}(\mathcal I_*)=H_{r+1}(\mathcal I_*)=0$ implies that $\iota_r$ is an isomorphism. Writing $\beta_q(Y)=\dim_{\Bbbk}H_q(Y)$, we obtain
\[
 \begin{aligned}
 \beta_r(G-T)-\beta_r(G)
 &=\dim\ker\iota_r-\dim\operatorname{coker}\iota_r\\
 &=\dim H_{r+2}(\mathcal I_*)-\rank a_{r+1}-\rank a_r.
 \end{aligned}
\]

Suppose that each inclusion $G-v_i\hookrightarrow G$ induces a homology isomorphism in each degree and that $d=\min_{i\ne j}d_G(v_i,v_j)$ is finite and $d\ge2$. Theorem~\ref{thm:separated-additivity} gives $\mathcal I_p=0$ for $p<d$. Hence, $\iota_r$ is an isomorphism for $0\le r<d-2$, while $\iota_{d-2}$ is surjective and
\[
 \ker\iota_{d-2}\cong
 \operatorname{coker}\!\left(H_{d-1}(G)
 \xrightarrow{a_{d-1}}H_d(\mathcal I_*)\right).
\]
Subsection~\ref{sec:deletion-sharpness} shows that these ranges are optimal.

For a matching gluing, we use the layers $A,B$, the common interface $K$, and the embedding $j_A:K\to A$ of Section~\ref{sec:matching-gluing}. Marked endpoints in $A$ are identified with the corresponding vertices of $K$.

\begin{corollary}\label{cor:matching-parallel-criterion}
Let $T=\{v_1,\ldots,v_n\}$ consist of matched endpoints in $A$. Suppose that each inclusion $G-v_i\hookrightarrow G$ induces a homology isomorphism in each degree. Then the inclusion $G-T\hookrightarrow G$ induces a homology isomorphism in each degree if and only if the chain map
\[
 j_{A*}:\mathcal I_*(K;\{v_1\},\ldots,\{v_n\})
 \longrightarrow\mathcal I_*(A;\{v_1\},\ldots,\{v_n\})
\]
is a quasi-isomorphism.
\end{corollary}
\begin{proof}
By Theorems~\ref{thm:parallel-interaction-criterion} and~\ref{thm:matching-interaction-cone}, the inclusion $G-T\hookrightarrow G$ induces a homology isomorphism in each degree exactly when $\operatorname{Cone}(\eta)$ is acyclic. By the cone sequence this holds exactly when $\eta$ is a quasi-isomorphism. Since $\eta=j_{A*}\circ\varepsilon$ and the orientation signs give $\varepsilon^2=\mathrm{id}$, this is equivalent to $j_{A*}$ being a quasi-isomorphism.
\end{proof}

In Example~\ref{ex:matching-interaction-cancellation}, $j_{A*}$ is the identity on $\Bbbk[2]$. The cone of $j_{A*}$ is contractible, although the source and target of $j_{A*}$ each have homology $\Bbbk$ in degree two.

For $S\subseteq V(G)$, let $i_S:G-S\hookrightarrow G$, and put
\[
 \mathcal A_r(G)=\{S\subseteq V(G):
 H_q(i_S)\text{ is an isomorphism for }0\le q\le r\}.
\]

\begin{example}
Let $G$ have edges $a\to x$, $x\to b$, $a\to y$, and $y\to b$. The inclusion $G-\{a,x,y\}=G[\{b\}]\hookrightarrow G$ induces a homology isomorphism in each degree. Deleting only $\{x,y\}$ leaves two isolated vertices, and the inclusion $G-\{x,y\}\hookrightarrow G$ induces
\[
 \Bbbk^2\longrightarrow\Bbbk,\quad(\lambda,\mu)\longmapsto\lambda+\mu
\]
on $H_0$. Thus, for each $r\ge0$, we have $\{a,x,y\}\in\mathcal A_r(G)$ and $\{x,y\}\notin\mathcal A_r(G)$. In particular, $\mathcal A_r(G)$ need not be closed under taking subsets.
\end{example}

\begin{definition}\label{def:robust-deletion-complex}
The \textit{robust deletion complex} of degree $r$ is
\[
 \mathcal D_r(G)=\{T\subseteq V(G):S\in\mathcal A_r(G)
 \text{ for each }S\subseteq T\}.
\]
\end{definition}

The family $\mathcal D_r(G)$ is closed under taking subsets and contains $\varnothing$, so it forms an abstract simplicial complex. Its vertices are exactly those $v\in V(G)$ for which $H_q(G-v)\to H_q(G)$ is an isomorphism for $0\le q\le r$.

For $A\subseteq B\subseteq T$, let $i_{B,A}:G-B\hookrightarrow G-A$. If $T\in\mathcal D_r(G)$, then the factorization $i_B=i_A\circ i_{B,A}$ yields
\[
 H_q(i_{B,A})=H_q(i_A)^{-1}H_q(i_B),\quad 0\le q\le r.
\]
Both factors on the right are isomorphisms through degree $r$. Conversely, taking $A=\varnothing$ gives the defining condition for $T\in\mathcal D_r(G)$.

\begin{theorem}\label{thm:deletion-higher-obstruction}
For integers $r\ge0$ and $m\ge3$, there exist a finite digraph $G$ and an $m$-element vertex set $T$ such that every inclusion $G-S\hookrightarrow G$ with $S\subsetneq T$ induces homology isomorphisms in all degrees, while
\[
 \widetilde H_*(G)=0,\quad
 \widetilde H_q(G-T)\cong
 \begin{cases}
 \Bbbk,&q=r,\\
 0,&q\ne r.
 \end{cases}
\]
In particular, $\mathcal D_r(G)|_T=\partial\Delta(T)$, where $\Delta(T)$ denotes the simplex with vertex set $T$.
\end{theorem}
\begin{proof}
Put $G=F_m\Box J^{\Box r}$ and $T=T_0\times\{0\}^r$. By Corollary~\ref{cor:higher-interaction-degree}, we have $\LH_*(G,S)=0$ for each $S\subsetneq T$. The only nonzero local group at $T$ is $\LH_{r+1}(G,T)\cong\Bbbk$. Since $\widetilde H_*(F_m)=\widetilde H_*(J)=0$, the K\"unneth formula gives $\widetilde H_*(G)=0$. For $S\subsetneq T$, exactness for $(G,G-S)$ gives the deletion isomorphisms. For $S=T$, the connecting map identifies $\LH_{q+1}(G,T)$ with $\widetilde H_q(G-T)$ for $q\ge0$. Hence, $\widetilde H_q(G-T)\cong\Bbbk$ for $q=r$ and zero otherwise. In degree $r$, deletion of $T$ changes reduced homology from $0$ to $\Bbbk$, whereas deletion of any proper subset preserves homology in every degree. Thus every proper subset of $T$ belongs to $\mathcal D_r(G)$ and $T$ does not.
\end{proof}

For fixed $r$, this construction gives a minimal nonface of size $m$ for each $m\ge3$. When $m=3$, every pair of vertices in $T$ spans an edge of $\mathcal D_r(G)$, but $T$ does not span a simplex. Hence $\mathcal D_r(G)$ need not be flag.

By Theorem~\ref{thm:robust-parallel-deletion}, $T\in\mathcal D_r(G)$ if each $v\in T$ satisfies the local deletion condition of degree $r$ and $d_G(v,w)\ge r+3$ for all distinct $v,w\in T$.

\begin{definition}
A nonempty set $C\subseteq V(G)$ is a \textit{connected dominating set} if $G[C]$ is weakly connected and $N_1^G(C)=V(G)$.
\end{definition}

\begin{proposition}\label{prop:deletion-connected-domination}
Let $G$ be a nonempty weakly connected digraph. A set $T\subseteq V(G)$ belongs to $\mathcal D_0(G)$ if and only if $V(G)\setminus T$ is a connected dominating set.
\end{proposition}
\begin{proof}
Put $C=V(G)\setminus T$. Since $G$ is nonempty and weakly connected, $T\in\mathcal D_0(G)$ exactly when $G-S$ is nonempty and weakly connected for each $S\subseteq T$.

Suppose $T\in\mathcal D_0(G)$. Taking $S=T$ shows that $G[C]$ is nonempty and weakly connected. For each $v\in T$, taking $S=T\setminus\{v\}$ shows that $G[C\cup\{v\}]$ is weakly connected. Thus $v$ has a neighbor in $C$, and $N_1^G(C)=V(G)$.

Conversely, suppose $C$ is a connected dominating set and let $S\subseteq T$. The graph $G-S$ contains $G[C]$, and each vertex of $T\setminus S$ has a neighbor in $C$. Hence $G-S$ is nonempty and weakly connected, as required.
\end{proof}

\subsection{Sharp bounds and limitations}\label{sec:deletion-sharpness}

The separation bounds $r+3$ for isomorphisms and $r+2$ for surjectivity are sharp even when deleting any single marked vertex induces a homology isomorphism in each degree.

Fix $r\ge0$. Let $o$ and $t$ be the source and sink of $G=C_{r+2}$, so $d_G(o,t)=r+2$.

\begin{corollary}\label{cor:parallel-iso-sharp}
For each $v\in\{o,t\}$, the inclusion $G-v\hookrightarrow G$ induces a homology isomorphism in each degree. For $i:G-\{o,t\}\hookrightarrow G$, the kernel of $H_r(i)$ has dimension one.
\end{corollary}
\begin{proof}
By Example~\ref{ex:separation-cube} with $m=r+2$, we have $\LH_*(G,o)=\LH_*(G,t)=0$. The long exact sequence of $(G,G-v)$ gives the assertion for each $v\in\{o,t\}$.

For $T=\{o,t\}$, the same example gives $\LH_{r+1}(G,T)\cong\Bbbk$. Since $G$ is a product of directed intervals, the K\"unneth formula gives $H_{r+1}(G)=0$. The exact sequence
\[
 0\longrightarrow\LH_{r+1}(G,T)\xrightarrow{\partial_*}H_r(G-T)
 \xrightarrow{H_r(i)}H_r(G)
\]
therefore identifies $\ker H_r(i)$ with $\Bbbk$. Since $d_G(o,t)=r+2$, the bound $r+3$ in Theorem~\ref{thm:robust-parallel-deletion} is sharp.
\end{proof}

For the surjectivity bound, take $Z=Z_4$ from Example~\ref{ex:matching-interface-choice}. Its edges are $0\to1$, $1\to2$, $2\to3$, and $0\to3$. We have
\[
 \Omega_p(Z)=0\quad(p\ge2),\quad H_0(Z)\cong H_1(Z)\cong\Bbbk.
\]
The cycle $z=e_{01}+e_{12}+e_{23}-e_{03}$ represents a generator of $H_1(Z)$.

Set
\[
 G_r=Z^{\Box r}\Box I,\quad I=(0\to1),\quad
 u=(0,\dots,0),\quad w=(1,\dots,1),
\]
with $G_0=I$. Distances add in the Cartesian product, so $d_{G_r}(u,w)=r\,d_Z(0,1)+d_I(0,1)=r+1$.

\begin{theorem}\label{thm:parallel-surj-sharp}
For each $v\in\{u,w\}$, the inclusion $G_r-v\hookrightarrow G_r$ induces a homology isomorphism in each degree. The inclusion $G_r-\{u,w\}\hookrightarrow G_r$ induces the map
\[
 H_r(G_r-\{u,w\})=0\longrightarrow H_r(G_r)\cong\Bbbk.
\]
\end{theorem}
\begin{proof}
The endpoint complexes $Q_*(I,0)$ and $Q_*(I,1)$ are acyclic. By the local K\"unneth formula, $\LH_*(G_r,u)=\LH_*(G_r,w)=0$. The long exact sequence of $(G_r,G_r-v)$ therefore gives the claimed homology isomorphism for each $v\in\{u,w\}$.

The shuffle products of vertex and edge chains form a basis of $\Omega_*(G_r)$ with pairwise disjoint elementary-path supports. For $x\in V(G_r)$, the images $[b]_x=b+\Omega_p(G_r-x)$ of the degree-$p$ basis chains $b$ whose supports contain a path through $x$ form a basis of $Q_p(G_r,x)$. A basis chain meeting exactly one of $u,w$ contributes the identity $\Bbbk\to\Bbbk$ to the comparison map, and one meeting both contributes the diagonal map $\Bbbk\to\Bbbk^2$. The unique basis chain meeting both marked vertices is
\[
 h=e_{01}^{\times(r+1)}.
\]
Here $e_{01}^{\times(r+1)}$ is the signed shuffle of $r+1$ copies of $e_{01}$. The quotient by the image of the comparison is
\[
 \mathcal I_*(G_r;\{u\},\{w\})\cong\Bbbk[r+1],\quad
 \LH_q(G_r,\{u,w\})\cong
 \begin{cases}
 \Bbbk,&q=r,\\
 0,&q\ne r.
 \end{cases}
\]
The identification $\LH_q(G_r,\{u,w\})\cong H_{q+1}(\mathcal I_*)$ follows from Theorem~\ref{thm:parallel-interaction-criterion}.

Put $y=z^{\times r}$, with $y=1$ for $r=0$. The cycle $c=y\times e_0$ represents a generator of $H_r(G_r)\cong\Bbbk$. To compute its image in $\LH_r(G_r,\{u,w\})$, take
\[
 b=\bigl((-1)^{r+1}[y\times e_{01}]_u,\ 0\bigr)\in Q_{r+1}(G_r,u)\oplus Q_{r+1}(G_r,w).
\]
The comparison map in degree $r$ is
\[
 \rho:Q_r(G_r,\{u,w\})\longrightarrow Q_r(G_r,u)\oplus Q_r(G_r,w),
 \quad [d]_{\{u,w\}}\longmapsto([d]_u,[d]_w).
\]
Since $\partial y=0$, $[y\times e_1]_u=0$, and $[c]_w=0$, we obtain
\[
 \partial b=\bigl([y\times e_0-y\times e_1]_u,0\bigr)
 =([c]_u,[c]_w)=\rho([c]_{\{u,w\}}).
\]
The coefficients of $[h]_u$ and $[h]_w$ in $b$ are $(-1)^{r+1}$ and $0$, respectively. Since this pair lies outside the diagonal, the image $\bar b$ generates $H_{r+1}(\mathcal I_*)\cong\Bbbk$. The connecting isomorphism takes $[\bar b]$ to the relative class of $c$. Since $[c]$ also generates $H_r(G_r)$, the map
\[
 H_r(G_r)\longrightarrow\LH_r(G_r,\{u,w\})
\]
is an isomorphism. Since $\LH_{r+1}(G_r,\{u,w\})=0$, exactness identifies $H_r(G_r-\{u,w\})$ with the kernel of this isomorphism. Thus $H_r(G_r-\{u,w\})=0$.
\end{proof}

The inclusion $i:G_r-\{u,w\}\hookrightarrow G_r$ induces the map
\[
 H_r(i):H_r(G_r-\{u,w\})\longrightarrow H_r(G_r),
 \quad 0\longrightarrow\Bbbk.
\]
Thus, $\operatorname{coker}H_r(i)\cong\Bbbk$ at distance $d_{G_r}(u,w)=r+1$, so the surjectivity bound $r+2$ in Proposition~\ref{prop:parallel-fixed-degree} is optimal.

By Proposition~\ref{prop:deletion-certificate}, the local deletion condition of degree $r$ requires surjectivity in degree $r+1$ as well as isomorphisms through degree $r$.

\begin{theorem}\label{thm:no-distance-from-single-deletion}
For any integers $r,D\ge0$, there exists a finite digraph $G$ with vertices $u,w$ satisfying $d_G(u,w)\ge D$ such that both homomorphisms
\[
 H_q(G-u)\longrightarrow H_q(G),\quad
 H_q(G-w)\longrightarrow H_q(G)
\]
are isomorphisms for $0\le q\le r$, but
\[
 \beta_r(G-\{u,w\})=\beta_r(G)+1.
\]
\end{theorem}
\begin{proof}
For $N\ge4$, let $Z_N$ be the digraph with edges
\[
 0\longrightarrow1\longrightarrow\cdots\longrightarrow N-1,
 \quad 0\longrightarrow N-1.
\]
As in Example~\ref{ex:matching-interface-choice}, $\Omega_p(Z_N)=0$ for $p\ge2$ and $H_0(Z_N)\cong H_1(Z_N)\cong\Bbbk$. For each vertex $x$, the underlying graph of $Z_N-x$ is a path, so $\widetilde H_*(Z_N-x)=0$. The exact sequence of $(Z_N,Z_N-x)$ gives
\[
 \LH_q(Z_N,x)\cong
 \begin{cases}
 \Bbbk,&q=1,\\
 0,&q\ne1,
 \end{cases}
\]
and $H_1(Z_N)\xrightarrow{\cong}\LH_1(Z_N,x)$.

Set $G=Z_N\Box Z_4^{\Box r}$. By the Künneth formula, for any vertex $x$, one has
\[
 H_{r+1}(G)\cong\Bbbk,\quad
 \LH_q(G,x)\cong
 \begin{cases}
 \Bbbk,&q=r+1,\\
 0,&q\ne r+1.
 \end{cases}
\]
The map $H_{r+1}(G)\to\LH_{r+1}(G,x)$ is the tensor product of the isomorphisms on each factor, so it is itself an isomorphism. The relative exact sequence then shows that the inclusion $G-x\hookrightarrow G$ induces homology isomorphisms in degrees $0,\dots,r$.

Put $L=\max\{D,r+3\}$ and choose $N=2L$. Take $u,w$ with first coordinates $0,L$ and identical coordinates in every $Z_4$ factor. The underlying graph of $Z_N$ is a cycle of length $2L$, so
\[
 d_G(u,w)=d_{Z_N}(0,L)=L.
\]
Since $L\ge r+3$, Theorem~\ref{thm:separated-additivity} yields
\[
 \LH_{r+1}(G,\{u,w\})\cong\Bbbk^2,\quad
 \LH_r(G,\{u,w\})=0.
\]
Let $b:H_{r+1}(G)\to\LH_{r+1}(G,\{u,w\})$ be induced by the quotient. Its components at $u$ and $w$ are isomorphisms $\Bbbk\to\Bbbk$, so $b$ has rank one and $\operatorname{coker}b\cong\Bbbk$. Exactness for $(G,G-\{u,w\})$ identifies $\operatorname{coker}b$ with the kernel of $H_r(G-\{u,w\})\to H_r(G)$ and gives
\[
 0\longrightarrow\Bbbk\longrightarrow H_r(G-\{u,w\})
 \longrightarrow H_r(G)\longrightarrow0,
\]
which establishes the claimed Betti-number relation.
\end{proof}

A bounded neighborhood also fails to determine whether $G-v\hookrightarrow G$ induces a homology isomorphism.

\begin{proposition}\label{prop:deletion-not-local}
For each pair of integers $r,s\ge0$, there exist digraphs $G,G'$ and vertices $v\in V(G)$, $v'\in V(G')$ such that $(B_s^G(v),v)\cong(B_s^{G'}(v'),v')$ as rooted digraphs.
The inclusion $i_v:G-v\hookrightarrow G$ induces an isomorphism on $H_r$, whereas the inclusion $i_{v'}:G'-v'\hookrightarrow G'$ induces a surjective, noninjective map on $H_r$.
\end{proposition}
\begin{proof}
Choose $L>s+1$ and set $N=2L+2$. Let $Z_N$ be as in the proof of Theorem~\ref{thm:no-distance-from-single-deletion}, and put $P_{2L}=(0\to1\to\cdots\to2L)$. For $Y=Z_4^{\Box r}$ and $a\in V(Y)$, set
\[
 G=Z_N\Box Y,\quad G'=P_{2L}\Box Y,\quad v=v'=(L,a),
\]
with $Y$ a single vertex if $r=0$. Since
\[
 B_s^{Z_N}(L)=B_s^{P_{2L}}(L)=(L-s\to\cdots\to L+s),
\]
the product distance formula gives $B_s^G(v)=B_s^{G'}(v')$ as induced rooted digraphs.

For $G=Z_N\Box Y$, the proof of Theorem~\ref{thm:no-distance-from-single-deletion} gives
\[
 H_{r+1}(G)\xrightarrow{\cong}\LH_{r+1}(G,v),
 \quad \LH_r(G,v)=0.
\]
By exactness for $(G,G-v)$, surjectivity of $H_{r+1}(G)\to\LH_{r+1}(G,v)$ gives $\ker H_r(i_v)=0$, and $\LH_r(G,v)=0$ gives $\operatorname{coker}H_r(i_v)=0$.

For $G'=P_{2L}\Box Y$, deleting $L$ splits the interval into two intervals. The exact sequence of $(P_{2L},P_{2L}-L)$ gives $\LH_1(P_{2L},L)\cong\Bbbk$ and $\LH_q(P_{2L},L)=0$ for $q\ne1$. By the local K\"unneth formula, $\LH_q(G',v')\cong\LH_{q-1}(Y,a)$, since $\LH_*(P_{2L},L)$ is concentrated in degree one. Thus, $\LH_{r+1}(G',v')\cong\Bbbk$ and $\LH_r(G',v')=0$. The ordinary K\"unneth formula gives $H_{r+1}(G')\cong H_{r+1}(Y)=0$. The exact sequence of $(G',G'-v')$ now reads, in degree $r$,
\[
 0\longrightarrow\LH_{r+1}(G',v')\xrightarrow{\partial}
 H_r(G'-v')
 \xrightarrow{H_r(i_{v'})}H_r(G')\longrightarrow0.
\]
Hence, $H_r(i_{v'})$ is surjective with kernel $\im\partial\cong\Bbbk$.
\end{proof}

If $s\ge r+2$, finite-radius excision identifies the local homology groups of the two rooted graphs through degree $r+1$. For the inclusion $i_v:G-v\hookrightarrow G$, exactness gives
\[
 \ker H_r(i_v)\cong\operatorname{coker}\!\left(H_{r+1}(G)\longrightarrow\LH_{r+1}(G,v)\right).
\]
In Proposition~\ref{prop:deletion-not-local}, this cokernel is zero for $(G,v)$ and one-dimensional for $(G',v')$. Under the local deletion condition of degree $r$, both $\LH_{r+1}(G,v)$ and $\LH_r(G,v)$ vanish, so $H_r(i_v)$ is an isomorphism.

\subsection{Deletion in lattice digraphs}\label{sec:grid-deletion}

Let $G$ be a lattice digraph in $\mathbb Z^n$, with $n\ge1$. By Theorem~\ref{thm:grid-geometric-link}, the inclusion $G-v\hookrightarrow G$ induces a homology isomorphism in each degree if and only if the cubical link $\Lambda_v(G)$ has zero reduced homology over $\Bbbk$ in all degrees $j\ge-1$.

For $v\in V(G)$ and $S\subseteq V(G)$, put
\[
 \begin{aligned}
 W_v&=V(G)\cap\bigl(v+\{-1,0,1\}^n\bigr),\\
 U_G(S)&=
 \left(\bigcup_{\substack{\sigma\in X(G)\\
 \operatorname{Vert}(\sigma)\cap S\ne\varnothing}}
 \operatorname{Vert}(\sigma)\right)\setminus S.
 \end{aligned}
\]
Thus, $|W_v|\le3^n$, and $U_G(S)$ consists of the vertices outside $S$ that share a cube with a vertex of $S$. Write $N_n^G(S)$ for the set of vertices within undirected distance $n$ of $S$ in $G$. For $W\subseteq V(G)$, let $G[W]$ denote the induced subdigraph, so $B_n^G(v)=G[N_n^G(\{v\})]$.

Let $T\subseteq V(G)$. Any two vertices of a cube in $X(G)$ are at distance at most $n$. Hence, if distinct vertices of $T$ are at distance at least $n+1$, each cube contains at most one vertex of $T$.

\begin{theorem}\label{thm:grid-all-degree-deletion}
\begin{enumerate}[label=\textup{(\roman*)},leftmargin=2.4em]
 \item For each $v\in V(G)$, the inclusions $B_n^G(v)\hookrightarrow G$ and $G[W_v]\hookrightarrow G$ induce the chain isomorphisms
 \[
 Q_*(B_n^G(v),v)\xrightarrow{\cong}Q_*(G,v)
 \xleftarrow{\cong}Q_*(G[W_v],v).
 \]
 \item If $|\operatorname{Vert}(\sigma)\cap T|\le1$ for each $\sigma\in X(G)$, the comparison map is a chain isomorphism
 \[
 \rho:Q_*(G,T)\xrightarrow{\cong}\bigoplus_{v\in T}Q_*(G,v),
 \quad [c]_T\longmapsto([c]_v)_{v\in T}.
 \]
 If also $\widetilde H_*(\Lambda_v(G);\Bbbk)=0$ for each $v\in T$, the inclusion $G-S'\hookrightarrow G-S$ induces a homology isomorphism for each $S\subseteq S'\subseteq T$.
 \item For each $S\subseteq V(G)$, we have
 \[
 \begin{aligned}
 \Lambda_v(G-S)&=\Lambda_v(G),\quad v\notin S\cup U_G(S),\\
 U_G(S)&\subseteq\bigl(N_n^G(S)\setminus S\bigr)\cap\bigcup_{u\in S}W_u.
 \end{aligned}
 \]
\end{enumerate}
\end{theorem}
\begin{proof}
By Theorem~\ref{thm:grid-cubical-chains}, the classes $[\Phi(\sigma)]_T$ of $p$-cubes intersecting $T$ form a basis for $Q_p(G,T)$.

For \textup{(i)}, let $\sigma$ be a $p$-cube containing $v$. Note that every vertex $w$ of $\sigma$ satisfies
\[
 d_G(v,w)\le p\le n,\quad w-v\in\{-1,0,1\}^n.
\]
Every $p$-cube of $G$ containing $v$ therefore lies in $B_n^G(v)$ and in $G[W_v]$. For $H=B_n^G(v)$ or $H=G[W_v]$, the inclusion $H\hookrightarrow G$ induces a map $Q_p(H,v)\to Q_p(G,v)$ that sends the class of each such cube to the class of the same cube. Since these cubes index a basis on each side, the map is an isomorphism.

For \textup{(ii)}, let $S\subseteq T$. A $p$-cube cannot contain two vertices of $S$, since their distance within the cube would be at most $p\le n$, while distinct marks have distance at least $n+1$. Thus each basis cube meeting $S$ contains exactly one $v\in S$. The comparison map $\rho_S$ sends $[\Phi(\sigma)]_S$ to $[\Phi(\sigma)]_v$ in the $v$-summand. On the $v$-summand, $\rho_S^{-1}$ sends $[\Phi(\sigma)]_v$ to $[\Phi(\sigma)]_S$. If each marked link has zero reduced homology, the cubical link formula gives
\[
 \LH_q(G,S)\cong\bigoplus_{v\in S}
 \widetilde H_{q-1}(\Lambda_v(G);\Bbbk)=0,
 \quad q\ge0.
\]
The exact sequence of $(G,G-S)$ gives an isomorphism $H_q(i_S)$ for $i_S:G-S\hookrightarrow G$. For $S\subseteq S'\subseteq T$, the inclusion $G-S'\hookrightarrow G-S$ induces $H_q(i_S)^{-1}H_q(i_{S'})$.

For \textup{(iii)}, compatibility with vertex deletion yields
\[
 X(G-S)=\{\sigma\in X(G):\operatorname{Vert}(\sigma)\cap S=\varnothing\}.
\]
A simplex survives in $\Lambda_v(G-S)$ exactly when its cube avoids $S$. For $v\notin S\cup U_G(S)$, all cubes through $v$ survive, so $\Lambda_v(G-S)=\Lambda_v(G)$. For $v\in U_G(S)$, some cube contains both $v$ and a vertex $u\in S$. For these two vertices, we have $d_G(u,v)\le n$ and $v-u\in\{-1,0,1\}^n$.
\end{proof}

The ball $B_n^G(v)$ contains all cubes defining $\Lambda_v(G)$. After deleting $S$, only the links at vertices in $U_G(S)$ need to be checked again.

\begin{theorem}\label{thm:grid-optimal-bounds}
For lattice digraphs embedded in $\mathbb Z^n$, the following uniform bounds are sharp:
\begin{enumerate}[label=\textup{(\roman*)},leftmargin=2.4em]
 \item Radius $n$ for computing all local homology groups at a vertex.
 \item Radius $n$ for updating the condition $\widetilde H_*(\Lambda_v(G);\Bbbk)=0$ after vertex deletion.
 \item Separation $n+1$ for parallel deletion of vertices $v$ with $\widetilde H_*(\Lambda_v(G);\Bbbk)=0$.
\end{enumerate}
\end{theorem}

\begin{proof}
The upper bounds follow from Theorem~\ref{thm:grid-all-degree-deletion}. Let $o$ and $t$ be the source and sink of $C_n=I^{\Box n}$. Each vertex other than $t$ can be reached from $o$ in at most $n-1$ steps without passing through $t$. Hence, we have
\[
 B_{n-1}^{C_n}(o)=B_{n-1}^{C_n-t}(o)=C_n-t.
\]

Write $\Delta=\Lambda_o(C_n)$, an $(n-1)$-simplex. The full cube is the only face containing both $o$ and $t$, so $\Lambda_o(C_n-t)=\partial\Delta$. The cubical link formula in Theorem~\ref{thm:grid-geometric-link} gives
\[
 \LH_*(C_n,o)=0,\quad
 \LH_{n-1}(C_n-t,o)\cong\widetilde H_{n-2}(\partial\Delta;\Bbbk)\cong\Bbbk.
\]
For $n=1$, use $\widetilde H_{-1}(\varnothing;\Bbbk)=\Bbbk$. The rooted $(n-1)$-balls at $o$ coincide, and their ambient graphs have different local homology in degree $n-1$. This proves \textup{(i)}. Deleting $t$ changes $\widetilde H_{n-2}(\Lambda_o)$ at distance $n$, proving \textup{(ii)}.

For \textup{(iii)}, both $o$ and $t$ have simplex links. Let $i:C_n-\{o,t\}\hookrightarrow C_n$ be the inclusion. For $n\ge2$, Corollary~\ref{cor:parallel-iso-sharp} with $r=n-2$ gives
\[
 \dim_{\Bbbk}\ker H_{n-2}(i)=1.
\]
For $n=1$, the map $H_0(i)$ is $0\to\Bbbk$. Thus, separation $n$ is insufficient in each dimension $n\ge1$.
\end{proof}

\begin{remark}\label{rem:grid-improvement}
By Corollary~\ref{cor:grid-high-vanishing}, we have $H_q(G-T)=H_q(G)=0$ for $q\ge n$. It suffices to check the deletion maps in degrees $0,\ldots,n-1$. The general bounds with $r=n-1$ give the first row of the table, and the cubical model gives the optimal bounds in the second.
\begin{center}
\begin{tabular}{lccc}
\hline
Method & Computation radius & Update radius & Separation\\
\hline
General deletion criteria & $n+1$ & $n+1$ & $n+2$\\
Cubical chain model & $n$ & $n$ & $n+1$\\
\hline
\end{tabular}
\end{center}
For $n=1$, each cube of positive dimension is an edge, so $|\operatorname{Vert}(\sigma)\cap T|\le1$ for each cube is equivalent to $d_G(v,w)\ge2$ for distinct $v,w\in T$. For $n\ge2$, the cube condition is strictly weaker than pairwise separation by $n+1$. For example, take the directed path $0\to e_1\to2e_1$ in $\mathbb Z^n$ and mark its endpoints. They have distance two, and each cube contains at most one of them.
\end{remark}

If $\widetilde H_*(\Lambda_v(G);\Bbbk)=0$ for any $v\in T$, then Theorem~\ref{thm:parallel-interaction-criterion} gives a criterion for deleting $T$ without a separation assumption: the inclusion $G-T\hookrightarrow G$ induces a homology isomorphism if and only if the interaction complex $\mathcal I_*(G;\{v\}_{v\in T})$ is acyclic.

\section{Further questions}\label{sec:future}

In this section, we list several further questions remaining from this work.

\paragraph{Realizing the completed link.}
Let $P$ be a path complex, and let $v$ be a vertex of $P$. When the complex $D_*(P,v)$ defined in \ref{def:defect} is acyclic, Proposition~\ref{prop:completed-link-comparison} says that the link $L_v(P)$ is the path complex realization of the completed link up to quasi-isomorphism.

\begin{question}
Is there a path complex realization of the augmented complex of $\widehat{\mathcal L}_*(P,v)$? Ideally, such a realization should be natural with respect to inclusions preserving $v$.
\end{question}

\paragraph{The stabilization radius.}
Let $P_s=P|_{N_s(T)}$. Recall that the inclusion $P_s\hookrightarrow P$ of path complexes induces a local homology isomorphism in degree $r$ for $r\le s-1$. When $s=r$, the induced map is always surjective, yet it may admit a nontrivial kernel.
\begin{question}
Can one characterize, using allowed paths and their faces, when the map $\LH_r(P_r,T)\to\LH_r(P,T)$ is injective?
\end{question}

\paragraph{Gluing beyond matchings.}
Theorem~\ref{thm:matching-local-cone} gives the mapping-cone formula when the edges joining two induced layers form a matching.

\begin{question}
How can one extend this chain model for $Q_*(G,T)$ to interlayer edges with shared endpoints?
\end{question}

\paragraph{Vertex deletion in filtrations.}
Let $(G_t)_t$ be a finite increasing filtration on a fixed vertex set with a nonempty marked set $T$. Suppose that for each $v\in T$, the inclusion $G_t-v\hookrightarrow G_t$ induces a homology isomorphism in each degree. By Theorem~\ref{thm:parallel-interaction-criterion}, the inclusions $G_t-T\hookrightarrow G_t$ then induce isomorphisms of persistence modules precisely when the interaction complex is acyclic at each parameter.

\begin{question}
Under what conditions can changes near $T$ identify a proper subset of parameters at which testing acyclicity suffices for the entire filtration?
\end{question}

\paragraph{Localization at an allowed path.}
For an allowed path $\alpha$, let $P_{\not\supset\alpha}$ consist of the allowed paths that do not contain $\alpha$ as an ordered subsequence. If $\alpha$ represents a $d$-simplex $\sigma$ of a simplicial complex $K$, then
\[
 H_r(P(K),P(K)_{\not\supset\alpha})\cong\widetilde H_{r-d-1}(\Lk_K\sigma;\Bbbk).
\]

\begin{question}
When does an allowed $d$-path $\alpha$ in a regular path complex $P$ admit a link path complex $L_\alpha(P)$ giving the same link formula and agreeing with $L_v(P)$ when $\alpha=(v)$?
\end{question}

\section*{Acknowledgments}

This work was supported in part by the National Natural Science Foundation of China (NSFC Grant No. 12401080), the Science and Technology Research Program of Chongqing Municipal Education Commission (Grant No. KJQN202601151), the Scientific Research Foundation of Chongqing University of Technology, the National Key Research and Development Program of China (NKPs Grant No. 2024YFA1013201), the Natural Science Foundation of Chongqing (NSFCQ Grant No. CSTB2024NSCQ-LZX0040), and the Special Project of Chongqing Municipal Science and Technology Bureau (Grant No. 2025CCZ015).

\end{document}